\documentclass[10pt]{article}
\usepackage{amsfonts,latexsym,amstext}
\usepackage{amsmath}
\usepackage{amssymb}
\usepackage[english]{babel}
\usepackage[latin1]{inputenc}

\usepackage[usenames]{color}
\definecolor{red}{rgb}{1.0,0.0,0.0}

\definecolor{blu}{rgb}{0.0,0.0,1.0}

\definecolor{gre}{rgb}{0.03,0.50,0.03}

\newtheorem{theorem}{Theorem}[section]
\newtheorem{lemma}[theorem]{Lemma}
\newtheorem{proposition}[theorem]{Proposition}
\newtheorem{definition}[theorem]{Definition}

\newtheorem{hypothesis}[theorem]{Hypothesis}

\newtheorem{remark}[theorem]{Remark}
\newtheorem{corollary}[theorem]{Corollary}

\numberwithin{equation}{section}

\newcommand{\myref}[1]{(\ref {#1})}

\def\qed{{\hfill\hbox{\enspace${ \square}$}} \smallskip}
\def\sqr#1#2{{\vcenter{\vbox{\hrule height .#2pt \hbox{\vrule
 width .#2pt height#1pt \kern#1pt \vrule
width .#2pt} \hrule height .#2pt}}}}
\def\square{\mathchoice\sqr54\sqr54\sqr{4.1}3\sqr{3.5}3}

\def\qedo{\hbox{\hskip 6pt\vrule width6pt height7pt
depth1pt  \hskip1pt}\bigskip}

\def\eps{\varepsilon}

\def\ds{\begin{displaystyle}}
\def\eds{\end{displaystyle}}
\def\dis{\displaystyle }
\def\<{\left\langle }
\def\>{\right\rangle }

\def\dim{\noindent \hbox{{\bf Proof.} }}

\def\R{\mathbb R}
\def\N{\mathbb N}

\def\E{\mathbb E}
\def\P{\mathbb P}

\def\calc{{\cal C}}
\def\cald{{\cal D}}

\def\calf{{\cal F}}

\def\calh{{\cal H}}
\def\calk{{\cal K}}
\def\call{{\cal L}}
\def\caln{{\cal N}}
\def\calp{{\cal P}}

\def\calu{{\cal U}}

\def\call{{\cal L}}

\def\to{\rightarrow}

\begin{document}

\title{Stochastic Optimal Control with Delay in the Control:\\
solution through partial smoothing}
\date{}
 \author{Fausto Gozzi\\
Dipartimento di Economia e Finanza\\
Universit\`a LUISS - Guido Carli\\
Viale Romania 32,
00197 Roma,
Italy\\
e-mail: fgozzi@luiss.it\\
\\
Federica Masiero\\
Dipartimento di Matematica e Applicazioni\\ Universit\`a di Milano Bicocca\\
via Cozzi 55, 20125 Milano, Italy\\
e-mail: federica.masiero@unimib.it}

\maketitle
\begin{abstract}
Stochastic optimal control problems governed by delay equations with delay in the control are usually more
difficult to study than the the ones when the delay appears only in the state.
This is particularly true when we look at the associated Hamilton-Jacobi-Bellman (HJB) equation.
Indeed, even in the simplified setting
(introduced first by Vinter and Kwong \cite{VK} for the deterministic case)
the HJB equation is an infinite dimensional second order semilinear Partial Differential Equation (PDE)
that does not satisfy the so-called ``structure condition'' which substantially means
that ``the noise enters the system with the control''. The absence of such condition, together
with the lack of smoothing properties which is a common feature of problems with delay,
prevents the use of the known techniques (based on Backward Stochastic Differential Equations (BSDEs)
or on the smoothing properties of the linear part) to prove the existence of regular solutions of this HJB
equation and so no results on this direction have been proved till now.

In this paper we provide a result on existence of regular solutions of such
kind of HJB equations and we use it to solve completely the
corresponding control problem finding optimal feedback controls
also in the more difficult case of pointwise delay. The main tool used is a partial smoothing property that
we prove for the transition semigroup associated to the uncontrolled problem.
Such results holds for a specific class of equations and data which arises naturally in many applied problems.
\end{abstract}

\tableofcontents

\section{Introduction}


Optimal control problems governed by delay equations with delay in the control are usually
harder to study than the ones when the delay appears only in the state (see e.g.
\cite[Chapter 4]{BDDM07} and \cite{GM,GMSJOTA}). This is true already in the deterministic case
but things get worse in the stochastic case.
When one tries to apply the dynamic programming method the main difficulty is the fact that, even in the simplified setting introduced first by Vinter and Kwong \cite{VK} in the deterministic case (see e.g. \cite{GM} for the stochastic case), the associated HJB equation is an infinite dimensional second order semilinear PDE that
does not satisfy the so-called ``structure condition'', which substantially
means that the noise affects the system  only ``through the control''.

The absence of such condition, together with the lack of smoothing properties which is a common feature of problems with delay,
prevents the use of the known techniques, based on BSDE's
(see e.g. \cite{FT2})
or on fixed point theorems in spaces of continuous functions
(see e.g. \cite{CDP1,CDP2,DP3,G1,G2})
or in Gauss-Sobolev spaces
(see e.g. \cite{ChowMenaldi,GGSPA}),
to prove the existence of regular solutions of this HJB equation:
hence no results in this direction have been proved till now.
The viscosity solution technique can still be used (see e.g. \cite{GMSJOTA})
but to prove existence (and possibly uniqueness) of solutions that are merely continuous.
This is an important drawback in this context, since, to prove the existence of optimal feedback control strategies
through the dynamic programming approach, one needs at least the differentiability of the solution in the ``space-like'' variable.

The main aim of this paper is twofold: first to provide a new result of existence of regular solutions of such HJB equations
that holds when the state equation depends linearly on the history of the control and when the cost functional does not depend
on such history;
second to exploit such result to solve the corresponding stochastic optimal control problem finding optimal feedback
control strategies.
This allows to treat satisfactorily a specific class of state equations and data
which arises naturally in many applied problems (see e.g. \cite{FabbriGozzi08,FedTacSicon,GM,GMSJOTA,Kol-Sha,PhamBruder}). 
Note that, differently to what often happens in the literature (see e.g.
\cite{FedTacSicon,GM,GMSJOTA}) here we are able to treat also the case of ``pointwise delay''
which gives rise to an unbounded control operator in the state equation.

The key tool to prove such results is the proof of a ``partial'' smoothing property
for the transition semigroup associated to the uncontrolled equation which
we think is interesting in itself and is presented in Section 3.

We believe that such tool may allow to treat also examples where the state equation depends
on the history of the state variable, too. To keep things simpler, here we choose to develop and present the result
when this does not happen leaving the extension to a subsequent paper.


\subsection{Our results in a simple motivating case}
\label{SS:motiv}

To be more clear we now briefly describe our setting and our main result in a special case. Let $(\Omega, \calf,  \P)$
be a complete probability space
and consider the following linear controlled Stochastic
Differential Equation (SDE) in $\R$:
\begin{equation}
\left\{
\begin{array}
[c]{l}
dy(s)  =a_0 y(s) dt+b_0 u(s) ds +\int_{-d}^0b_1(\xi)u(s+\xi)d\xi+\sigma dW_s
,\text{ \ \ \ }s\in[t,T] \\
y(t)  =y_0,\\
u(\xi)=u_0(\xi), \quad \xi \in [-d,0).
\end{array}
\right.  \label{eq-contr-ritINTRO}
\end{equation}
Here $W$ is a standard Brownian motion in $\R$, and $(\calf_t)_{t\geq 0}$ is the
augmented filtration generated by $W$.
We assume $a_0,b_0\in \R$, $\sigma >0$. The parameter $d>0$ represents the maximum delay
the control takes to affect the system while $b_1$ is the density function
taking account the aftereffect of the control on the system.
The easier case is when $b_1\in L^2([-d,0],\R)$ (``distributed delay'')
while a more difficult, yet common, case is when $b_1$ is a measure, e.g. a Dirac delta in $-d$
(``pointwise delay''): both will be treated here.

The initial data are the initial state $y_0$ and the past history $u_0$ of the control.
The control $u$ belongs to $L^2_{\calf}(\Omega\times [0,T], U)$,
the space of predictable square integrable processes with values in
$U\subseteq \R$, closed.

Such kind of equations is used e.g. to model the effect of advertising on
the sales of a product (see e.g. \cite{GM,GMSJOTA}), the effect of investments
on growth (see e.g. \cite{FabbriGozzi08} in a deterministic case),
or, in a more general setting, to model optimal portfolio problems with execution delay,
(see e.g. \cite{PhamBruder}) or to model the interaction of drugs with tumor cells
(see e.g. \cite{Kol-Sha} p.17 in the deterministic case).

In many applied cases (like the ones quoted above) the goal of the problem is to minimize the following objective functional
\begin{equation}\label{costoconcretoINTRO}
J(t,x_0,u_0;u(\cdot))=\E \int_t^T \left(\bar\ell_0(s,y(s))+\bar\ell_1(u(s))\right)\;ds +\E  \bar\phi(y(T)).
\end{equation}
where $\bar\ell_0:[0,T]\times\R\rightarrow \R$, $\bar\ell_1:U\rightarrow \R$ and
$\bar\phi:\R\rightarrow \R$ are continuous functions satisfying suitable assumptions that will be introduced in Section \ref{section-statement}.
It is important to note that here $\bar\ell_0$, $\bar\ell_1$ and $\bar\phi$ do not depend on the past of the state and/or control.
This is a very common feature of such applied problems.

A standard way\footnote{It must be noted that, under suitable restrictions on the data, one can treat optimal control
problems with delay in the control also by a direct approach without transforming them in infinite dimensional problems.
However in the stochastic case such direct approach seems limited to a very special class of cases
(see e.g. \cite{LarssenRisebro03}) which does not include our model and models commonly used in applications
like the ones just quoted.} to approach these delayed control problems, introduced in \cite{VK} for the deterministic case
(see \cite{GM} for the stochastic case) is to reformulate the above linear delay equation as a linear SDE in the Hilbert space
$\calh:=\R \times L^2([-d,0],\R)$, with state variable $Y=(Y_0,Y_1)$ as follows.
\begin{equation}
\left\{
\begin{array}
[c]{l}
dY(s)  =AY(s) ds+Bu(s) ds+GdW_s
,\text{ \ \ \ }s\in[t,T] \\
Y(t)=x=(x_0,x_1),
\end{array}
\right.   \label{eq-astrINTRO}
\end{equation}
where $A$ generates a $C_0$-semigroup (see next section for precise definitions)
while, at least formally,\footnote{Note that, when $b_1$ is a measure, the above operator
$B$ is not bounded (see next section for the precise setting).}
\begin{equation}
 \label{BGINTRO}
B:\R\rightarrow \calh,\qquad Bu=(b_0 u, b_1(\cdot)u), \; u\in\R,
\qquad
G:\R\rightarrow \calh,\qquad Gx=(\sigma x, 0), \; x\in\R.
\end{equation}
Moreover $x_0=y_0$ while $x_1(\xi)=\int_{-d}^\xi b_1(\varsigma)u_0(\varsigma-\xi)d\varsigma$, ($\xi \in [-d,0]$) i.e.
the infinite dimensional datum $x_1$ depends on the past of the control.

The value function is defined as
$$
V(t,x):=\inf_{u(\cdot)\in L^2_{\calf}(\Omega\times [0,T], U)}
\E \left(\int_t^T \left[\bar\ell(s,Y_0(s))+\bar\ell_1(u(s))\right] ds +\bar\phi(Y_0(T))
\right)
$$
The associated HJB equation (whose candidate solution is the value function) is
\begin{equation}\label{HJBINTRO}
\left\{\begin{array}{l}\dis
-\frac{\partial v(t,x)}{\partial t}=
\frac{1}{2}Tr \;GG^*\nabla^2v(t,x)
+ \< Ax,\nabla v(t,x)\>_\calh
+ \bar H_{min} (\nabla v(t,x)) +\bar\ell_0(t,x),\qquad t\in [0,T],\,
x\in \calh,\\
\\
\dis v(T,x)=\bar\phi(x_0),
\end{array}\right.
\end{equation}
where, defining the current value Hamiltonian $\bar H_{CV}$ as
$$
\bar H_{CV}(p;u):=\<p,Bu\>_\calh+\bar\ell_1(u)=\<B^*p,u\>_\R +\bar\ell_1(u)
$$
we have
\begin{equation}\label{HminINTRO}
\bar H_{min} (p):=
\inf_{u\in U} \bar H_{CV}(p\,;u).
\end{equation}
The ultimate goal in studying this HJB equation is to find a solution $v$ with enough
regularity to prove a verification theorem
(i.e. that $v=V$ plus a sufficient condition for optimality) and to find
an optimal feedback map.

It is well known, see e.g. \cite{YongZhou99} Section 5.5.1, that, if the value function
$V$ is smooth enough, and if the current value Hamiltonian $\bar H_{CV}$ always admits at least a minimum point, a natural candidate
optimal feedback map is given by $(t,x) \mapsto u^*(t,x)$ where
$u^*(t,x)$ satisfies
\[
 \<\nabla V(t,x),Bu^*(t,x)\>_{\R}+\bar\ell_1(u^*(t,x))= \bar H_{min}(\nabla V(t,x)).
\]
i.e. where $u^*(t,x)$ is a minimum point of the function
$u\mapsto \bar H_{CV}(\nabla v(t,x);u)$, $\R\to \R$.
To take account of the possible unboundedness of $B$ it is convenient to define,
for $z \in \R$,
\begin{equation}\label{HminINTRObis}
H_{min} (z):=
\inf_{u\in U}  \{\<z,u\>_\R+\bar\ell_1(u) \}=:\inf_{u\in U} H_{CV}(z\,;u) \}
\end{equation}
so that
\begin{equation}\label{HminINTROter}
\bar H_{min} (p)=
\inf_{u\in U} \bar H_{CV}(p\,;u).
=\inf_{u\in U} \{\<B^*p,u\>_\R+\bar\ell_1(u ) \}
=:H_{min}(B^*p).
\end{equation}
Since $u^*(t,x)$ is a minimum point of the function
$u\mapsto H_{CV}(B^*\nabla v(t,x);u)$, $\R\to \R$,
the minimal regularity required to give sense to such term is the existence
of $B^*\nabla v(t,x)$ which we will call $\nabla^B v(t,x)$ according with the definition
and notation used e.g. in \cite{FTGgrad,Mas} in the case of bounded $B$ and generalized here in
Subsection \ref{subsection-C-directionalderivatives}
to the case of possibly unbounded $B$.
From now on we will use $H_{CV}$ and $H_{min}$ in place of $\bar H_{CV}$ and $\bar H_{min}$
writing $H_{min}(\nabla^B v(t,x))$ in place of $\bar H_{min}(\nabla v(t,x))$ in \myref{HJBINTRO}.

\bigskip

As announced, the main results of this paper are exactly to prove the existence of a
solution of the HJB equation \myref{HJBINTRO} with such regularity (Section 5) and
then to prove a verification theorem and the existence of optimal feedback controls (Section 6).
Such results are difficult to get since the known techniques do not apply to this problem.
We explain why.

To prove the existence and uniqueness of regular solutions of
HJB equations\footnote{Sometimes such equations are called semilinear Kolmogorov equations, see e.g. \cite{FT2}}
like (\ref{HJBINTRO}) in the simplest case when $H$ is Lipschitz continuous there are two main approaches in the literature:
a BSDEs based approach and a fixed point approach.

The BSDEs approach (see e.g. the paper \cite{FT2}, which is the infinite dimensional extension of results in \cite{PaPe}) is based
on representing the solution $v$ using a suitable BSDE. The regularity of $v$ is inherited from the data so all the coefficients, included the Hamiltonian $H$,
are assumed to be at least Gateaux differentiable.
To apply such approach the controlled state equation must satisfy the ``structure condition''
which, in our case, would mean that ${\rm Im} B\subseteq {\rm Im} G$.
Our HJB equation \myref{HJBINTRO} is out of such literature since the structure condition does not hold here. This is an intrinsic feature of such problem since the fact that $\operatorname{Im} B $ is not contained in $\operatorname{Im} G$ depends on the presence of the delay in the control. If the delay in the control disappears, then the structure condition hold, even if delay in the state is present (see e.g. \cite{GM,GMSJOTA,FT2,MasBanach}).

The fixed point approach consists in rewriting \myref{HJBINTRO} in a suitable integral form as follows.
\begin{equation}
  v(t,x) =R_{T-t}[\bar\phi](x_0)+\int_t^T R_{s-t}[H_{min}(
\nabla^B v(s,\cdot))+\bar\ell_0(s,\cdot)
](x)\; ds,\qquad t\in [0,T],\
x\in \calh,\label{solmildHJBINTRO}
\end{equation}
where by $R_{t},\,0\leq t\leq T$ we denote the transition semigroup
associated to the uncontrolled version of equation (\ref{eq-astrINTRO}).
In the literature such equation is treated by looking for a fixed point in suitable spaces
of functions (at least differentiable, in some sense, in $x$) using, as a key assumption,
a smoothing property of the transition semigroup $R_t$: i.e. that $R_t$ transforms any bounded measurable (or continuous) function $\varphi$ into a differentiable one with an integrable singularity of $R_t\varphi$ when $t \to 0^+$.
In this direction we mention the papers \cite{CDP1,CDP2,G1,G2} which treat the problem in spaces of continuous functions and the papers \cite{ChowMenaldi,GGSPA} that work in spaces $L^2$ with respect to a suitably defined invariant measure for the uncontrolled system.
In both cases, to get the key smoothing assumption on $R_t$ described just above, the authors need to assume strong properties of the data (a null controllability assumption when working in
spaces of continuous functions and again the structure condition when the live in $L^2$)
that cannot be verified in the case under study.

It may be possible to attack such HJB equations
with the viscosity solution approach (see e.g. \cite{GMSJOTA}
for an existence result in a similar case and \cite{FabbriGozziSwiech}
for a survey on the subject).
However this approach would not give us any regularity.
In the deterministic case there are papers
\cite{FedGolGozSICON1,FedGolGozSICON2,FedTacSicon}
where, as it is done in finite dimension (see \cite{BardiDolcetta,CannarsaSoner}),
a sort of $C^1$ regularity for the viscosity solution is proved through convexity
methods. However such approach does not seem to work, up to now, in the stochastic case.

\bigskip

In the present paper, in the spirit of the fixed point approach,
we start from the observation that, even if the semigroup $R_t$ in this case does not possess the above required smoothing property, it still regularizes some classes of functions in certain directions. In contrast to the standard smoothing property we call this last one {\em partial smoothing}.
More precisely, in the motivating case just described, we are able to prove that
for every bounded function $f:\calh\rightarrow\R$, of the form
\[
f\left(x_0, x_1\right)=\bar f(x_0),\quad\text{for some bounded function }\bar f:\R\rightarrow\R,
\]
the function $R_{t}[f]$ is continuous, and for every $x\in\calh$
the derivative $\nabla^B R_{t}[f](x)$ exists and blows up like $\sqrt{t}$
even if $B$ is unbounded
(see Section \ref{section-smoothOU}: Theorem \ref{lemma-der-gen} for the general case and
Proposition \ref{cor-der} and \ref{lemmaderhpdeb} for the application to
our specific case).

With this result at hand we act as in the fixed point approach, looking for solutions satisfying the integral equality \myref{solmildHJBINTRO}.
In the definition of the cost functional $\bar\phi$ and $\bar\ell_0$ depend on the state variable, which means that, in the abstract formulation, they depend only on the first component, so the partial smoothing property just stated fits to them. Then
we can carefully apply a fixed point technique (in doing which the choice of the spaces
must be done in a subtle way) to find a unique solution $v$ to the integral form \myref{solmildHJBINTRO} of the HJB equation (Section \ref{sec-HJB}, Theorem \ref{esistenzaHJB}).

On this solution $v$ (which we call mild solution according to most of the literature) we require the minimal regularity such that formula (\ref{solmildHJBINTRO}) makes sense.
Such minimal regularity is not enough to prove a verification theorem (which needs
to apply Ito's formula to $v$). However we are able to prove that $v$ is the limit
of classical solutions of approximating equations to which Ito's formula applies
(Section 7.1, Lemma \ref{lm:approximation}).
This is enough to prove the verification theorem (Section 7.2, Theorem \ref{teorema controllo}).

Finally adding some further regularity assumptions we prove that $v=V$ and we characterize the optimal control by a feedback law, which will be of the form
\begin{equation}\label{u in gammazeroIntro}
u(s)=\gamma\left(\nabla^B V(s,Y(s))\right),
\text{ \ \ \ \ }\mathbb{P}\text{-a.s.\ for almost every
}s\in\left[ t,T\right],
\end{equation}
for a suitable function $\gamma$ (Section 7.3, Theorem \ref{teo su controllo feedback}).

\bigskip

Notice that in the present paper we deal with a finite dimensional control delayed equation
(\ref{eq-contr-ritINTRO}), that here in the introduction we have presented in dimension one for the sake of simplicity. The same arguments apply if we consider the case of a controlled stochastic differential equation in an infinite dimensional Hilbert space $\calh_0$ with delay in the control as follows.
\begin{equation}
\left\{
\begin{array}
[c]{l}%
dy(t)  =A_0 y(t) dt+B_0 u(t) dt +\int_{-d}^0B_1(\xi)u(t+\xi)d\xi+\sigma dW_t
,\text{ \ \ \ }t\in[0,T] \\
y(0)  =y_0,\\
u(\xi)=u_0(\xi), \quad \xi \in [-d,0).
\end{array}
\right.  \label{eq-contr-rit-infINTRO}
\end{equation}
Here $W$ is a cylindrical Wiener process in another Hilbert space $\Xi$, $A_0$ is the generator of a strongly continuous semigroup in $\calh_0$, $\sigma \in \call (\Xi,\calh_0)$, and we have to assume some smoothing properties for the Ornstein Uhlenbeck transition
semigroup with drift term given by $A_0$ and diffusion equal to $\sigma$, see Remarks \ref{rm:diminf},  \ref{remarkinfdim1} and \ref{remark-derinfdim} for more details.

\subsection{Plan of the paper}

The plan of the paper is the following:
\begin{itemize}
  \item in Section \ref{section-statement} we give some notations and we present the problem and the main assumptions;
  \item in Section \ref{section-smoothOU} we prove the partial smoothing property for the Ornstein-Uhlenbeck transition semigroup, and we explain how to adapt it to an infinite dimensional setting;
  \item in Section \ref{section-smooth-conv} we introduce some spaces of functions where we will perform the fixed point argument and we prove regularity of some convolutions type integrals;
  \item  in Section \ref{sec-HJB}
we solve the HJB equation (\ref{HJBINTRO}) in mild sense;
  \item Section \ref{sec-verifica} is devoted to solve the optimal control problem. In Subsection 7.1 we approximate the mild solution of the HJB equation, in Subsection 7.2 we prove a verification theorem, and finally in Subsection 7.3
  we identify the value function of the control problem with the solution of the HJB equation and we characterize the optimal control by a feedback law.
\end{itemize}

\section{Preliminaries}\label{section-prel}

\subsection{Notation}\label{subsection-notation}

Let $H$ be a Hilbert space. The norm of an element $x$ in $H$ will be
denoted by $\left|  x\right|_{H}$ or simply $\left|  x\right|  $, if no
confusion is possible, and by $\left\langle \cdot,\cdot\right\rangle _H$, or simply
by $\left\langle \cdot,\cdot\right\rangle $ we denote the scalar product in $H$.
We denote by $H^{\ast}$ the dual
space of $H$.
If $K$ is another Hilbert space, $\call(H,K)$ denotes the
space of bounded linear operators from $H$ to $K$ endowed with the usual
operator norm. All Hilbert
spaces are assumed to be real and separable.

\noindent In what follows we will often meet inverses of operators which are not
one-to-one. Let $Q\in \call\left(H,K\right)  $. Then $H_{0}=\ker Q$ is a closed subspace of $H$. Let
$H_{1}$ be the orthogonal complement of $H_{0}$ in $H$: $H_{1}$ is
closed, too. Denote by $Q_{1}$ the restriction of $Q$ to $H_{1}:$ $Q_{1}$ is
one-to-one and $\operatorname{Im}Q_{1}=\operatorname{Im}Q$. For $k\in
\operatorname{Im}Q$, we define $Q^{-1}$ by setting
\[
Q^{-1}\left(  k\right)  :=Q_{1}^{-1}\left(  k\right)  .
\]
The operator $Q^{-1}:\operatorname{Im}Q\rightarrow H$ is called the
pseudoinverse of $Q$. $Q^{-1}$ is linear and closed but in general not
continuous. Note that if $k\in\operatorname{Im}Q$, then
$ Q_{1}^{-1}\left(  k\right)$ is the unique element of
\(\left\{ h  :Q\left(  h\right)  =k\right\}
\)
with minimal norm (see e.g. \cite{Z}, p.209).

In the following, by $(\Omega, \calf, \P)$ we denote a complete probability
space, and by $L^2_\calp(\Omega\times[0,T],H)$
the Hilbert space of all predictable processes
$(Z_t)_{t\in[0,T]}$ with values in $H$, normed by
$\Vert Z\Vert^2 _{L^2_\calp(\Omega\times[0,T],H)}=\E\int_0^T\vert Z_t\vert^2\,dt$.

Next we introduce some spaces of functions. We let $H$ and $K$ be Hilbert spaces.
By $B_b(H,K)$ (respectively $C_b(H,K)$, $UC_b(H,K)$) we denote the space of all functions
$f:H\rightarrow K$ which are Borel measurable and bounded (respectively continuous
and bounded, uniformly continuous and bounded).

Given an interval $I\subseteq \R$ we denote by
$C(I\times H,K)$ (respectively $C_b(I\times H,K)$)
the space of all functions $f:I \times H\rightarrow K$
which are continuous (respectively continuous and bounded).
$C^{0,1}(I\times H,K)$ is the space of functions
$ f\in C(I\times H)$ such that for all $t\in I$
$f(t,\cdot)$ is Fr\'echet differentiable.
By $UC_{b}^{1,2}(I\times H,K)$
we denote the linear space of the mappings $f:I\times H \to K$
which are uniformly continuous and bounded
together with their first time derivative $f_t$ and its first and second space
derivatives $\nabla f,\nabla^2f$.

If $K=\R$ we do not write it in all the above spaces.

\subsection{$C$-derivatives}\label{subsection-C-directionalderivatives}

We now introduce the $C$-directional derivatives following e.g.
\cite{FTGgrad} or \cite{Mas}, Section 2. Here $H$, $K$, $Z$ are Hilbert spaces.
\begin{definition}
\label{df4:Gder}
Let $C:K \rightarrow H$ be a bounded linear operator and let $f:H\rightarrow Z$.
\begin{itemize}
\item The $C$-directional
derivative $\nabla^{C}$ at a point $x\in H$ in the direction $k\in K$ is defined
as:
\begin{equation}
\nabla^{C}f(x;k)=\lim_{s\rightarrow 0}
\frac{f(x+s Ck)-f(x)}{s},\text{ }s\in\mathbb{R},
\label{Cderivata}
\end{equation}
provided that the limit exists.
\item We say that a continuous
function $f$ is $C$-G\^ateaux differentiable at a point $x\in H$
if $f$ admits the $C$-directional derivative in every direction $k\in K$ and
there exists a linear operator, called the $C$-G\^ateaux differential, $\nabla^{C}f(x)\in\call(K,Z)$, such that
$\nabla^{C}f(x;k)=\nabla^{C}f(x)k$ for $x \in H$, $k\in K$.
The function $f$ is $C$-G\^ateaux differentiable on $H$ if it is
$C$-G\^ateaux differentiable at every point $x\in H$.

\item We say that $f$ is $C$-Fr\'echet differentiable
at a point $x\in H$ if it is $C$-G\^ateaux differentiable and if the limit
in (\ref{Cderivata}) is uniform for $k$ in the unit ball of $K$. In this case
we call $\nabla^C f(x)$ the $C$-Fr\'echet derivative (or simply the $C$-derivative) of $f$ at $x$.
We say that $f$ is $C$-Fr\'echet differentiable on $H$ if it is $C$-Fr\'echet differentiable at every point $x\in H$.
\end{itemize}
\end{definition}

Note that, in doing the $C$-derivative, one considers only the directions
in $H$ selected by the image of $C$.
%
%
When $Z=\R$ we have $\nabla^C f(x) \in K^*$. Usually we will identify $K$ with its dual
$K^*$ so $\nabla^C f(x)$ will be treated as an element of $K$.

If $f:H\to \R$ is G\^ateaux (Fr\'echet) differentiable on $H$ we have that, given any $C$
as in the definition above, $f$ is $C$-G\^ateaux (Fr\'echet) differentiable on $H$
and
$$
\<\nabla^{C}f(x),k\>_{K}  =\<\nabla f(x),Ck\>_{H}
$$
i.e. the $C$-directional derivative is just the usual directional derivative at a point $x\in H$ in direction $Ck\in H$. Anyway the $C$-derivative,
as defined above, allows us to deal also with functions that are not G\^ateaux differentiable in every direction.

For our purposes (to treat the case of pointwise delay) we need to extend the concept of $C$-derivative to the case when $C:K\rightarrow H$ is a closed linear
and (possibly) unbounded operator with dense domain.
\begin{definition}
\label{df4:Gderunbounded}
Let $f:H\rightarrow Z$. Let $C:D(C)\subseteq K \to H$ be a closed linear (possibly)
unbounded operator with dense domain.
\begin{itemize}
\item
The $C$-directional derivative $\nabla^{C}f(x;k)$ at a point $x\in H$ in the direction
$k\in D(C)$ is defined exactly as in (\ref{Cderivata}).

  \item We say that $f$ is $C$-G\^ateaux differentiable
at a point $x\in H$ if $f$ admits the $C$-directional derivative in every
direction $k\in  D(C)$ and there exists a {\bf bounded} linear operator,
the $C$-G\^ateaux derivative $\nabla^C f(x)\in K^*$, such that $\nabla^{C}f(x;k)  =\nabla^{C}f(x)k$
for $x\in H$ and $k \in D(C)$. We say that $f$ is $C$-G\^ateaux
differentiable on $H$ if it is $C$-G\^ateaux differentiable at every point $x\in
H$.
  \item We say that $f$ is $C$-Fr\'echet differentiable
at a point $x\in H$ if it is $C$-G\^ateaux differentiable and if the limit
in (\ref{Cderivata}) is uniform for $k$ in the unit ball of $K$ intersected with $D(C)$. In this case
we call $\nabla^C f(x)$ the $C$-Fr\'echet derivative (or simply the $C$-derivative) of $f$ at $x$.
We say that $f$ is $C$-Fr\'echet differentiable on $H$ if it is $C$-Fr\'echet differentiable at every point $x\in H$.
\end{itemize}
\end{definition}

\medskip

\begin{remark}
\label{rm:Gderunbounded1}
In this last case, even if $f$ is Fr\'echet differentiable at $x \in H$, the $C$-derivative may not exist in such point.
Indeed consider the following case. Let $H$ be a Hilbert space and
$C:D(C)\subseteq H \to H$ be a closed linear operator on $H$ with dense domain and with unbounded adjoint $C^*$ on $H$ whose domain is $D(C^*)$.
Let $f:H \to \R$ be Fr\'echet differentiable at all $x \in H$ (e.g. $f(x)=|x|^2$).
By definition of $C$-directional derivative we have, for every $x \in H$ and $h\in D(C)$, that
$$\nabla^{C}f\left(  x;h\right)  =\<\nabla f\left(  x\right),Ch\>_H.
$$
On the other hand, if the $C$-derivative of $f$ exists at $x\in H$ then we should have $D^C f(x) \in \call(H;\R)=H^*$ (that we identify with $H$).
Hence, if $f$ was $C$-differentiable in all $H$ this would imply that, for any $x \in H$,
$$
|\nabla^{C}f(x;k)|=|\<\nabla^C f(x),k\>_H |\le c|k|,\quad \forall k\in D(C).
$$
This would mean that $\nabla f(x)\in D(C^*)$ for all $x \in H$, which cannot be true in our example.
\end{remark}

Now we define suitable spaces of $C$-differentiable functions.

\begin{definition}
\label{df4:Gspaces}
Let $I$ be an interval in $\R$ and let $H$, $K$, $Z$ be
suitable real Hilbert spaces.
\begin{itemize}
\item
We call $C^{1,C}_{b}(H,Z)$ the space of all functions in $C_b(H,Z)$ which
admit continuous and bounded $C$-Fr\'echet derivative. Moreover we call
$C^{0,1,C}_b(I\times H,Z)$ the space of functions
$f \in C_b(I\times H,Z)$ such that for every $t\in I$,
$f(t,\cdot )\in C^{1,C}_b(H,Z)$. When $Z=\R$ we omit it.

\item We call $C^{2,C}_{b}(H,Z)$ the space of all functions $f$ in $C^1_b(H,Z)$ which
admit continuous and bounded directional second order derivative $\nabla^C \nabla f$;
by $C^{0,2,C}_b(I\times H,Z)$ we denote the space of functions
$f \in C_b(I\times H,K)$ such that for every $t\in I$,
$f(t,\cdot )\in C^{2,C}_b(H,Z)$. When $Z=\R$ we omit it.

\item
For any $\alpha\in(0,1)$ and $T>0$ we denote by $C^{0,1,C}_{\alpha}([0,T]\times H)$ the space of functions
$ f\in C_b([0,T]\times H,Z)\cap C^{0,1,C}_b((0,T]\times H)$ such that
the map $(t,x)\mapsto t^{\alpha} \nabla^C f(t,x)$
belongs to $C_b((0,T]\times H,K^*)$.
The space $C^{0,1,C}_{\alpha}([0,T]\times H)$
is a Banach space when endowed with the norm
\[
 \left\Vert f\right\Vert _{C^{0,1,C}_{\alpha}([0,T]\times H)  }=\sup_{(t,x)\in[0,T]\times H}
\vert f(t,x)\vert+
\sup_{(t,x)\in (0,T]\times H}  t^{\alpha }\left\Vert \nabla^C f(t,x)\right\Vert_{K^{\ast}}.
\]
When clear from the context we will write simply
$\left\Vert f\right\Vert _{C^{0,1,C}_{\alpha}}$.
\item For any $\alpha\in(0,1)$ and $T>0$ we denote by $C^{0,2,C}_{\alpha}([0,T]\times H)$
the space of functions
$ f\in C_b([0,T]\times H)\cap C^{0,2,C}((0,T]\times H)$ such that for all $t\in(0,T],\,x\in H$ the map
$ (t,x)\mapsto t^{\alpha} \nabla^C\nabla f(t,x)$ is bounded and continuous as a map
from $(0,T]\times H$ with values in $H^*\times K^*$.
The space $C^{0,2,C}_{\alpha}([0,T]\times H)$ turns out to be a Banach space if it is endowed with the norm
\begin{align*}
\Vert  & f \Vert _{C^{0,2,C}_{\alpha}([0,T]\times H)  }\\
  &=\sup_{(t,x)\in[0,T]\times H}
  \vert f(t,x)\vert+
  \sup_{(t,x)\in[0,T]\times H}  \left\Vert \nabla f\left(  t,x\right)  \right\Vert _{H^*}
  +\sup_{(t,x)\in[0,T]\times H}  t^{\alpha}
  \left\Vert \nabla^C\nabla f\left(  t,x\right)  \right\Vert _{H^\ast\times K^{\ast}}.
\end{align*}
\end{itemize}
\end{definition}


%

\section{Setting of the problem and main assumptions}\label{section-statement}

\subsection{State equation}\label{subsection-stateequation}

In a complete probability space $(\Omega, \calf,  \P)$
we consider the following controlled stochastic
differential equation in $\R^n$ with delay in the control:
\begin{equation}
\left\{
\begin{array}
[c]{l}%
dy(t)  =a_0 y(t) dt+b_0 u(t) dt +\int_{-d}^0b_1(\xi)u(t+\xi)d\xi+\sigma dW_t
,\text{ \ \ \ }t\in[0,T] \\
y(0)  =y_0,\\
u(\xi)=u_0(\xi), \quad \xi \in [-d,0),
\end{array}
\right.  \label{eq-contr-rit}
\end{equation}
where we assume the following.

\begin{hypothesis}\label{ipotesibasic}
\begin{itemize}
\item[]
  \item[(i)] $W$ is a standard Brownian motion in $\R^k$, and $(\calf_t)_{t\geq 0}$ is the
augmented filtration generated by $W$;
  \item[(ii)] $a_0\in \call(\R^n;\R^n)$, $\sigma$ is in $\call(\R^k;\R^n)$;
  \item[(iii)] the control strategy $u$ belongs to $\calu$ where
$$\calu:=\left\lbrace z\in L^2_{\calp}(\Omega\times [0,T], \R^m):
u(t)\in U \;a.s.\right\rbrace $$
where $U$ is a closed subset of $\R^n$;
  \item[(iv)] $d>0$ (the maximum delay the control takes to affect the system);
  \item[(v)] $b_0 \in \call(\R^m;\R^n)$;
  \item[(vi)]
$b_1$ is an $m\times n$ matrix of signed Radon measures on $[-d,0]$
(i.e. it is an element of the dual space of $C([-d,0],\call(\R^m;\R^n))$
and the integral
$\int_{-d}^0b_1(\xi)u(t+\xi)d\xi$ is understood in the sense
$\int_{-d}^0 b_1(d\xi) u(t+\xi)$
($b_1$ is the density of the time taken by the control to affect the system).
\end{itemize}
\end{hypothesis}

\medskip

Notice that  assumption (vi) on $b_1$ covers the
very common case of pointwise delay but it is technically complicated to
deal with: indeed it gives rise, as we are going to see in next subsection,
to an unbounded control operator $B$.
In some cases we will use the following more restrictive
assumption on the term $b_1$ which leaves aside the pointwise delay case:
\begin{equation}\label{eq:b1restrictive}
b_1\in L^2([-d,0],\call(\R^m;\R^n)).
\end{equation}

\begin{remark}\label{rm:diminf}
Our results can be generalized to the case when the process $y$ is infinite dimensional.
More precisely, let $y$ be the solution of the following controlled stochastic
differential equation in an infinite dimensional Hilbert space $H$, with delay in the control:
\begin{equation}
\left\{
\begin{array}
[c]{l}%
dy(t)  =A_0 y(t) dt+B_0 u(t) dt +\int_{-d}^0B_1(\xi)u(t+\xi)d\xi+\sigma dW_t
,\text{ \ \ \ }t\in[0,T] \\
y(0)  =y_0,\\
u(\xi)=u_0(\xi), \quad \xi \in [-d,0).
\end{array}
\right.  \label{eq-contr-rit-inf}
\end{equation}
Here $W$ is a cylindrical Wiener process in another Hilbert space $\Xi$, and $(\calf_t)_{t\geq 0}$ is the augmented filtration generated by $W$.
$A_0$ is the generator of a strongly continuous semigroup in $H$.
The diffusion term $\sigma$ is in $\call(\Xi;H)$ and is such that for every $t>0$
the covariance operator
\[
\int_{0}^{t}e^{sA_0}\sigma\sigma^*e^{sA_0^{\ast}}ds
\]
of the stochastic convolution
\[
\int_0^t e^{(t-s)A_0}\sigma\,dW_s
\]
is of trace class and, for some $\gamma \in (0,1)$,
$$
\int_{0}^{t}s^{-\gamma}{\rm Tr}e^{sA_0}\sigma\sigma^*e^{sA_0^{\ast}}ds < + \infty.
$$
The control strategy $u$ belongs to $L^2_{\calp}(\Omega\times [0,T], U_1)$,
where $U_1$ is another Hilbert space, and the space of admissible controls $\calu$ is built in analogy with the finite dimensional case requiring control strategies to take values in a given closed subset $U$ of $U_1$.
On the control operators we assume $B_0 \in \call(U_1;H)$, $B_1: [-d,0] \to\call(U_1,H))$
such that $B_1 u \in L^2([-d,0],H)$ for all $u \in U$.\footnote{It seems possible also in this case, to cover the case when $B_1$ is a measure, on the line of Hypothesis \ref{ipotesibasic}-(vi), but we do not do this here for simplicity.}
In this case, following again \cite{GM,VK}, the problem can be reformulated as an abstract evolution equation
in the Hilbert space $\calh$ that this time turns out to be $H\times L^2([-d,0],H)$.
All the results of this paper hold true in this case, under suitable minor changes that will be clarified along the way.
\hfill\qedo
\end{remark}

\subsection{Infinite dimensional reformulation}
\label{subsection-infdimref}

Now, using the approach of \cite{VK} (see \cite{GM} for the stochastic case), we reformulate
equation (\ref{eq-contr-rit}) as an abstract stochastic
differential equation in the Hilbert space $\calh=\R^n\times L^2([-d,0],\R^n)$.
To this end we introduce the operator $A : \cald(A) \subset \calh
\rightarrow \calh$ as follows: for $(y_0,y_1)\in \calh$
\begin{equation}\label{A}
A(y_0 ,y_1 )=( a_0 y_0 +y_1(0), -y_1'), \qquad
\cald(A)=\left\lbrace(y_0,y_1)\in \calh:y_1\in W^{1,2}([-d,0],\R^n), y_1(-d)=0 \right\rbrace.
\end{equation}
We denote by $A^*$ the adjoint operator of $A$:
\begin{equation}
 \label{Astar}
A^{*}(y_0 ,y_1 )=( a_0 y_0, y_1'), \qquad
\cald(A^{*})=\left\lbrace (y_0,y_1)\in \calh:y_1\in W^{1,2}([-d,0],\R^n), y_1(0)=y_0 \right\rbrace .
\end{equation}
We denote by $e^{tA}$ the $C_0$-semigroup generated by $A$: for
$y=(y_0,y_1)\in \calh$,
\begin{equation}
e^{tA} \left(\begin{array}{l}y_0 \\y_1\end{array}\right)=
\left(
\begin{array}
[c]{ll}%
e^{ta_0 }y_0+\int_{-d}^{0}1_{[-t,0]} e^{(t+s)a_0 } y_1(s)ds \\[3mm]
y_1(\cdot-t)1_{[-d+t,0]}(\cdot).
\end{array}
\right)  \label{semigroup}
\end{equation}
We will use, for $N\in \N$ big enough, the resolvent operator
$(N-A)^{-1}$ which can be computed explicitly
giving
\begin{equation}
(N-A)^{-1} \left(\begin{array}{l}y_0 \\y_1\end{array}\right)=
\left(
\begin{array}
[c]{ll}%
(N-a_0)^{-1}\left[y_0+\int_{-d}^{0} e^{N s } y_1(s)ds\right]
\\[3mm]
\int_{-d}^{\cdot} e^{N(s-\cdot)} y_1(s)ds.
\end{array}
\right)  \label{resolvent}
\end{equation}
Similarly, denoting by $e^{tA^*}=(e^{tA})^*$ the $C_0$-semigroup generated by $A^*$,
we have for
$z=\left(z_0,z_1\right)\in \calh $
\begin{equation}
e^{tA^*} \left(\begin{array}{l}z_0 \\z_1\end{array}\right)=
\left(
\begin{array}[c]{ll}
e^{t a_0^* }z_0 \\[3mm]
e^{(\cdot+t) a_0^* }z_0 1_{[-t,0]}(\cdot) +z_1(\cdot+t)1_{[-d,-t)}(\cdot).
\end{array}
\right)  \label{semigroupadjoint}
\end{equation}
The infinite dimensional noise operator is defined as
\begin{equation}
 \label{G}
G:\R^{k}\rightarrow \calh,\qquad Gy=(\sigma y, 0), \; y\in\R^k.
\end{equation}
The control operator $B$, when \myref{eq:b1restrictive} holds, is
bounded and defined as
\begin{equation}
 \label{B}
B:\R^{m}\rightarrow \calh,\qquad Bu=(b_0 u, b_1(\cdot)u), \; u\in\R^m
\end{equation}
and its adjoint is
\begin{equation}
 \label{B*}
B^*:\calh^* \rightarrow \R^{m},\qquad B^*(x_0,x_1)=
b^*_0 x_0+\int_{-d}^0 b_1(\xi)^*x_1(\xi)d\xi, \; (x_0,x_1)\in\calh.
\end{equation}
In the general case $B$ is possibly unbounded
and defined as (here we denote by $C^*([-d,0],\R^n)$
the dual space of $C([-d,0],\R^n)$)
\begin{equation}
 \label{Bnotbdd}
\begin{array}{c}
B:\R^{m}\rightarrow \R^n \times C^*([-d,0],\R^n),
   \\[2mm]
(Bu)_0=b_0 u,
\quad \<f,(Bu)_1\>_{C,C^*} = \int_{-d}^0 f(\xi)b_1(d\xi )u , \quad u\in\R^m, \quad
f \in C([-d,0],\R^n).
\end{array}
\end{equation}
In this case the adjoint $B^*$ is (here we denote by $C^{**}([-d,0],\R^n)$
the dual space of $C^*([-d,0],\R^n)$, which contains $C([-d,0],\R^n)$)
\begin{equation}
 \label{B*notbdd}
\begin{array}{c}
B^*:\R^n \times C^{**}([-d,0],\R^n) \rightarrow \R^{m},
   \\[2mm]
B^*(x_0,x_1)=
b^*_0 x_0+\int_{-d}^0 b_1(d\xi)^*x_1(\xi)d\xi,
\; (x_0,x_1)\in \R^n \times C([-d,0],\R^n).
\end{array}
\end{equation}
It will be useful to write the explicit expression of the first component
of the operator $e^{tA}B$ as follows
\begin{equation}\label{eq:etAB}
\left(e^{tA}B\right)_0:\R^m \to \R^n,\qquad    \left(e^{tA}B\right)_0 u=
    e^{ta_0}b_0u+ \int_{-d}^0 1_{[-t,0]}e^{(t+r)a_0}b_1(dr)u, \quad u \in \R^m.
\end{equation}
%
%

Given any initial datum $(y_0,u_0)\in \calh$ and any admissible control $u\in \calu$ we call $y(t;y_0,u_0,u)$ (or simply $y(t)$ when clear from the context) the unique solution (which comes from standard results on SDE's, see e.g.  \cite{IkedaWatanabe} Chapter 4, Sections 2 and 3)
of (\ref{eq-contr-rit}).

Let us now define the process
$Y=(Y_0,Y_1)\in L^2_\calp(\Omega \times [0,T],\calh)$ as
$$
Y_0(t)=y(t), \qquad Y_1(t)(\xi)=\int_{-d}^\xi b_1(d\zeta)u(\zeta+t-\xi),
$$
where $y$ is the solution of equation (\ref{eq-contr-rit}), $u$ is the control process in (\ref{eq-contr-rit}) and $b_1(d\zeta)$ is understood as
$b_1(\zeta)d\zeta$ when \myref{eq:b1restrictive} holds.
By Proposition 2 of \cite{GM},
the process $Y$
is the unique solution of the abstract evolution equation
in $\calh$
\begin{equation}
\left\{
\begin{array}
[c]{l}
dY(t)  =AY(t) dt+Bu(t) dt+GdW_t
,\text{ \ \ \ }t\in[ 0,T] \\
Y(0)  =y=(y_0,y_1),
\end{array}
\right.   \label{eq-astr}%
\end{equation}
where $y_0=x_0$ and $y_1(\xi)=\int_{-d}^\xi b_1(d\zeta)u_0(\zeta-\xi)$.
Note that we have $y_1\in L^2([-d,0];\R^n)$\footnote{This can be seen, e.g.,
by a simple application of Jensen inequality and Fubini theorem. }.
Taking the integral (or mild) form of (\ref{eq-astr}) we have
\begin{equation}
Y(t)  =e^{tA}y+\int_0^te^{(t-s)A}B u(s) ds +\int_0^te^{(t-s)A}GdW_s
,\text{ \ \ \ }t\in[ 0,T]. \\
  \label{eq-astr-mild}%
\end{equation}

\subsection{Optimal Control problem}


The objective is to minimize, over all controls in $\calu$,
the following finite horizon cost:
 \begin{equation}\label{costoconcreto}
J(t,x,u)=\E \int_t^T \left(\bar\ell_0(s,y(s))+\bar\ell_1(u(s))\right)\;ds +\E  \bar\phi(x(T)).
\end{equation}
where $\bar\ell_0:[0,T]\times\R^n\rightarrow \R$ and
$\bar\phi:\R^n\rightarrow \R$ are continuous and bounded while
$\bar\ell_1:U\rightarrow\R$ is measurable and bounded from below.
Referring to the abstract formulation (\ref{eq-astr}) the cost in (\ref{costoconcreto}) can be rewritten also as
\begin{equation}\label{costoconcreto1}
J(t,x;u)=\E \left(\int_t^T \left[\ell_0(s,Y(s))+\ell_1(u(s))\right]\,ds + \phi(Y(T))\right),
\end{equation}
where
$\ell_0:[0,T]\times \calh\rightarrow \R$, $\ell_1:U\to \R$ are defined by setting
\begin{equation}\label{l_0}
\ell_0(t,x):=\bar\ell_0(t,x_0) \quad
\forall x=(x_0,x_1)\in \calh
\end{equation}
\begin{equation}\label{l_1}
\ell_1:=\bar\ell_1
\end{equation}
(here we cut the bar only to keep the notation homogeneous)
while $\phi :\calh\rightarrow \R$ is defined as
\begin{equation}\label{fi0}
\phi(x):=\bar\phi(x_0) \quad
\forall x=(x_0,x_1)\in \calh.
\end{equation}

Clearly, under the assumption above,
$\ell_0$ and $\phi$
are continuous and bounded while $\ell_1$ is measurable and bounded from below.
The value function of the problem is
\begin{equation}\label{valuefunction}
 V(t,x):= \inf_{u \in \calu}J(t,x;u).
\end{equation}
As done in Subsection \ref{SS:motiv}, we define the Hamiltonian in a modified way
(see \myref{HminINTRObis}); indeed, for $p\in \calh$, $u \in U$,
we define the current value Hamiltonian $H_{CV}$ as
$$
H_{CV}(p\,;u):=\<p,u\>_{\R^m}+\ell_1(u)
$$
and the (minimum value) Hamiltonian by
\begin{equation}\label{psi1}
H_{min}(p)=\inf_{u\in U}H_{CV}(p\,;u),
 \end{equation}
The associated HJB equation with unknown $v$ is then formally written as
\begin{equation}\label{HJBformale1}
  \left\{\begin{array}{l}\dis
-\frac{\partial v(t,x)}{\partial t}=\frac{1}{2}Tr \;GG^*\nabla^2v(t,x)
+ \< Ax,\nabla v(t,x)\>_\calh +\ell_0(t,x)+ {H}_{min} (\nabla^Bv(t,x)),\qquad t\in [0,T],\,
x\in D(A),\\
\\
\dis v(T,x)=\phi(x).
\end{array}\right.
\end{equation}

To get existence of mild solutions of \myref{HJBformale1} we will need the following assumption.

\begin{hypothesis}\label{ipotesicostoconcreto}
\begin{itemize}
\item[]
  \item[(i)] $\phi\in C_b(\calh)$ and it is given by \myref{fi0} for
  a suitable $\phi \in  C_b(\R^n)$;
  \item[(ii)] $\ell_0\in C_b([0,T] \times\calh)$ and it is given by \myref{l_0} for
  a suitable $\bar\ell_0 \in  C_b([0,T] \times \R^n)$;
  \item[(iii)] $\ell_1:U\rightarrow\R$ is measurable and bounded from below;
  \item[(iv)] the Hamiltonian $H_{min}:\R^m \to \R$ is Lipschitz continuous so
  there exists $L>0$ such that
  \begin{equation}\label{eq:Hlip}
    \begin{array}{c}
\vert H_{min }(p_1)-H_{min }(p_2)\vert\leq L \vert p_1-p_2\vert
\quad \forall\,p_1,\,p_2\in\R^m;
       \\[1.5mm]
\vert H_{min }(p)\vert\leq L(1 + \vert p\vert )
\quad \forall\,p\in\R^m.
    \end{array}
\end{equation}
\end{itemize}
\end{hypothesis}
To get more regular solutions (well defined second derivative $\nabla^B\nabla$,
which will be used to prove existence of optimal feedback controls) we will need
the following further assumption.

\begin{hypothesis}\label{ipotesicostoconcretobis}
\begin{itemize}
\item[]
\item[(i)] $\ell_0$ is continuously differentiable in the variable $x$
  with bounded derivative.
  \item[(ii)] the Hamiltonian $H_{min}:\R^m \to \R$ is continuously differentiable
  and, for a given $L>0$, we have, beyond \myref{eq:Hlip},
  \begin{equation}\label{eq:Hlipder}
    \begin{array}{c}
\vert \nabla H_{min }(p_1)-\nabla H_{min }(p_2)\vert\leq L \vert p_1-p_2\vert
\quad \forall\,p_1,\,p_2\in\R^m;
    \end{array}
\end{equation}
\end{itemize}
\end{hypothesis}

\begin{remark}
The assumption \myref{eq:Hlip} of Lipschitz continuity of $H_{min }$ is satisfied e.g. if the set $U$ is compact.
Indeed, for every $p_1,\,p_2\in\R^m$
\[
 \vert H_{min}(p_1)-H_{min}(p_2)\vert \leq \vert \<p_1,u\>-\<p_2,u\>\vert,\quad u\in U
\]
and in the case of $U$ compact the Lipschitz property immediately follows.
The Lipschitz continuity of $H_{min }$ is satisfied also in the case
when $U$ is unbounded, if the current cost has linear growth at infinity.

Moreover the assumption \myref{eq:Hlipder} of Lipschitz continuity of $\nabla H_{min }$ is verified e.g. if the function $\ell_1$ is convex, differentiable with invertible derivative and with $(\ell_1')^{-1}$ Lipschitz continuous since in this case
$(\ell_1')^{-1}(p)=\nabla H_{min}(p)$.
\hfill\qedo
\end{remark}

\begin{remark}\label{remark:costoconcreto}
We list here, in order of increasing difficulty,
some possible generalization of the above assumptions and of the consequent results.
\begin{itemize}
  \item[(i)]
All our results on the HJB equation and on the control problem
could be extended without difficulties to the case
when the boundedness assumption on $\bar \phi$ and $\bar \ell_0$
(and consequently on $\phi$ and $\ell_0$) can be replaced by a polynomial growth assumption:
namely that, for some $N\in \N$, the functions
\begin{equation}\label{eq:polgrowthphil0}
 x\mapsto \dfrac{\phi(x)}{1+\vert x\vert ^N},
 \qquad
 (t,x)\mapsto \dfrac{\ell_0(t,x)}{1+\vert x\vert ^N},
\end{equation}
are bounded.
The generalization of Theorem \ref{lemma-der-gen} to this case can be achieved
by straightforward changes in the proof, on the line of what is done, in a different context, in
\cite{Ce95} or in \cite{Mas-inf-or}.
  \item[(ii)]
Since our results on the HJB equation are based on smoothing properties
(proved in Section \ref{section-smoothOU}) which holds also for
measurable functions, we could consider current cost and final cost only measurable
instead of continuous. The proofs would be very similar but using different
underlying spaces.
  \item[(iii)]
Using the approach of \cite{G2} it seems possible to relax the Lipschitz assumptions on the Hamiltonian function asking only local Lipschitz continuity of the Hamiltonian function,
but paying the price of requiring differentiability of the data.
\end{itemize}
In this paper we do not perform all such generalizations
since we want to concentrate on the main point: {\em the possibility of solving
the HJB equation and the control problem without
requiring the so-called {\bf structure condition}.}
\end{remark}

\section{Partial smoothing for the Ornstein-Uhlenbeck
semigroup}
\label{section-smoothOU}
This section is devoted to what we call the ``partial'' smoothing property
of Ornstein-Uhlenbeck transition semigroup in a general Hilbert space.
First, in Subsection \ref{section-smoothOU-generalsetting}, we give two
general results on it (Theorem \ref{lemma-der-gen} for the first $C$-derivative
and Proposition \ref{lemma-reg-R_t} for the second derivative); then in Subsection \ref{section-smoothOU-particular}, we prove two specific results for our problem
(Propositions \ref{cor-der} and \ref{lemmaderhpdeb}).

\subsection{Partial smoothing in a general setting}
\label{section-smoothOU-generalsetting}
Let $H, \Xi$ be two real and separable Hilbert spaces and let us consider
the Ornstein-Uhlenbeck process $X^x(\cdot)$ in $H$ which solves the
following SDE in $H$:
\begin{equation}
\left\{
\begin{array}
[c]{l}%
dX(t)  =AX(t) dt+GdW_t
,\text{ \ \ \ }t\ge 0\\
X(0)  =x,
\end{array}
\right.\label{ornstein-gen}
\end{equation}
where $A$ is the generator of a strongly continuous semigroup in
$H$, $(W_t)_{t\geq 0}$ is a cylindrical Wiener process in $\Xi$ and
$G:\Xi\rightarrow K$.
In mild form, the Ornstein-Uhlenbeck process $X^x$ is given by
\begin{equation}
X^x(t)  =e^{tA}x +\int_0^te^{(t-s)A}GdW_s
,\text{ \ \ \ }t\ge 0. \\
  \label{ornstein-mild-gen}
\end{equation}
$X$ is a Gaussian process, namely for every $t>0$, the law of
$X(t)$ is $\caln (e^{tA}x,Q_t)$, the Gaussian measure with mean $e^{tA}x$ and
covariance operator $Q_t$,
where
\[
 Q_t=\int_0^t e^{sA}GG^*e^{sA^*}ds.
\]
The associated Ornstein-Uhlenbeck transition semigroup $R_t$, is defined by setting,
for every $f\in B_b(H)$ and $x\in H$,
\begin{equation}
 \label{ornstein-sem-gen}
R_t[f](x)=\E f(X^x(t))
=\int_K f(z+e^{tA}x)\caln(0,Q_t)(dz).
\end{equation}
where by $X^x$ we denote the Ornstein-Uhlenbeck process above with initial datum
given by $x\in H$.

It is well known (see e.g. \cite[Section 9.4]{DP1}),
that $R_t$ has the strong Feller property (i.e. it transforms
bounded measurable functions in continuous ones) if and only if
\begin{equation}\label{eq:inclusionDZ}
\operatorname{Im}e^{tA}\subseteq \operatorname{Im} Q_t^{1/2},
\end{equation}
and that such property is equivalent to the so-called null-controllability
of the linear control system identified by the couple of operators
$(A,G)$ (here $z(\cdot)$ is the state and $a(\cdot)$
is the control):
$$
z'(t)=Az(t)+Ga(t),\qquad z(0)=x.
$$
(see again \cite[Appendix B]{DP1}).
Under \myref{eq:inclusionDZ} $R_t$ also transforms any bounded measurable
function $f$ into a Fr\'echet differentiable one,
the so-called ``smoothing'' property,
and
$$
\|\nabla R_t[f]\|_\infty \le \|\Gamma(t)\|_{\call(H)}\|f\|_\infty
$$
where $\Gamma(t):=Q_t^{-1/2}e^{tA}$.

Here we extend this property in two directions: searching for $C$-derivatives and
applying $R_t$ to a specific class of bounded measurable functions.
(see \cite{Lu} for results in this direction in finite dimension).


Now, let $P:H\rightarrow H$ be a bounded linear operator;
given any $\bar\phi : \operatorname{Im}(P)\to \R$
measurable and bounded we define a function $\phi \in B_b(H)$, by setting
\begin{equation}\label{fi-gen}
\phi(x)=\bar\phi(Px) \quad
\forall x\in H.
\end{equation}
We now prove that, under further assumptions on the operators $A$, $G$, $C$ and $P$,
the semigroup $R_t$ maps functions $\phi$, defined as in (\ref{fi-gen}),
into $C$-Fr\'echet differentiable functions.
%
%
\begin{theorem}\label{lemma-der-gen} Let $A$ be the generator of a strongly continuous semigroup in $H$ and $G:\Xi\rightarrow H$.
Let $K$ be another real and separable Hilbert space and let
 $C:D(C)\subseteq K\rightarrow H$ be a linear (possibly) unbounded operator with
dense domain. Let
$\bar\phi:\operatorname{Im}(P)\rightarrow\R$ be measurable and bounded
and define $\phi :H\rightarrow \R$ as in (\ref{fi-gen}).
Fix $t>0$ and assume that $Pe^{tA}C: K\rightarrow H$ is well defined.

Then,
$R_t[\phi]$ is $C$-Fr\'echet differentiable if
\begin{equation}\label{eq:inclusion-iff-P}
 \operatorname{Im} \left(P e^{tA}C\right)\subseteq \operatorname{Im} Q_t^{1/2}.
\end{equation}
In this case we have, for every $k\in K$,
\begin{align}
\langle\nabla^C(R_{t}\left[\phi\right]  )(x)  ,k \rangle_K =
\int_{H}
\bar\phi\left(Pz+Pe^{ tA}x\right)
\<Q_{t}^{-1/2}P e^{tA}C k, Q_t^{-1/2}z\>_H\caln(0,Q_t)(dz).
\label{eq:formulaDRT1-gen-P}
\end{align}
Moreover for every $k\in K$
we have the estimate
\begin{equation*}
\vert \<\nabla^C(R_{t}\left[\phi\right]  )(x)  k\>_K\vert
\leq \Vert \bar\phi\Vert_\infty
\left\Vert Q_{t}^{-1/2}P e^{tA} C
\right\Vert_{\call(K;H)} \;\vert k\vert_K,
\end{equation*}
\end{theorem}


\dim
Let $k\in D(C) \subseteq K$. We compute the directional derivative
in the direction $Ck$.
\begin{align*}
& \lim_{\alpha\rightarrow0}\frac{1}{\alpha}
\left[R_{t}[\phi](x+\alpha Ck)-R_{t}[\phi](x)\right]
\\[2mm]
&  =\lim_{\alpha\rightarrow 0}\frac{1}{\alpha}
\left[\int_{H}\bar\phi\left(Pz+Pe^{tA}(x+\alpha Ck)\right)
\caln\left(0,Q_{t}\right)(dz)
-\int_{H}\bar\phi\left(Pz+Pe^{tA}x\right)\caln\left(0,Q_{t}\right)(dz)\right]
\\[2mm]
&  =\lim_{\alpha\rightarrow 0}\frac{1}{\alpha}
\left[\int_{H}\bar\phi\left(Pz+Pe^{tA}x\right)
\caln\left(Pe^{tA}\alpha Ck ,Q_{t}\right)(dz)
-\int_{H}\bar\phi\left(Pz+Pe^{tA}x\right)\caln\left(0,Q_{t}\right)(dz)\right].
\end{align*}
By the Cameron-Martin theorem, see
e.g. \cite{DP3}, Theorem 1.3.6, the Gaussian measures
$\caln\left(Pe^{tA}\alpha Ck,Q_t\right)$ and
$\mathcal{N}\left(0,Q_{t}\right)$ are equivalent if and only if
$Pe^{tA}\alpha Ck\in\operatorname{Im}Q_t^{1/2}$.
In such case,
setting, for $y \in \operatorname{Im}Q_t^{1/2}$,
\begin{align}
&d(t,y,z)     =\frac{d\caln\left(y,Q_{t}\right)}
{d\mathcal{N}\left(0,Q_{t}\right)  }(z)
  =\exp\left\{  \left\langle Q_{t}^{-1/2}
y,Q_{t}^{-1/2}z\right\rangle_H
-\frac{1}{2}\left|Q_{t}^{-1/2}y\right|_H^{2}\right\}  ,
\label{eq:density1}
\end{align}
we have, arguing exactly as in \cite{DP1}, proof of Theorem 9.26,
\begin{align*}
&  \nabla^C(R_{t}[\phi])(x)k=
\lim_{\alpha\rightarrow 0}\frac{1}{\alpha}
\int_{H}\bar\phi\left(Pz+Pe^{tA}x\right)\frac{d(t,Pe^{tA}\alpha Ck,z)-1}{\alpha}
\caln(0,Q_t)(dz)
\\[2mm]
&  =\int_{H}\bar\phi\left(Pz+Pe^{tA}x\right)
\<Q_{t}^{-1/2}P e^{tA}C k, Q_t ^{-1/2}z\>_H\caln(0,Q_t)(dz)
\end{align*}
which gives (\ref{eq:formulaDRT1-gen-P}).
Consequently
\begin{align*}
 \vert \<\nabla^C(R_{t}\left[\phi\right]  )(x) , k\>_K\vert
&\leq
\Vert \bar\phi\Vert_\infty
\left(\int_{H}\<Q_{t}^{-1/2}P e^{tA} Ck,
Q_{t}^{-1/2}z\>^2\caln(0,Q_t)(dz)\right)^{1/2}
= \Vert \bar\phi\Vert_\infty
\Vert Q_{t}^{-1/2} P e^{tA}Ck \Vert_{\call(K;H)}. \nonumber
\end{align*}
This gives the claim.
\qed

\begin{remark}\label{remark-gen0}
In \cite{DP1}, Remark 9.29, it is showed that the analogous of condition \ref{eq:inclusion-iff-P}
is also a necessary condition for the Fr\'echet differentiability of $R_t[\phi]$ for any bounded Borel $\phi$. In this case this is not obvious as we deal with a special class of $\phi$ and so the counterexample provided in such Remark may not belong to this class.
\end{remark}

\begin{remark}\label{remark-gen1}
We consider two special cases of the previous Theorem \ref{lemma-der-gen} that will be useful in next section.
\begin{itemize}
  \item[(i)]
Let $K=H$ and $C=I$. In this case Theorem \ref{lemma-der-gen} gives Fr\'echet differentiability:
for $t>0$ $R_t[\phi]$ is Fr\'echet differentiable if
\begin{equation}\label{eq:inclusion-iff-tutta}
 \operatorname{Im}\left( P e^{tA}\right)\in  \operatorname{Im} Q_t^{1/2}
\end{equation}
 and we have,
for every $h\in H$,
\begin{align}
\langle\nabla(R_{t}[\phi])(x),h \rangle_H =
\int_{H}\bar\phi\left(Pz+Pe^{ tA}x\right)
\<Q_{t}^{-1/2}Pe^{tA} h, Q_t^{-1/2}z\>_H\caln(0,Q_t)(dz).
\label{eq:formulaDRT1-gen-tutta}
\end{align}
Moreover for every $h\in H$
we have the estimate
\begin{equation*}
\vert\< \nabla(R_{t}[\phi])(x),h\>_H\vert
\leq \Vert \bar\phi\Vert_\infty\left\Vert Q_{t}^{-1/2}Pe^{tA}
\right\Vert_{\call(H;H)} \;\vert h\vert_H,
\end{equation*}
\item[(ii)]
Let $K_0$ and $K_1$ be two real and separable Hilbert spaces and let
$K=K_0\times K_1$ be the product space.
Now, given any $\bar\phi\in B_b(K_0)$, we define, in the same way as in (\ref{fi0}),
a function $\phi \in B_b(K)$, by setting
\begin{equation}\label{fi-genbis}
\phi(k)=\bar\phi(k_0) \quad
\forall k=\left(k_0, k_1 \right)\in K.
\end{equation}
Let $P:K\rightarrow K_0$ be the projection on the first component of $K$: for every $k=\left(k_0, k_1 \right)\in K$, $Pk=k_0$.
Theorem \ref{lemma-der-gen} says that
$R_t[\phi]$ is $C$-Fr\'echet differentiable for every $t>0$ if
\begin{equation}\label{eq:inclusion-iff}
 \operatorname{Im}\left(\left( e^{tA}C\right)_0,0\right)\in  \operatorname{Im} Q_t^{1/2}, \quad \forall t>0
\end{equation}
 and we have, for every $k\in K$,
\begin{align}
\langle\nabla^C(R_{t}\left[\phi\right]  )(x)  ,k \rangle_K =
\int_{K}
\bar\phi\left(z_0+(e^{ tA}x)_0\right)
\<Q_{t}^{-1/2}\left(\left( e^{tA} Ck\right)_0,0\right),Q_t^{-1/2}z\>_K\caln(0,Q_t)(dz).
\label{eq:formulaDRT1-gen}
\end{align}
Moreover for every $k\in K$
we have the estimate
\begin{equation}\label{eq:stimaDRT1-gen}
\vert \<\nabla^C(R_{t}\left[\phi\right])(x),k\>_K\vert
\leq \Vert \bar\phi\Vert_\infty
\left\Vert Q_{t}^{-1/2}
\left(\left( e^{tA} C\right)_0,0
\right)
\right\Vert_{\call(K;K)} \;\vert k\vert_K,
\end{equation}
\end{itemize}
\end{remark}

\begin{remark}\label{rm:crescitapoli-gen}
In Theorem \ref{lemma-der-gen} we prove the partial smoothing for functions $\phi$
defined as in \myref{fi-gen} for functions $\bar \phi$
bounded and measurable. The boundedness assumption on $\bar \phi$
(and consequently on $\phi$) can be replaced by a polynomial growth assumption:
namely that, for some $N\in \N$,
\[
 x\mapsto \dfrac{\bar\phi(x)}{1+\vert x\vert ^N}
\]
is bounded.
The generalization of Theorem \ref{lemma-der-gen} to this case can be achieved
by straightforward changes in the proof, on the line of what is done in
\cite{Ce95} or in \cite{Mas-inf-or}.
\end{remark}

\subsubsection{Second derivatives}

We now prove that, if $\phi$ is more regular,
also $\nabla^C R_t[\phi]$ and $\nabla R_t[\phi]$ have more regularity.
This fact, in the context of our model (see Subsection \ref{section-smoothOU-particular}),
will be used in Section \ref{sec-HJB} to prove the Verification Theorem
\ref{teorema controllo} and the existence of optimal controls in feedback form in
Theorem \ref{teo su controllo feedback}.


\begin{proposition}\label{lemma-reg-R_t} Let $A$ be the generator of a strongly continuous semigroup in $H$ and $G:\Xi\rightarrow H$. Let $K$ be another real and separable Hilbert space and let
 $C:D(C)\subseteq K\rightarrow H$ be a linear (possibly) unbounded operator with
dense domain. Let
$\bar\phi:\operatorname{Im}(P)\rightarrow\R$ be measurable and bounded
and define $\phi :H\rightarrow \R$ as in (\ref{fi-gen}).
Fix $t>0$ and assume that $Pe^{tA}C: K\rightarrow H$ is well defined. Assume that (\ref{eq:inclusion-iff-P})
holds true.
If $\bar\phi$ is such that $\phi\in C^1_b(H)$, then
for every $t>0$ the first order derivatives $\nabla^{C}R_{t}\left[\phi\right]$ and
$\nabla R_{t}\left[\phi\right]$ exist and are bounded, with the second one
given by
\begin{equation}\label{eq:derivateprimelisce}
\<\nabla(R_{t}\left[\phi\right]  )(x) , h\>_H=R_{t}\left[\<\nabla\phi, e^{tA}h\>_H\right](x),\
 \quad \forall\, h\in\calh.
\end{equation}
Moreover the second order derivatives $\nabla\nabla^{C}R_{t}\left[\phi\right]$,
$\nabla^{C}\nabla R_{t}\left[\phi\right]$ exist, coincide, and we have
\begin{align}
&\< \nabla \nabla^{C}(R_{t}\left[\phi\right]  )(x)  k,h\>_H=
\nonumber
\\[2mm]
&=\int_{H} \<\nabla\bar\phi\left(Pz+Pe^{tA}x\right),
Pe^{tA}h\>_H
\<(Q_{t})^{-1/2}\left( Pe^{tA} Ck\right),
(Q_t)^{-1/2}z\>_H\caln(0,Q_t)(dz).
\label{eq:derivateseconde}
\end{align}
Moreover for every $k\in K,\,h\in H$
we have the estimate
\begin{equation}\label{eq:stimaDRT2-gen}
\vert\< \nabla\nabla^C(R_{t}\left[\phi\right])(x)k,h\>_H\vert
\leq \Vert \nabla\bar\phi\Vert_\infty
\Vert Q_{t}^{-1/2}P e^{tA} C
\Vert_{\call(K;H)} \;\vert k\vert_K\,\vert h\vert_H.
\end{equation}
\end{proposition}
\dim
We first prove \myref{eq:derivateprimelisce}.
Let $\bar\phi:\operatorname{Im}P\rightarrow \R$ be such that, defining $\phi$ as in (\ref{fi-gen}), $\phi\in C^1_b(H)$.
For any
$h\in H$ we have,
applying the dominated convergence theorem,
\begin{align*}
 \<\nabla R_{t}\left[\phi\right]  (x),h\>_H
&
=\lim_{\alpha\rightarrow 0}\frac{1}{\alpha}
 \left[\int_H\phi\left(z+e^{tA}(x+\alpha h)\right)
 \caln\left(0,Q_{t}\right)  (dz)
 -\int_H
 \phi\left(z+e^{tA}x\right)  \caln\left(  0,Q_{t}\right)  (dz) \right]
 \\[2mm]
 &  =\int_H\lim_{\alpha\rightarrow 0}\frac{1}{\alpha}
 \left[\phi\left(z+e^{tA}(x+\alpha h)\right)-\phi\left(z+e^{tA}x\right)\right]
 \caln\left(0,Q_{t}\right)  (dz)\\[2mm]
 &  =\int_H\<\nabla \phi\left(z+e^{tA}x\right), e^{tA}h\>_H
 \caln\left(0,Q_{t}\right)  (dz)=R_{t}\left[\<\nabla \phi , e^{tA}h\>_H\right ]  (x).
\end{align*}
The boundedness
of $\<\nabla(R_{t}\left[\phi\right]  )(x) , h\>$ easily follows.
\noindent We compute the second order derivatives starting from
$\nabla \nabla^C R_t[\phi]$. Using \myref{eq:formulaDRT1-gen-P}
and the Dominated Convergence Theorem
we get, for $h\in H\,,k\in K$,
\begin{align*}
&\lim_{\alpha\rightarrow 0}\dfrac{1}{\alpha}
\left[\<\nabla^{C}(R_{t}\left[\phi\right]  )(x+\alpha h)  k\>_K
-\<\nabla^{C}(R_{t}\left[\phi\right]  )(x),  k\>_K
\right]=\\[2mm]
&=\lim_{\alpha\rightarrow 0}\dfrac{1}{\alpha}
 \int_{H}\left( \bar\phi\left(Pz+P e^{tA}(x+\alpha h)\right) -
 \bar\phi\left(Pz+Pe^{tA}x\right)\right)
 \<(Q_{t})^{-1/2}\left( Pe^{tA} Ck\right), (Q_t)^{-1/2}z\>_H
 \caln(0,Q_t)(dz)\\[2mm]
 &= \int_{H} \<\nabla\bar\phi\left(Pz+Pe^{tA}x\right),
 Pe^{tA}h\>_H
 \<(Q_{t})^{-1/2}\left( Pe^{tA} Ck\right), (Q_t)^{-1/2}z\>_H
 \caln(0,Q_t)(dz),
\end{align*}
Similarly, using (\ref{eq:derivateprimelisce}) and \myref{fi-gen},
we get, for $h\in H\,,k\in K$,
\begin{align*}
&\lim_{\alpha\rightarrow 0}\dfrac{1}{\alpha}
\left[\<\nabla(R_{t}\left[\phi\right]  )(x+\alpha Ck) , h\>_H
-\<\nabla(R_{t}\left[\phi\right]  )(x),  h\>_H
\right]=\\[2mm]
 &=\lim_{\alpha\rightarrow 0}\dfrac{1}{\alpha}
 \int_{H}\left(\<\nabla \bar\phi\left(Pz+Pe^{tA}(x+\alpha Ck)\right),Pe^{tA}h\>_H -
\<\nabla\bar\phi\left(Pz+Pe^{tA}x\right),Pe^{tA}h\>_H\right)
 \caln(0,Q_t)(dz)\\[2mm]
 &= \int_{H} \<\nabla\bar\phi\left(Pz+Pe^{tA}x\right),
 Pe^{tA}h\>_H
 \<(Q_{t})^{-1/2}\left( Pe^{tA}C k\right), (Q_t)^{-1/2}z\>_H
 \caln(0,Q_t)(dz).
\end{align*}
The above immediately implies \myref{eq:derivateseconde} and the estimate \myref{eq:stimaDRT2-gen}.
\qed
%

\subsection{Partial smoothing in our model}
\label{section-smoothOU-particular}
In the setting of Section \ref{section-statement}
we assume that Hypothesis \ref{ipotesibasic} holds true.
We take $\calh=\R^n \times L^2(-d,0;\R^n)$, $\Xi=\R^k$,
$(\Omega, \calf,  \P)$ a complete probability space, $W$ a standard
Wiener process in $\Xi$, $A$ and $G$ as in (\ref{A}) and (\ref{G}).
Then, for $x\in \calh$, we take the Ornstein-Uhlenbeck process $X^x(\cdot)$
given by \myref{ornstein-mild-gen}.
The associated Ornstein-Uhlenbeck transition semigroup $R_t$ is defined
as in \myref{ornstein-sem-gen} for all $f\in B_b(\calh)$.


The operator $P$ of the previous subsection here is the projection
$\Pi_0$ on the first component of the space
$\calh$, similarly to Remark \ref{remark-gen1}- (ii).
Hence, given any $\bar\phi\in B_b(\R^n)$, we define, as in (\ref{fi0})
a function $\phi \in B_b(\calh)$, by setting
\begin{equation}\label{fi}
\phi(x)=\bar\phi(\Pi_0 x)=\bar\phi(x_0) \quad
\forall x=\left(
x_0 ,x_1  \right)\in \calh.
\end{equation}
For such functions, the Ornstein-Uhlenbeck semigroup $R_t$ is written as
\begin{equation}
 \label{ornstein-sem-spec}
R_t[\phi](x)=\E \phi (X^x(t))=\E \bar\phi ((X^x(t))_0)=\int_\calh\bar\phi((z+e^{tA}x)_0)\caln(0,Q_t)(dz).
\end{equation}
Concerning the covariance operator $Q_t$ we have the following.
\begin{lemma} \label{lemmaQ_t}
Let $A$ and $G$ be defined respectively by (\ref{A}) and by (\ref{G})
and let $t \ge 0$.
Let $Q^0_t$ be the selfadjoint operator in $\R^n$ defined as
\begin{equation}\label{eq:Q^0def}
Q^0_t:=\int_0^t e^{sa_0}\sigma\sigma^*
e^{sa_0^*}\, ds.
\end{equation}
Then for every $(x_0,x_1) \in \calh$ we have
\begin{equation}\label{eq:Q^0Q}
Q_t\left(
x_0, x_1
\right) =\left(
Q^0_t x_0, 0
 \right)
\end{equation}
and so
$$
\operatorname{Im} Q_t
=\operatorname{Im} Q_t^0\times \left\lbrace 0\right\rbrace\subseteq \R^n
\times \left\lbrace 0\right\rbrace$$
Hence, for every $\bar\phi\in B_b(\R^n)$ and for the corresponding
$\phi :\calh\rightarrow\R$ defined in (\ref{fi}) we have
\begin{equation}
 \label{eq:ornstein-sem-specbis}
R_t[\phi](x)=\int_{\R^n}\bar\phi \left(z_0+(e^{tA}x)_0 \right)\caln(0,Q_t^0)(dz_0).
\end{equation}
\end{lemma}
\dim Let $\left(x_0,x_1\right)\in \calh$ and $t\ge 0$.
By direct computation we have
\begin{align*}
Q_t\left( \begin{array}{l} x_0 \\x_1
\end{array}
 \right) &=\int_0^t e^{sA}GG^*e^{sA^*}\left( \begin{array}{l}
x_0 \\x_1
\end{array}
  \right) \, ds\\[1.5mm]
 &=\int_0^te^{sA}\left(\begin{array}{ll}
  \sigma\sigma^*
 &0
 \\0&0\end{array}\right)
 e^{sA^*}\left( \begin{array}{l}
 x_0 \\x_1
 \end{array}
  \right) \, ds\\[1.5mm]
 &=\int_0^te^{sA}\left(\begin{array}{l}
  \sigma\sigma^*
 e^{sa_0^*}x_0\\ 0
 \end{array}
  \right) \, ds=\int_0^t\left(\begin{array}{l}e^{s a_0 }\sigma\sigma^*e^{sa_0^*}x_0\\
 0  \end{array}
\right) \, ds
\end{align*}
from which the first claim (\ref{eq:Q^0Q}) follows.
The second claim (\ref{eq:ornstein-sem-specbis}) is immediate.
\hfill\qed
\begin{remark}\label{remarkinfdim1}
The statement of the above lemma holds true (substituting $\R^n$ with the Hilbert space $K_0$
introduced below) also in the following more general setting.
Let $\calh=K_0\times K_1$ where $K_0$ and $K_1$ are both real separable Hilbert spaces.
Let $\Xi$ be another separable Hilbert space (the noise space) and consider
the Ornstein-Uhlenbeck process
\begin{equation}
X(t)  =e^{tA}x +\int_0^te^{(t-s)A}GdW_s
,\text{ \ \ \ }t\ge 0, \\
  \label{ornstein-mildbis}
\end{equation}
where $A$ generates a strongly continuous semigroup on $\calh$ and $G\in \call(\Xi,\calh)$.
Assuming that
\begin{itemize}
  \item $G=(\sigma,0)$ for $\sigma \in \call(\Xi,K_0)$ so
$GG^*=\left(\begin{array}{ll}
  \sigma\sigma^*
 &0
 \\0&0\end{array}\right)$ with $\sigma\sigma^* \in \call(K_0)$;
  \item for every $k_0\in K_0$, $t\ge 0$,
\begin{equation}\label{eq:invarianceeta}
e^{tA}(k_0,0)=  \left(e^{tA_0}k_0, 0\right),
    \end{equation}
where $A_0$ generates a strongly continuous semigroup in $K_0$;
\end{itemize}
then the claim still hold.
Indeed in such case we have, for $t\ge 0$, $k_0\in K_0$, $k_1 \in K_1$,
\begin{equation}\label{eq:invarianceetastar}
  \left(e^{tA^*}(k_0,k_1)\right)_0=e^{tA_0^*}k_0,
    \end{equation}
where $A_0^*$ is the adjoint of $A_0$.\footnote{Indeed once we know that
$e^{tA}(k_0,0)_1=0$ then (\ref{eq:invarianceeta}) is equivalent to ask
(\ref{eq:invarianceetastar}).} So, for $t\ge 0$,
$$
Q_t^0k_0=\int_0^t e^{sA_0} \sigma\sigma^*e^{sA_0^*}k_0  ds
$$
and
$$
 Q_t(k_0,k_1)=(Q^0_t k_0 ,0).
$$
This works, in particular, in the case described in Remark \ref{rm:diminf}.
\hfill\qedo
\end{remark}

We now analyze when Theorem \ref{lemma-der-gen} can be applied
in the cases $C=I$ or $C=B$ concentrating on the cases when
the singularity at $t=0^+$ of $\|Q_t^{-1/2}\Pi_0e^{tA}C\|$ is integrable,
as this is needed to solve the HJB equation (\ref{HJBINTRO}).

\subsubsection{$C=I$}

By Theorem \ref{lemma-der-gen} we have our partial smoothing
(namely (\ref{eq:formulaDRT1-gen}) and (\ref{eq:stimaDRT1-gen}))
for $C=I$ if
\[
\operatorname{Im}\Pi_0 e^{tA}\subseteq \operatorname{Im}Q^{1/2}_t.
\]
By Lemma \ref{lemmaQ_t} and \myref{semigroup} this implies
\begin{equation}\label{eq:inclusionsmoothingI}
\operatorname{Im}e^{ta_0}\subseteq \operatorname{Im}(Q^0_t)^{1/2}.
\end{equation}
Since, clearly, $e^{ta_0}$ is invertible and
$\operatorname{Im}(Q^0_t)^{1/2}=\operatorname{Im}Q^0_t $,
then \myref{eq:inclusionsmoothingI} is true if and only if the operator
$Q^0_t$ is invertible. On this we have the following result, taken from
\cite{Z}[Theorem 1.2, p.17] and \cite{Seidman}.

%

\begin{lemma}\label{lm:Q0inv}
The operator $Q_t^0$ defined in (\ref{eq:Q^0def})
is invertible for all $t>0$ if and only if
$$
\operatorname{Im} (\sigma,a_0\sigma, \dots , a_0^{n-1}\sigma)=  \R^n.
$$
This happens if and only if
the linear control system identified by the couple $(a_0,\sigma)$
is null controllable.
In this case, for $t \to 0^+$,
$$
\|(Q_t^0)^{-{1/2}}\| \sim   t^{-r-{1/2}}$$
where $r$ is the Kalman exponent, i.e. the minimum $r$ such that
$$
\operatorname{Im} (\sigma,a_0\sigma, \dots , a_0^{r}\sigma)=  \R^n.
$$
Hence $r=0$ if and only if $\sigma$ is onto.
\end{lemma}

We now pass to the smoothing property.

\begin{proposition}\label{cor-der}
Let $A$ and $G$ be defined respectively by (\ref{A}) and (\ref{G}).
Let $\bar\phi:\R^n\rightarrow\R$ be measurable and bounded and define,
as in (\ref{fi}), $\phi :\calh\rightarrow\R$, by setting $\phi(x)=\bar\phi(x_0)$
for every $x=(x_0,x_1)\in \calh$.
Then, if $Q^0_t$ is invertible, we have the following:
\begin{itemize}
  \item [(i)] the function $(t,x) \mapsto R_t[\phi](x)$
  belongs to $C_b((0,+\infty)\times \calh)$. Moreover it is Lipschitz continuous in $x$
  uniformly in $t\in [t_0,t_1]$ for all $0<t_0<t_1<+\infty$.
\item[(ii)] Fix any $t>0$.
$R_t[\phi]$ is Fr\'echet differentiable and we have,
for every $h\in \calh$,
\begin{align}
\langle\nabla(R_{t}\left[\phi\right]  )(x)  ,h \rangle_\calh =
\int_{\R^n}
\bar\phi\left(z_0+(e^{ tA}x)_0\right)
\<(Q^0_{t})^{-1/2}
\left( e^{tA} h\right)_0, (Q^0_t)^{-1/2}z_0\>_{\R^n}\caln(0,Q^0_t)(dz_0).
\label{eq:formulaDRT1}
\end{align}
where $\left( e^{tA} x\right)_0$, $\left( e^{tA} h\right)_0$
are given by (\ref{semigroup}).
Moreover for every $h\in \calh$
we have the estimate
\begin{equation*}
\vert \<\nabla(R_{t}\left[\phi\right]  )(x),  h\>_\calh\vert
\leq \Vert \bar\phi\Vert_\infty
\left\Vert (Q^0_{t})^{-1/2}
\left( e^{tA} \right)_0
\right\Vert_{\call(\calh;\R^n)} \;\vert h\vert_\calh.
 \end{equation*}
Hence for all $T>0$ there exists $C_T$ such that
\begin{equation}\label{stimadertutta}
 \vert  \<\nabla R_{t}\left[\phi\right]  (x),h\>_\calh \vert
 \leq C_T t^{-r -{1/2}} \Vert \bar\phi \Vert_\infty\vert h\vert_\calh,
 \quad t\in [0,T],
 \end{equation}
 where $r$ is the Kalman exponent which is $0$ if
 and only if $\sigma$ is onto.
\item [(iii)] Fix any $t>0$.
$R_t[\phi]$ is $B$-Fr\'echet differentiable and we have, for every $k\in \R^m$,
\begin{align}
&\<\nabla^B(R_{t}\left[\phi\right]  )(x),  k\>_{\R^m}
=\int_{\R^n}
\bar\phi\left(z_0+\left( e^{tA} x\right)_0\right)
\<(Q^0_{t})^{-1/2}
\left( e^{tA} Bk\right)_0, (Q^0_t)^{-1/2}z_0\>_{\R^n}\caln(0,Q^0_t)(dz_0).
\label{eq:formulaDRT1B}
\end{align}
Moreover, for every $k\in \R^m$,
\begin{equation}\label{eq:formulastimaDRT1B}
 \vert \<\nabla^{B}(R_{t}\left[\phi\right]  )(x),k\>_{\R^m}\vert
\leq \Vert \bar\phi\Vert_\infty
\left\Vert (Q^0_{t})^{-1/2}
\left( e^{tA} B\right)_0
\right\Vert_{\call(R^m,\R^n)} \; \vert k\vert_{\R^m}.
\end{equation}
Hence for all $T>0$ there exists $C_T$ such that
 \begin{equation}\label{stimader}
 \vert  \<\nabla^{B}R_{t}\left[\phi\right]  (x),k\>_{\R^m} \vert
 \leq C_T t^{-r-{1/2}} \Vert \bar\phi \Vert_\infty\vert k\vert_{\R^m},
 \quad t\in [0,T],
 \end{equation}
 where $r$ is the Kalman exponent which is $0$ if
 and only if $\sigma$ is onto.
\end{itemize}
\end{proposition}




\dim Point (ii) immediately follows from the
invertibility of $Q_t^0$, the discussion just before Lemma \ref{lm:Q0inv},
and Theorem \ref{lemma-der-gen}.
Point (i) follows from point (ii) and form the continuity of trajectories
of the Ornstein Uhlenbeck process \myref{ornstein-gen} with $A$ and $G$ given by
(\ref{A}) and (\ref{G}).
Point (iii) follows observing that the operator $\Pi_0 e^{tA}B:\R^m\to\R^n$, given in
\myref{eq:etAB} is well defined and hence, thanks to the
invertibility of $Q_t^0$, Theorem \ref{lemma-der-gen} can be applied.
\hfill\qedo


%

\subsubsection{$C=B$}

By Theorem \ref{lemma-der-gen} we have the partial smoothing
(\ref{eq:formulaDRT1-gen} and (\ref{eq:stimaDRT1-gen})
for $C=B$ if
\begin{equation}\label{eq:inclusionsmoothingB}
\operatorname{Im}\Pi_0e^{tA}B\subset \operatorname{Im}Q^{1/2}_t
=\operatorname{Im}Q^0_t
\end{equation}
Since, as proved e.g. in \cite{Z} (Lemma 1.1, p. 18 and formula (2.11) p. 210),
$$
\operatorname{Im}(Q^0_t)^{1/2}=\operatorname{Im}(\sigma,a_0\sigma, \dots a_0^{n-1}\sigma ),
$$
then, using \myref{eq:etAB},
\myref{eq:inclusionsmoothingB} is verified if and only if
\begin{equation}\label{eq:inclusionsmoothingBbis}
\operatorname{Im}\left(e^{ta_0}b_0 +\int_{-d}^0 1_{[-t,0]}e^{(t+r)a_0}b_1(dr)
\right)\subseteq \operatorname{Im}(\sigma,a_0\sigma, \dots a_0^{n-1}\sigma ).
\end{equation}
We now provide conditions, possibly weaker than the invertibility
of $Q_t$, under which \myref{eq:inclusionsmoothingBbis} is verified and the singularity at $t=0^+$ of $\|Q_t^{-1/2}\Pi_0 e^{tA}B\|$
is integrable.
We first recall the following result (see \cite{Z}, Proposition 2.1, p. 211).
\begin{proposition}\label{propimzab}
If $F_1$ and $F_2$ are linear bounded operators acting between
separable Hilbert spaces $X$, $Z$ and $Y$, $Z$ such that
$\|F_1^* f \| = \| F_2^* f \|$  for any
$f \in Z^*$, then $Im F_1 = Im F_2$ and $\|F_1^{-1}z\| = \|F_2^{-1}z\|$
for all  $z\in Im F_1$.
\end{proposition}


\begin{proposition}\label{lemmaderhpdeb}
Assume that Hypothesis \ref{ipotesibasic} holds. Assume moreover that, either
\begin{equation}\label{eq:hpdebreg}
\operatorname{Im}(e^{ta_0}b_0)\subseteq
\operatorname{Im}\sigma, \; \forall t > 0;
\qquad b_1 \in L^2(-d,0,\call(\R^m;\R^n)), \; \operatorname{Im}b_1(s)\in\operatorname{Im}\sigma,
\quad  a.e.\, \forall s\in[-d,0]
\end{equation}
or
\begin{equation}\label{eq:hpdebregbis}
\operatorname{Im}\left(e^{ta_0}b_0 +\int_{-d}^0 1_{[-t,0]}e^{(t+r)a_0}b_1(dr)
\right)
\subseteq\operatorname{Im}\sigma,
\quad \forall t>0.
\end{equation}
Then, for any bounded measurable $\phi$ as in (\ref{fi}),
$R_{t}\left[\phi\right]$ is $B$-Fr\'echet differentiable for every
$t>0$, and, for every $h\in \R^m$,
$\<\nabla^{B}(R_{t}\left[\phi\right])(x),k\>_{\R^m}$ is given by \myref{eq:formulaDRT1B}
and satisfies the estimate \myref{eq:formulastimaDRT1B}.
Moreover for all $T>0$ there exists $C_T$ such that
\begin{equation}
\vert\<  \nabla^{B}(R_{t}\left[\phi\right])  (x),k\>_{\R^m}\vert
\leq C_T t^{-{1/2}} \Vert \bar\phi \Vert_\infty \;\vert k\vert_{\R^m}.
\label{eq:stimaderBnew}
\end{equation}
\end{proposition}
\dim
Consider the following linear deterministic controlled system in $\calh$:
\begin{equation}\label{eq:detcontrsyst}
\left\{
\begin{array}{l}
 dX(t)=AX(t)dt+Gu_1(t)dt \\
X(0)=Bh,
\end{array}
\right.
\end{equation}
where the state space is $\calh$, the control space is $U_1=\R^k$, the control strategy is
$u_1\in L^2_{loc}([0,+\infty);U_1)$, the initial point is $Bh$ with $h\in\R^m$.
Define the linear operator
$$
{\mathcal L}^0_t: L^2([0,t];U_1)\to \R^n, \quad  u_1(\cdot)\mapsto
\int_0^te^{a_0 (t-s)}\sigma u_1(s)ds.
$$
Then the first component of the state trajectory is
\begin{equation}\label{sistcontrfin}
X^0(t)=\Pi_0 e^{tA}Bk + {\mathcal L}^0_t u_1
\end{equation}
Hence $X^0$ can be driven to $0$ in time $t$ if and only if
$$
\Pi_0 e^{tA}Bk \in \operatorname{Im}\mathcal L^0_t
$$
In such case, by the definition of pseudoinverse
(see Subsection \ref{subsection-notation}),
we have that the control which brings $X^0$ to $0$ in time $t$ with minimal $L^2$ norm is $(\mathcal L^0_t)^{-1}\Pi_0 e^{tA}Bk$
and the corresponding minimal square norm is
\begin{equation}
\mathcal{E}\left(  t,Bk\right):  =\min\left\{ \int_{0}
^{t}\left\vert u_{1}(s)\right\vert ^{2}ds:X\left(  0\right)
=B k,\text{ }X^0\left(  t\right)  =0\right\}
=\|(\mathcal L^0_t)^{-1}\Pi_0 e^{tA}Bk\|^2_{L^2(0,t;U_1)}. \label{energyB}
\end{equation}
Since for all $z \in \R^n$ we have
$$
\|(Q_t^0)^{1/2}z\|^2_{\R^n}= |\<Q_t^0z,z\>_{\R^n}|=
\|(\mathcal L^0_t)^* z\|^2_{L^2(0,t;U_1)}
$$
then, by Proposition \ref{propimzab}, we get
$$
\operatorname{Im}\left((Q^0_{t})^{{1/2}}\right)=\operatorname{Im}\mathcal L^0_t .
$$
and
$$
\|(\mathcal L^0_t)^{-1}\Pi_0 e^{tA}Bk\|_{L^2(0,t;U_1)}
=\|(Q^0_{t})^{-{1/2}}\Pi_0 e^{tA}Bk\|_{\R^n}.
$$
Hence, by \myref{energyB}, to estimate
$\|(Q^0_{t})^{-{1/2}}\Pi_0 e^{tA}Bk\|_{\R^n}$
it is enough to estimate the minimal energy to steer $X^0$ to $0$ in time $t$.
When (\ref{eq:hpdebreg}) holds we see, by simple computations, that the control
\begin{equation}\label{control}
\bar u_1(s)=-\dfrac{1}{t}\sigma^{-1}e^{sa_0}b_0 k-\sigma^{-1}b_1(-s)k 1_{[-d,0]}(-t),
\quad s \in [0,t],
\end{equation}
where $\sigma^ {-1}$ is the pseudoinverse of $\sigma$,
brings $X^0$ to $0$ in time $t$.
Hence, for a suitable $C>0$ we get
\[
\mathcal{E}\left(  t,Bk\right)\leq
\int_0^t \bar u_1^2(s) ds
\leq C\left(\dfrac{1}{t}+
\Vert b_1 \Vert^2_{L^2 \left([-d,0];\call(\R^m;\R^n)\right)}\right)|k|^2_{\R^m}.
\]
So, for a, possibly different constant $C$, we get $\Vert (Q^0_{t})^{-{1/2}}
\left( e^{tA} Bk\right)_0
\Vert_{\call(\R^m;\R^n)}\leq C t^{-\frac{1}{2}}|k|_{\R^n}$ and the estimate is proved.
If we assume (\ref{eq:hpdebregbis}) we can take as a control,
on the line of \cite{Z}, Theorem 2.3-(iii), p.210,
\begin{equation}\label{controlbis}
\hat u_1(s)=
-\sigma^{*}e^{(t-s)a_0^*}(Q_t^0)^{-1}\left(e^{ta_0}b_0 k+
\int_{-d}^0 1_{[-t,0]}e^{(t+r)a_0}b_1(dr)k\right),
\quad s \in [0,t],
\end{equation}
and use that the singularity of the second term as $t\to 0^+$ is still of order $\dfrac{1}{t}$ since (\ref{eq:hpdebregbis}) holds (see e.g. \cite{SY}, Theorem 1).
Once this estimate is proved, the proof of the $B$-Fr\'echet differentiability
is the same as the one of Proposition \ref{cor-der}-(iii).
%
\qed

\begin{remark}\label{remark-derinfdim}
The above results can be generalized to the case, introduced in Remark \ref{rm:diminf} above, when the first component of the space $\calh$ is infinite dimensional.
\begin{itemize}
\item For the case $C=I$ the required partial smoothing
holds if we ask, in place of the invertibility of $Q_t^0$,
that, for every $t> 0$,
\begin{equation}\label{eq:hpdiff}
  \operatorname{Im}\left( e^{tA} \right)_0 \subseteq   \operatorname{Im}(Q_t^0)^{{1/2}}
\end{equation}
which would imply that the linear operator $(Q_t^0)^{-{1/2}}\left( e^{tA} \right)_0$
is continuous from $K_1$ into itself.
\item For the case $C=B$, the required partial smoothing
holds if we ask that, for every $t> 0$,
\begin{equation}\label{eq:hpBdiff}
  \operatorname{Im}\left( e^{tA} B \right)_0 \subseteq   \operatorname{Im}(Q_t^0)^{{1/2}}
\end{equation}
which would imply that the operator $(Q_t^0)^{-{1/2}}\left( e^{tA} B\right)_0$ is continuous from $U$ to $K_1$.
\end{itemize}
Clearly, in this generalized setting the estimates (\ref{stimadertutta}) and
(\ref{stimader}) do not hold any more and they depend
on the specific operators $A$, $B$, $\sigma$.
%
\hfill\qedo
\end{remark}

\section{Smoothing properties of the convolution}
\label{section-smooth-conv}
By Proposition
\ref{lemmaderhpdeb}
we know that if \myref{eq:hpdebreg} or \myref{eq:hpdebregbis} hold and $\phi$ is as in (\ref{fi}) with $\bar\phi$ measurable and bounded then $\nabla^{B}(R_t \left[\phi\right]  )(x)$ exists and its norm blows up like
$t^{-{1/2}}$ at $0^+$. Moreover if
$\bar\phi\in C_b(\R^n)$, then $R_t[\phi]\in C_b([0,T]\times\calh)$, see e.g. \cite{PriolaTesi} Proposition 6.5.1 (or the discussion at the end of \cite{PriolaStudia}).

We now prove that, given $T>0$,
for any element $f$ of a suitable family of functions
in $C_b([0,T]\times \calh)$,
a similar smoothing property for the convolution integral
$\int_{0}^{t}R_{t-s} [f(s,\cdot)](x) ds$ holds.
This will be a crucial step to prove the existence and uniqueness
of the solution of our HJB equation in next section.

For given $\alpha \in (0,1)$ we define now a space designed for our purposes.

\begin{definition}\label{df:Sigma}
Let $T>0$, $\alpha \in (0,1)$.
A function $g\in C_b([0,T]\times \calh)$
belongs to $\Sigma^1_{T,\alpha}$ if there exists a function
$f\in C^{0,1}_{\alpha}([0,T]\times \R^n)$ such that
$$g(t,x)=f\left(t,(e^{tA}x)_0\right),
\qquad \forall (t,x) \in [0,T]\times \calh.
$$
\end{definition}
If $g\in \Sigma^1_{T,\alpha}$, for any $t\in(0,T]$ the function $g(t,\cdot)$ is
both Fr\'echet differentiable and $B$-Fr\'echet differentiable.
Moreover, for $(t,x)\in [0,T]\times \calh$, $h \in \calh$, $k\in \R^m$,
$$
\<\nabla g(t,x),h\>_{\calh}=\<\nabla f\left(t,(e^{tA}x)_0\right),(e^{tA}h)_0\>_{\R^n},
\quad and \quad
\<\nabla^B g(t,x),k\>_{\R^m}=\<\nabla f\left(t,(e^{tA}x)_0\right),(e^{tA}Bk)_0\>_{\R^n}.
$$
This in particular imply that, for all $k\in \R^m$
\begin{equation}
\label{eq:nablaperSigma}
\<\nabla^B g(t,x),k\>_{\R^m}=\<\nabla g(t,x),Bk\>_{
\<\R^n \times C([-d,0];\R^n),\R^n \times C^*([-d,0];\R^n)\>},
\end{equation}
which also means $B^*\nabla g = \nabla^B g$.
For later notational use we call $\bar f\in C_b((0,T]\times \R^n;\R^m)$
the function defined by
$$
\<\bar f(t,y),k\>_{\R^m}=
t^\alpha\<\nabla f\left(t,y\right),(e^{tA}Bk)_0\>_{\R^n},
\qquad (t,y)\in (0,T]\times \R^n, \quad k \in \R^m,
$$
which is such that
$$
t^\alpha\nabla^B g(t,x)=\bar f\left(t,(e^{tA}x)_0\right).
$$
We also notice that if $g\in \Sigma^1_{T,\alpha}$, then in order to have $g$ $B$-Fr\'echet differentiable it suffices to require $(e^{tA}B)_0$
bounded and continuous.

When (\ref{eq:hpdebreg}) or
(\ref{eq:hpdebregbis}) hold we know, by Proposition \ref{lemmaderhpdeb},
that the function $g(t,x)=R_t[\phi](x)$
for $\phi$ given by (\ref{fi}) with $\bar\phi$ bounded and continuous,
belongs to $\Sigma^1_{T,1/2}$.

\begin{lemma}\label{lemma:Sigma}
The set $\Sigma^1_{T,\alpha}$ is a closed subspace of $C^{0,1,B}_{\alpha}([0,T]\times \calh)$.
\end{lemma}
\dim
It is clear that $\Sigma^1_{T,\alpha}$ is a vector subspace of
$C^{0,1,B}_{\alpha}([0,T]\times \calh)$. We prove now that it is closed.
Take any sequence $g_n \to g$ in $C^{0,1,B}_{\alpha}([0,T]\times \calh)$.
Then to every $g_n$ we associate the corresponding $f_n$ and $\bar f_n$.
The sequence $\{f_n\}$ is a Cauchy sequence in
$C_b([0,T]\times \R^n)$. Indeed for any $\epsilon >0$ take
$(t_\epsilon,y_\epsilon)$ such that
$$
\sup_{(t,y) \in [0,T]\times \R^n}|f_n(t,y)-f_m(t,y)|
<
\epsilon + |f_n(t_\epsilon,y_\epsilon)-f_m(t_\epsilon,y_\epsilon)|
$$
Then choose $x_\epsilon \in \calh $ such that
$y_\epsilon=(e^{t_\epsilon A}x_\epsilon)_0$
(this can always be done choosing e.g.
$x_\epsilon=(e^{-t_\epsilon a_0}y_\epsilon,0)$).
Hence we get
$$
\sup_{(t,y) \in [0,T]\times \R^n}|f_n(t,y)-f_m(t,y)|<
\epsilon + |g_n(t_\epsilon,x_\epsilon)-g_m(t_\epsilon,x_\epsilon)|
\le
\epsilon + \sup_{(t,x) \in [0,T]\times \calh}
|g_n(t,x)-g_m(t,x)|.
$$
Since $\{g_n\}$ is Cauchy, then $\{f_n\}$ is Cauchy, too. So there exists a function
$f \in C_b([0,T]\times \R^n)$ such that $f_n \to f$ in $C_b([0,T]\times \R^n)$.
This implies that $g(t,x)=f(t,(e^{tA}x)_0)$ on $[0,T]\times \calh$.
With the same argument we get that there exists a function
$\bar f \in C_b((0,T]\times \R^n;\R^m)$ such that $\bar f_n \to \bar f$
in $C_b((0,T]\times \R^n;\R^m)$.
This implies that $t^\alpha \nabla^B g(t,x)=\bar f(t,(e^{tA}x)_0)$
on $(0,T]\times \calh$.
\qed

Next, in analogy to what we have done defining $\Sigma^1_{T,\alpha}$,
we introduce a subspace $\Sigma^2_{T,\alpha}$ of functions
$g \in C^{0,2,B}_{\alpha}([0,T]\times \calh)$
that depends in a special way on the variable $x\in\calh$.
\begin{definition}\label{df:Sigma2}
A function $g\in C_b([0,T]\times \calh)$ belongs to $\Sigma^2_{T,\alpha}$
if there exists a function
$f\in C^{0,2}_{\alpha}([0,T]\times \R^n)$ such that
for all $(t,x) \in [0,T]\times \calh$,
\begin{align*}
&g(t,x)=f\left(t,(e^{tA}x)_0\right).
\end{align*}
\end{definition}
If $g\in \Sigma^2_{T,\alpha}$ then for any $t\in(0,T]$ the function $g(t,\cdot)$
is Fr\'echet differentiable and
$$
\<\nabla g(t,x),h\>_\calh=
\<\nabla f\left(t,(e^{tA}x)_0\right),(e^{tA}h)_0\>_{\R^n},
\quad \hbox{for $(t,x)\in [0,T]\times \calh$, $h \in \calh$.}
$$
Moreover also $\nabla g(t,\cdot)$ is $B$-Fr\'echet differentiable and
$$
\<\nabla^B\left(\nabla g(t,x)h\right)),k\>_{\R^m}=
\<\nabla^2 f\left(t,(e^{tA}x)_0\right)(e^{tA}h)_0,(e^{tA}Bk)_0 \>_{\R^n},
\quad \hbox{for $(t,x)\in [0,T]\times \calh$, $h \in \calh$, $k \in \R^m$.}
$$
We also notice that, since the function $f$ is twice continuously
Fr\'echet differentiable the second order derivatives
$\nabla^B\nabla g$ and $\nabla\nabla^B g$ both exist and coincide:
$$
 \<\nabla^B\<\nabla g(t,x),h\>_{\calh},k\>_{\R^m}
 =\<\nabla\<\nabla^B g(t,x),k\>_{\R^m},h\>_{\calh}.
$$
Again for later notational use we call
$\bar f_1\in C_b([0,T]\times \R^n;\R^m)$
the function defined by
$$
\<\bar f_1(t,y),h\>_{\R^m}=
\<\nabla f\left(t,y\right),(e^{tA}Bh)_0\>_{\R^n},
\qquad (t,y)\in [0,T]\times \R^n, \quad h \in R^m,
$$
which is such that
$$
\nabla^B g(t,x)=\bar f_1\left(t,(e^{tA}x)_0\right).
$$
Similarly we call
$\bar{\bar f}\in C_b\left((0,T]\times \R^n;\call(\calh,\R^m)\right)$
the function defined by
$$
\<\<\bar{\bar f}(t,y),h\>_\calh,k\>_{\R^m}=
t^\alpha\<\nabla^2 f\left(t,y\right)(e^{tA}h)_0,(e^{tA}Bk)_0 \>_{\R^n}
\qquad (t,y)\in [0,T]\times \R^n, \quad h \in \calh, \; k\in R^m,
$$
which is such that
$$
t^\alpha\nabla^B\nabla g(t,x)=t^\alpha\nabla\nabla^B g(t,x)=
\bar{\bar f}\left(t,(e^{tA}x)_0\right).
$$

We now pass to the announced smoothing result.

\begin{lemma}\label{lemma reg-convpercontr}
Let
(\ref{eq:hpdebreg}) or (\ref{eq:hpdebregbis}) hold true.
Let $T>0$ and let $\psi:\R^m\rightarrow\R$ be a continuous
function satisfiying Hypothesis (\ref{ipotesicostoconcreto}), estimates
(\ref{eq:Hlip}).
Then
\begin{itemize}
 \item [i)]
for every $g \in \Sigma^1_{T,{1/2}}$, the function
$\hat g:[0,T]\times \calh\rightarrow \R$ belongs to $\Sigma^1_{T,{1/2}}$ where
\begin{equation}\label{iterata-primag}
\hat g(t,x) =\int_{0}^{t}
R_{t-s} [\psi(\nabla^{B}g(s,\cdot))](x)  ds.
\end{equation}
 Hence, in particular, $\hat g(t,\cdot)$
is $B$-Fr\'echet differentiable for every $t\in (0,T]$ and, for all $x\in \calh$,
\begin{equation}\label{stimaiterata-primag}
\left\vert \nabla ^B(\hat g(t,\cdot))(x) \right\vert_{(\R^m)^*}    \leq
    C\left(t^{{1/2}}+\Vert g \Vert_{C^{0,1,B}_{{1/2}}}\right).
 \end{equation}
 If $\sigma$ is onto, then $\hat g(t,\cdot)$
is Fr\'echet differentiable for every $t\in (0,T]$ and, for all $h\in\calh$, $x\in \calh$,
\begin{equation}\label{stimaiterata-primagbis}
\left\vert\nabla(\hat g(t,\cdot))(x) \right\vert_{\calh^*}
   \leq
    C\left(t^{{1/2}}+\Vert g \Vert_{C^{0,1}_{{1/2}}}\right).
\end{equation}
\item [ii)] Assume moreover that $\psi\in C^1(\R^m)$.
For every $g \in \Sigma^2_{T,{1/2}}$, the function
$\hat g$ defined in (\ref{iterata-primag})
belongs to $\Sigma^2_{T,{1/2}}$. Hence, in particular, the second order derivatives
$\nabla\nabla^B\hat g(t,\cdot)$ and $\nabla^B\nabla\hat g(t,\cdot)$ exist, coincide and
for every $t\in (0,T]$ and, for all $x\in \calh$,
\begin{equation}\label{stimaiterata-secondag}
\left\vert \nabla ^B\nabla(\hat g(t,\cdot))(x) \right\vert_{\calh^*\times(\R^m)^*}    \leq
    C
\Vert g \Vert_{C^{0,2,B}_{{1/2}}}
 \end{equation}
 If $\sigma$ is onto, then $\hat g(t,\cdot)$
is twice Fr\'echet differentiable and for every $t\in(0,T]$, for all $h\in\calh$ and $x\in \calh$,
\begin{equation}\label{stimaiterata-secondagbis}
\left\vert\nabla^2(\hat g(t,\cdot))(x) \right\vert_{\calh^*\times \calh^*}
   \leq
    C
\Vert g \Vert_{C^{0,2}_{{1/2}}}
\end{equation}
\end{itemize}

\end{lemma}
\dim
We start by proving that (\ref{iterata-primag}) is $B$-Fr\'echet differentiable
and we exhibit its $B$-Fr\'echet derivative.
Recalling (\ref{ornstein-sem-gen}) we have
\begin{align*}
  \int_{0}^{t}
R_{t-s} \left[\psi\left(\nabla^{B}(g(s,\cdot))\right)\right](x)  ds
&=\int_{0}^{t}
\int_{\calh} \psi\left(\nabla^{B}(g(s,\cdot))
\left(z+e^{(t-s)A}x \right)\right)\caln(0,Q_{t-s})(dz)
\end{align*}
By the definition of $\Sigma^1_{T,{1/2}}$, we see that
\begin{align}\label{eq:nablaBshiftg}
 &s^{1/2} \nabla^{B}g(s,z+e^{(t-s)A}x)
 =
 \bar f\left(s, (e^{sA}z)_0+(e^{tA}x)_0\right)
 \qquad \forall t\ge s>0, \; \forall x,z \in \calh.
\end{align}
Hence the function $\hat f$ associated to $\hat g$ is
\begin{align*}
\hat f (t,y)&=\int_{0}^{t}
\int_{\calh} \psi\left(s^{-{1/2}}
 \bar f\left(s, (e^{sA}z)_0+y\right)
\right)\caln(0,Q_{t-s})(dz)
\end{align*}
with, by our assumptions on $\psi$,
$$
\Vert\hat f \Vert_\infty \le C\int_{0}^{t}\left( 1 +
s^{-{1/2}} \Vert \bar f \Vert_\infty \right)ds
$$
To compute the $B$-directional derivative we look at the limit
\begin{align*}
 &\lim_{\alpha\rightarrow 0}\dfrac{1}{\alpha}
 \left[\int_{0}^{t}
 R_{t-s}\left[\psi\left( \nabla^{B}(g(s,\cdot))\right)\right](x+\alpha Bk)ds -
\int_{0}^{t}
 R_{t-s}\left[\psi\left( \nabla^{B}(g(s,\cdot))\right)\right](x)  ds\right].
\end{align*}
From what is given above we get
\begin{align*}
&\int_{0}^{t}
 R_{t-s}\left[\psi\left( \nabla^{B}(g(s,\cdot))\right)\right](x+\alpha Bk)ds \\
&=\int_{0}^{t} \int_\calh  \psi\left( s^{-{1/2}}
\bar f\left(s, (e^{sA}z)_0+(e^{tA}(x+\alpha Bh))_0\right)
\right)\caln(0,Q_{t-s})(dz) ds.
\end{align*}
Since the last integrand only depends on $(e^{sA}z)_0+(e^{ tA}(x+\alpha Bk))_0$,
we have, arguing exactly as in the proof of Theorem \ref{lemma-der-gen},
\begin{align*}
&\int_{0}^{t}
 R_{t-s}\left[\psi\left( \nabla^{B}(g(s,\cdot))\right)\right](x+\alpha Bk)ds =
\\[2mm]
&=\int_{0}^{t} \int_\calh
\psi\left(s^{-{1/2}}
\bar f\left(s, (e^{sA}z)_0+(e^{tA}x)_0\right)
\right)
\caln\left(\left(
(e^{tA}\alpha Bk)_0, 0\right),Q_{t-s}\right)(dz)ds=
\\[2mm]
&=\int_{0}^{t} \int_\calh
\psi\left(s^{-{1/2}}
\bar f\left(s, (e^{sA}z)_0+(e^{tA}x)_0\right)
\right)
d(t,t-s,\alpha Bk,z)
\caln\left(0,Q_{t-s}\right)(dz)ds,
\end{align*}
where
\begin{align}
d( t_1,t_2,y,z)   &  =\frac{d\caln\left(
\left(\left(e^{t_1 A}y\right)_0 , 0\right),Q_{t_2}
\right)  }{d\mathcal{N}\left(  0,Q_{t_2}\right)  }(z)
\nonumber
\\[2mm]
&  =\exp\left\{  \left\langle Q_{t_2}^{-1/2}
\left(
\left(e^{t_1A}y\right)_0, 0\right),Q_{t_2}^{-{1/2}}z\right\rangle_\calh
-\frac{1}{2}\left|  Q_{t_2}^{-{1/2}}\left(
\left( e^{t_1A}y\right)_0, 0\right)
\right|_\calh  ^{2}\right\}  ,
\label{eq:density1bisg}
\end{align}
Hence
\begin{align*}
&\lim_{\alpha\rightarrow 0}\dfrac{1}{\alpha}
 \left[\int_{0}^{t}
 R_{t-s}\left[\psi\left( \nabla^{B}(g(s,\cdot))\right)\right](x+\alpha Bk)ds -
\int_{0}^{t}
 R_{t-s}\left[\psi\left( \nabla^{B}(g(s,\cdot))\right)\right](x)  ds\right]=
 \\[2mm]
&=\lim_{\alpha\rightarrow 0}\dfrac{1}{\alpha}
\int_{0}^{t} \int_\calh
\psi\left(s^{-{1/2}}
\bar f\left(s, (e^{sA}z)_0+(e^{tA}x)_0\right)
\right)
\frac{d(t,t-s,\alpha Bk,z)-1}{\alpha}
\caln\left(0,Q_{t-s}\right)(dz)ds
\\[2mm]
&=
\int_{0}^{t} \int_\calh
\psi\left(s^{-{1/2}}
\bar f\left(s, (e^{sA}z)_0+(e^{tA}x)_0\right)
\right)
\<Q_{t-s}^{-{1/2}}\left(
\left( e^{tA} Bk\right)_0,0
\right), Q_{t-s} ^{-{1/2}}z\>_\calh
\caln\left(0,Q_{t-s}\right)(dz)ds.
\end{align*}
Since the above limit is uniform for $k$ in the unit sphere, then
we get the required $B$-Fr\'echet differentiability and
\begin{align}
&\<\nabla ^B \left(\int_{0}^{t}
R_{t-s}\left[\psi\left( \nabla^{B}(g(s,\cdot)\right)\right]  ds
\right) (x),k \>_{\R^m}=
\label{eq:derBconvnew}
\\[3mm]
\nonumber
&=
\int_{0}^{t} \int_\calh
\psi\left(
s^{-{1/2}} \bar f\left(s, (e^{sA}z)_0+(e^{tA}x)_0\right)
\right)
\<(Q^0_{t-s})^{-{1/2}}
\left( e^{tA} Bk\right)_0, (Q^0_{t-s})^{-{1/2}}z_0\>_{\R^n}
\caln\left(0,Q_{t-s}\right)(dz)ds.
\end{align}
Finally we prove the estimate (\ref{stimaiterata-primag}). Using the above representation and the Holder inequality we have
\begin{align*}
&\left\vert\<\nabla ^B \left(\int_{0}^{t}
R_{t-s}\left[\psi\left( \nabla^{B}(g(s,\cdot))\right)\right]  ds
\right) (x),k \>_{\R^m}\right\vert \le
\\[3mm]
&\leq C\int_{0}^{t}
 \int_\calh  \left(1+\left\vert
s^{-{1/2}} \bar f\left(s, (e^{sA}z)_0+(e^{tA}x)_0\right)
\right\vert\right)
 \left \vert
\<(Q^0_{t-s})^{-{1/2}}
\left( e^{tA} Bk\right)_0, (Q^0_{t-s})^{-{1/2}}z_0\>_{\R^n}
 \right\vert
\caln(0,Q_{t-s})(dz) ds\\
&\leq C\int_{0}^{t}\left(1+s^{-{1/2}} \left\Vert g \right\Vert_{C^{0,1,B}_{{1/2}}}
\right)
\left\Vert (Q^0_{t-s})^{-{1/2}}(e^{tA}Bk)_0 \right\Vert_{\call(\R^m;\R^n)} ds \\
 &\leq C\int_{0}^{t}\left(1+s^{-{1/2}}\left\Vert g \right\Vert_{C^{0,1,B}_{{1/2}}}\right)
 (t-s)^{-{1/2}}\vert k\vert_{\R^m} \,ds \leq C\left(t^{{1/2}}
 +\left\Vert g \right\Vert_{C^{0,1,B}_{{1/2}}} \right)
\vert h\vert_{\R^m}.
\end{align*}
Observe that in the last step we have used
the estimate
$$
\left\Vert (Q^0_{t-s})^{-{1/2}}(e^{tA}Bh)_0 \right\Vert_{\call(\R^m;\R^n)}
\leq C({t-s})^{-{1/2}}
$$
which follows from the proof of Proposition \ref{cor-der} (or Proposition \ref{lemmaderhpdeb}).
Moreover we have also used that
\begin{equation}
\label{eq:betaeuler}
 \int_0^t(t-s)^{-{1/2}}s^{-{1/2}}ds=\int_0^1(1-x)^{-{1/2}}x^{-{1/2}}dx=\beta \left({1/2},{1/2} \right) ,
\end{equation}
where by $\beta(\cdot,\cdot)$ we mean the Euler beta function.

The Fr\'echet differentiability and the estimate
(\ref{stimaiterata-primagbis})
is proved exactly in the same way using the fact that $\sigma$ is onto
and Proposition \ref{cor-der}.

Now we consider the case of $g\in\Sigma^2_{T,{1/2}}.$
We start by proving that (\ref{iterata-primag}) is Fr\'echet differentiable
and, in order to compute the Fr\'echet derivative, we use \myref{eq:derivateprimelisce}
(which is true for every $\phi \in C^1_b(H)$, see its proof)
looking at the limit, for $h \in \calh$,
\begin{align}
&\<\nabla\hat g(t,x),h\>_\calh
=\lim_{\alpha\rightarrow 0}\dfrac{1}{\alpha}
 \left[\int_{0}^{t}
 R_{t-s}\left[\psi\left( \nabla^{B}(g(s,\cdot))\right)\right](x+\alpha h)ds -
\int_{0}^{t}
 R_{t-s}\left[\psi\left( \nabla^{B}(g(s,\cdot))\right)\right](x)  ds\right]
 \nonumber
 \\
&=\int_{0}^{t} R_{t-s}
\left[\<\nabla\left(\psi\left( \nabla^{B}(g(s,\cdot))\right)\right),
 e^{(t-s)A}h\>_\calh\right](x)ds
\nonumber
\\
&=\int_{0}^{t} R_{t-s}
\left[\<\nabla\psi\left( \nabla^{B}(g(s,\cdot))\right),
\nabla \nabla^{B}(g(s,\cdot)) e^{(t-s)A}h\>_{\R^m}\right](x)ds
\nonumber
\\
&=\int_{0}^{t} \int_\calh
\<\nabla\psi\left( \nabla^{B}(g(s,z+e^{(t-s)A}x))\right),
\nabla \nabla^{B}(g(s,z+e^{(t-s)A}x)) e^{(t-s)A}h\>_{\R^m}
\caln\left(0,Q_{t-s}\right)(dz)ds
\label{eq:perdersecvera}
\\
&=\int_{0}^{t} \int_\calh
\left[\<\nabla\psi\left(
\bar f_1\left(s, (e^{sA}z)_0+(e^{tA}x)_0\right)\right),
s^{-{1/2}}
\bar{\bar f}\left(s, (e^{sA}z)_0+(e^{tA}x)_0\right) e^{(t-s)A}h\>_{\R^m} \right]
\caln\left(0,Q_{t-s}(dz) \right)ds
\nonumber
\end{align}
Now, from calculations similar to the ones performed in the first part we arrive at
\begin{align}\label{eq:der-seconda-conv-Sigma}
\<\nabla^B\<\nabla\hat g(t,x),h\>_\calh,k\>_{\R^m}
&=\int_{0}^{t} \int_\calh
\left[\<\nabla\psi\left(
\bar f_1\left(s, (e^{sA}z)_0+(e^{tA}x)_0\right)\right),
s^{-{1/2}}
\bar{\bar f}\left(s, (e^{sA}z)_0+(e^{tA}x)_0\right) e^{(t-s)A}h\>_{\R^m}
\right.
\nonumber
\\[2mm]
&
\qquad\qquad\left.
\<(Q^0_{t-s})^{-{1/2}}
\left( e^{tA} Bk\right)_0, (Q^0_{t-s})^{-{1/2}}z_0\>_{\R^n}\right]ds.
\end{align}
Since $\nabla \psi$ is bounded (as it satisfies \myref{eq:Hlip})
then (\ref{stimaiterata-secondag}) easily follows
by the definition of $\bar{\bar f}$ and from \myref{eq:betaeuler}.
Second order differentiability and estimate
(\ref{stimaiterata-secondagbis}), when $\sigma$ is onto, follow
in the same way.
\qed



\section{Regular solutions of the HJB equation}
\label{sec-HJB}

We show first, in Subsection \ref{SS:EXUNMILD}, that the HJB equation (\ref{HJBINTRO}) admits a unique mild solution $v$ which is $B$-Fr\'echet
differentiable.
Then (Subsection \ref{SS:SECONDDERIVATIVE}) we prove a further regularity result
whose proof is more complicated than the previous one
and that will be useful to solve our control problem in Section \ref{sec:contr-feedback}.

\subsection{Existence and uniqueness of mild solutions}
\label{SS:EXUNMILD}

We start
showing how to rewrite (\ref{HJBformale1}) in its
integral (or ``mild'') form as anticipated in the introduction, formula \myref{solmildHJB}.
Denoting by $\call$ the generator of the Ornstein-Uhlenbeck semigroup
$R_{t}$, we know that, for all $f\in C^{2}_b(\calh)$ such that $\nabla f \in D(A^*)$
(see e.g. \cite{CeGo} Section 5 or also \cite{DaPratoZabczyk95} Theorem 2.7):

\begin{equation}\label{eq:ell}
 \call[f](x)=\frac{1}{2} Tr \;GG^*\; \nabla^2f(x)
+ \< x,A^*\nabla f(x)\>.
\end{equation}
The HJB equation (\ref{HJBINTRO}) can then be formally rewritten as
\begin{equation}\label{HJBformale}
  \left\{\begin{array}{l}\dis
-\frac{\partial v(t,x)}{\partial t}=\call [v(t,\cdot)](x) +\ell_0(t,x)+
H_{min} (\nabla^B v(t,x)),\qquad t\in [0,T],\,
x\in \calh,\\
\\
\dis v(T,x)=\bar\phi(x_0).
\end{array}\right.
\end{equation}
By applying formally the variation of
constants formula we then have
\begin{equation}
v(t,x) =R_{T-t}[\phi](x)+\int_t^T R_{s-t}\left[
H_{min}(\nabla^B v(s,\cdot))+\ell_0(s,\cdot)\right](x)\; ds,\qquad t\in [0,T],\
x\in H,\label{solmildHJB}
\end{equation}
We use this formula to give the notion of mild
solution for the HJB equation (\ref{HJBformale}).

\begin{definition}\label{defsolmildHJB}
We say that a
function $v:[0,T]\times \calh\rightarrow\mathbb{R}$ is a mild
solution of the HJB equation (\ref{HJBformale}) if the following
are satisfied:
\begin{enumerate}

\item $v(T-\cdot, \cdot)\in C^{0,1,B}_{{1/2}}\left([0,T]\times \calh\right)$;

\item  equality (\ref{solmildHJB}) holds on $[0,T]\times \calh$.
\end{enumerate}
\end{definition}

\begin{remark}\label{rm:crescitapoli-HJB}
Since $C^{0,1,B}_{{1/2}}\left([0,T]\times \calh\right)
\subset C_b([0,T]\times \calh)$ (see Definition \ref{df4:Gspaces})
the above Definition \ref{defsolmildHJB} requires, among other properties, that a mild solution is continuous and bounded up to $T$.
This constrains the assumptions on the data, e.g. it implies that the final datum $\phi$ must be continuous and bounded.
As recalled in Remark \ref{rm:crescitapoli-gen}-(i) and (ii)
we may change this requirement in the above definition asking only polynomial growth in $x$
and/or measurability of $\phi$. Most of our main results will remain true
with straightforward modifications.
\end{remark}

Since the transition semigroup $R_t$ is not even strongly Feller
we cannot study the existence and uniqueness of a mild solution of equation (\ref{HJBformale}) as it is done e.g. in \cite{G1}.
We then use the partial smoothing property studied in Sections \ref{section-smoothOU} and
\ref{section-smooth-conv}.
Due to Lemma \ref{lemma reg-convpercontr} the right space where
to seek a mild solution seems to be $\Sigma^1_{T,{1/2}}$;
indeed our existence and uniqueness result will be proved by a fixed point argument in such space.

\begin{theorem}\label{esistenzaHJB}
Let Hypotheses \ref{ipotesibasic} and \ref{ipotesicostoconcreto} hold and
let (\ref{eq:hpdebreg}) or (\ref{eq:hpdebregbis}) hold.
Then the HJB equation (\ref{HJBformale})
admits a mild solution $v$ according to Definition \ref{defsolmildHJB}.
Moreover $v$
is unique among the functions $w$ such that $w(T-\cdot,\cdot)\in\Sigma_{T,1/2}$ and it satisfies, for suitable $C_T>0$, the estimate
\begin{equation}\label{eq:stimavmainteo}
\Vert v(T-\cdot,\cdot)\Vert_{C^{0,1,B}_{{1/2}}}\le C_T\left(\Vert\bar\phi \Vert_\infty
+\Vert\bar\ell_0 \Vert_\infty \right).
\end{equation}
Finally if the initial datum $\phi$ is also continuously $B$-Fr\'echet
(or Fr\'echet) differentiable,
then $v \in C^{0,1,B}_{b}([0,T]\times \calh)$ and, for suitable $C_T>0$,
\begin{equation}\label{eq:stimavmainteobis}
\Vert v\Vert_{C^{0,1,B}_{b}}\le C_T\left(\Vert\phi \Vert_\infty
+\Vert\nabla^B\phi \Vert_\infty+\Vert\ell_0 \Vert_\infty \right)
\end{equation}
(substituting $\nabla^B\phi$ with $\nabla\phi$ if $\phi$ is Fr\'echet differentiable).
\end{theorem}
\dim
We use a fixed point argument in $\Sigma^1_{T,{1/2}}$. To this aim,
first we rewrite (\ref{solmildHJB}) in a forward way. Namely
if $v$ satisfies \myref{solmildHJB} then, setting $w(t,x):=v(T-t,x)$ for any
$(t,x)\in[0,T]\times \calh$, we get that $w$ satisfies
\begin{equation}
  w(t,x) =R_{t}[\phi](x)+\int_0^t R_{t-s}[
H_{min}(
\nabla^B w(s,\cdot))+\ell_0(s,\cdot)
](x)\; ds,\qquad t\in [0,T],\
x\in H,\label{solmildHJB-forward}
\end{equation}
which is the mild form of the forward HJB equation
\begin{equation}\label{HJBformaleforward}
  \left\{\begin{array}{l}\dis
\frac{\partial w(t,x)}{\partial t}=\call [w(t,\cdot)](x) +\ell_0(t,x)+
H_{min} (\nabla^B w(t,x)),\qquad t\in [0,T],\,
x\in \calh,\\
\\
\dis w(0,x)=\phi(x).
\end{array}\right.
\end{equation}
Define the map $\calc$ on $\Sigma^1_{T,{1/2}}$  by setting, for $g\in \Sigma^1_{T,{1/2}}$,
\begin{equation}\label{mappaC}
 \calc(g)(t,x):=R_{t}[\phi](x)+\int_0^t R_{t-s}[
H_{min}(
\nabla^B g(s,\cdot))+\ell_0(s,\cdot)
](x)\; ds,\qquad t\in [0,T],
\end{equation}
By Proposition \ref{lemmaderhpdeb}
and Lemma \ref{lemma reg-convpercontr}-(i) we deduce that $\calc$
is well defined in $\Sigma^1_{T,{1/2}}$ and takes its values in $\Sigma^1_{T,{1/2}}$. Since in Lemma \ref{lemma:Sigma}
we have proved that $\Sigma^1_{T,{1/2}}$ is a closed subspace of $C^{0,1,B}_{{1/2}}([0,T]\times\calh)$, once we have proved that $\calc$ is a contraction,
by the Contraction Mapping Principle there exists a unique (in $\Sigma^1_{T,{1/2}}$) fixed point of the map $\calc$, which gives a mild solution of (\ref{HJBformale}).

\noindent Let $g_1,g_2 \in \Sigma^1_{T,{1/2}}$. We evaluate
$\Vert \calc(g_1)-\calc (g_2)\Vert_{\Sigma_{T,{1/2}}}=\Vert \calc(g_1)-\calc (g_2)\Vert_{C^{0,1,B}_{{1/2}}}$. First of all, arguing as in the proof of
Lemma \ref{lemma reg-convpercontr} we have, for every $(t,x)\in [0,T]\times H$,
\begin{align*}
  \vert \calc (g_1)(t,x)- \calc(g_2)(t,x) \vert& =\left\vert \int_0^t R_{t-s}\left[H_{min}\left(\nabla^B g_1(s,\cdot)\right)
 -H_{min}\left(\nabla^B g_2(s,\cdot)\right)\right](x)ds\right\vert\\
 &\le \int_0^t s^{-{1/2}} L \sup_{y \in H}\vert s^{{1/2}}\nabla^B (g_1-g_2)(s,y)\vert ds
  \leq 2Lt^{{1/2}}\Vert g_1-g_2 \Vert_{C^{0,1,B}_{{1/2}}}.
\end{align*}
Similarly, arguing exactly as in the proof of (\ref{stimaiterata-primag}), we get
\begin{align*}
t^{{1/2}}\vert \nabla^B\calc (g_1)(t,x) &- \nabla^B\calc(g_2)(t,x) \vert =
t^{{1/2}}\left\vert \nabla^B\int_0^t  R_{t-s}\left[H_{min}
\left(\nabla^B g_1(s,\cdot)\right)-H_{min}\left(\nabla^B
g_2(s,\cdot)\right)\right](x)ds\right\vert\\
& \leq t^{{1/2}} L \Vert g_1-g_2 \Vert_{C^{0,1,B}_{{1/2}}} \int_0^t (t-s)^{-{1/2}}
s^{-{1/2}} ds
\le t^{{1/2}}L \beta\left({1/2} , {1/2}\right)\Vert g_1-g_2 \Vert_{C^{0,1,B}_{{1/2}}}.
\end{align*}
Hence, if $T$ is sufficiently small, we get
\begin{equation}\label{stima-contr}
 \left\Vert \calc (g_1)-\calc(g_2)\right\Vert _{C^{0,1,B}_{{1/2}}
  }\leq C
\left\Vert g_1-g_2\right\Vert _{C^{0,1,B}_{{1/2}}
}
\end{equation}
with $C<1$. So the map $\calc$ is a contraction in $\Sigma^1_{T,{1/2}}$
and, if we denote by $w$ its unique fixed point, then $v:=w(T-\cdot,\cdot)$
turns out to be a mild solution of the HJB equation (\ref{HJBformale}),
according to Definition \ref{solmildHJB}.

Since the constant $L$ is independent of $t$, the case of generic $T>0$ follows
by dividing the interval $[0,T]$
into a finite number of subintervals of length $\delta$ sufficiently small, or equivalently, as done in \cite{Mas},
by taking an equivalent norm with an adequate exponential weight, such as
\[
 \left\Vert f\right\Vert _{\eta,C^{0,1,B}_{{1/2}}
 }=\sup_{(t,x)\in[0,T]\times \calh}
\vert e^{\eta t}f(t,x)\vert+
\sup_{(t,x)\in (0,T]\times \calh}  e^{\eta t}t^{{1/2}}
\left\Vert \nabla^B f\left(  t,x\right)  \right\Vert _{(\R^m)^{\ast}},
\]


The estimate (\ref{eq:stimavmainteo}) follows from Proposition \ref{cor-der} and Lemma \ref{lemma reg-convpercontr}.

Finally the proof of the last statement follows observing that, if
$\phi$ is continuously $B$-Fr\'echet (or Fr\'echet) differentiable,
then $R_t[\phi]$ is continuously $B$-Fr\'echet differentiable
with $\nabla^B R_t[\phi]$ bounded in $[0,T]\times \calh$, see lemma \ref{lemma-reg-R_t},
formula (\ref{eq:derivateprimelisce}).
This allows to perform the fixed point, exactly as done in the first part of the proof,
in $C^{0,1,B}_b([0,T]\times \calh)$ and to prove estimate \myref{eq:stimavmainteobis}.
\qed

\begin{corollary}\label{diffle-corollario}
Let Hypotheses \ref{ipotesibasic} and \ref{ipotesicostoconcreto} hold and let $\sigma$ be onto.
Then the mild solution of equation (\ref{HJBformale})
found in the previous theorem
is also Fr\'echet differentiable,
and the following estimate holds true
\begin{equation}\label{stimadiffle-solHJB}
 \left\Vert v(T-\cdot,\cdot)\right\Vert _{C^{0,1}_{{1/2}}
  }\leq C_T\left(\Vert\bar\phi \Vert_\infty
+\Vert\bar\ell_0 \Vert_\infty \right)
\end{equation}
for a suitable $C_T>0$.
\end{corollary}
\dim
Let $v$ be the mild solution of equation (\ref{HJBformale}), and $\forall
t\in[0,T],\,x\in\calh$ define $w(t,x):=v(T-t,x)$, so that $w$ satisfies
(\ref{solmildHJB-forward}), so that by applying
the last statement of Lemma \ref{lemma reg-convpercontr}
it is immediate to see that $w\in C^{0,1}_{{1/2}}([0,T]\times \calh)$.
By differentiating (\ref{solmildHJB-forward}) we get
\begin{equation*}
 \nabla w(t,x) =\nabla R_{t}[\phi](x)+\nabla\int_0^t R_{t-s}[
H_{min}(
\nabla^B w(s,\cdot)+\ell_0(s,\cdot)
](x)\; ds,\qquad t\in [0,T],\
x\in H,
\end{equation*}
By Lemma \ref{lemma-reg-R_t}, the above recalled variation of Lemma
\ref{lemma reg-convpercontr} and estimate (\ref{stimaiterata-primagbis}), we get that
\begin{equation*}
 \vert\nabla w(t,x) \vert\leq C t^{-{1/2}}\Vert\phi\Vert_\infty+Ct^{{1/2}}\left(1+\Vert
 w\Vert_{C^{0,1,B}_{{1/2}}}+\Vert\ell_0\Vert_\infty\right),
\qquad t\in [0,T],\
x\in \calh,
\end{equation*}
which gives the claim using the estimate for
$\Vert w\Vert_{C^{0,1,B}_{{1/2}}}$ given in (\ref{eq:stimavmainteo}).
\qed

\subsection{Second derivative of mild solutions}
\label{SS:SECONDDERIVATIVE}

The further regularity result we are going to prove will be needed in Section \ref{sec:contr-feedback}, Theorem
\ref{teo su controllo feedback}.
A similar result can be found in \cite{G1}, Section 4.2. Here we use the same line of proof
but we need to argue in a different way to get the apriori estimates.

\begin{theorem}\label{lemma-stimev}
Let Hypotheses \ref{ipotesibasic}, \ref{ipotesicostoconcreto} and
\ref{ipotesicostoconcretobis} hold.
Let also (\ref{eq:hpdebreg}) or (\ref{eq:hpdebregbis}) hold.
Let $v$ be the mild solution of the HJB equation (\ref{HJBformale}) as from
Theorem \ref{esistenzaHJB}. Then we have the following.
\begin{itemize}
  \item[(i)]
If $\phi$ is continuously differentiable then we have $v \in \Sigma^2_{T,{1/2}}$,
hence the second order derivatives $\nabla^B\nabla v$ and $\nabla\nabla^B v$
exist and are equal.
Moreover there exists a constant $C>0$ such that
\begin{equation}\label{stimanablav}
\vert \nabla v(t,x)\vert \leq C\left( \Vert \nabla\bar\phi\Vert_\infty+
\Vert \nabla\bar\ell_0\Vert_\infty\right),
\end{equation}
\begin{equation}\label{stimanablav^2}
\vert \nabla^B\nabla v(t,x)\vert= \vert \nabla\nabla^B v(t,x)\vert
\leq C\left( (T-t)^{-{1/2}}
\Vert \nabla\bar\phi\Vert_\infty+(T-t)^{{1/2}}
\Vert \nabla\bar\ell_0\Vert_\infty\right).
\end{equation}
Finally, if $\sigma$ is onto, then also $\nabla^2 v$ exists and is continuous
and, for suitable $C>0$,
\begin{equation}\label{stimanablav^2bis}
\vert \nabla^2 v(t,x)\vert
\leq C\left( (T-t)^{-{1/2}}
\Vert \nabla\bar\phi\Vert_\infty+(T-t)^{{1/2}}
\Vert \nabla\bar\ell_0\Vert_\infty\right).
\end{equation}

  \item[(ii)]
If $\phi$ is only continuous then the function
$(t,x)\mapsto (T-t)^{1/2}v(t,x)$ belongs to $\Sigma^2_{T,{1/2}}$.
Moreover there exists a constant $C>0$ such that
\begin{equation}\label{stimanablavreg}
\vert \nabla v(t,x)\vert \leq C\left((T-t)^{-1/2} \Vert\bar\phi\Vert_\infty+
\Vert \nabla\bar\ell_0\Vert_\infty\right),
\end{equation}
\begin{equation}\label{stimanablav^2reg}
\vert \nabla^B\nabla v(t,x)\vert= \vert \nabla\nabla^B v(t,x)\vert
\leq C\left( (T-t)^{-1}
\Vert \bar\phi\Vert_\infty+(T-t)^{-1/2}
\Vert \nabla\bar\ell_0\Vert_\infty\right).
\end{equation}
Finally, if $\sigma$ is onto, then also $\nabla^2 v$ exists and is continuous
in $[0,T)\times \calh$
and, for suitable $C>0$,
\begin{equation}\label{stimanablav^2bisreg}
\vert \nabla^2 v(t,x)\vert
\leq C\left( (T-t)^{-1}
\Vert \nabla\bar\phi\Vert_\infty+(T-t)^{-1/2}
\Vert \nabla\bar\ell_0\Vert_\infty\right).
\end{equation}
\end{itemize}
\end{theorem}	
\dim
We start proving (i) by applying the Contraction Mapping Theorem in a closed ball $B_T(0,R)$
($R$ to be chosen later) of the space $ \Sigma^2_{T,{1/2}}$.
By Proposition \ref{lemma-reg-R_t} and Lemma \ref{lemma reg-convpercontr} we deduce that
the map $\calc$ defined in \myref{mappaC}
brings $\Sigma^2_{T,{1/2}}$ into $\Sigma^2_{T,{1/2}}$.
Moreover, for every $g \in \Sigma^2_{T,{1/2}}$, we get, first using \myref{eq:Hlip},
\begin{align*}
&  \vert \calc (g)(t,x) \vert\le
\left\vert
R_{t}[\phi](x)
\right\vert
+
\left\vert \int_0^t R_{t-s}\left[H_{min}\left(\nabla^B g(s,\cdot)\right)
+ \ell_0(s,\cdot)\right](x)ds
\right\vert
\leq
\Vert \phi\Vert_\infty
+t L\left(1+\Vert g \Vert_{C^{0,1,B}_{0}}\right) + t\|\ell_0\|_\infty;
\end{align*}
second by \myref{eq:derivateprimelisce}, \myref{eq:perdersecvera}, \myref{eq:Hlip}
(calling $M:=\sup_{[0,T]}\|e^{tA}\|$)
\begin{align}
\nonumber
  \vert\nabla \calc (g)(t,x) \vert &\le
\left\vert
\nabla R_{t}[\phi](x)
\right\vert
+
\left\vert \nabla\int_0^t R_{t-s}\left[H_{min}\left(\nabla^B g(s,\cdot)\right)
+ \ell_0(s,\cdot)\right](x)ds
\right\vert\\[2mm]
&\le M \Vert \nabla\phi\Vert_\infty+M
\int_0^t \left[\left\Vert\nabla H_{min}\left(\nabla^B g(s,\cdot)\right)
\nabla\nabla^B g(s,\cdot)\right\Vert_\infty + \|\nabla\ell_0(s,\cdot)\|_\infty\right]
ds
\nonumber
\\
& \leq M\left[ \Vert \nabla\phi\Vert_\infty + t\Vert \nabla\ell_0\Vert_\infty
 + L t^{{1/2}}\Vert g \Vert_{C^{0,2,B}_{{1/2}}}\right];
\label{eq:stimadersecnew1}
\end{align}
third by \myref{eq:stimaDRT2-gen} (with
(\ref{eq:hpdebreg}) or (\ref{eq:hpdebregbis})), and \myref{stimaiterata-secondag}
\begin{align}
\nonumber
&t^{{1/2}}\vert \nabla^B\nabla\calc (g)(t,x)\vert
\le
t^{{1/2}}
\left\vert
\nabla^B\nabla R_{t}[\phi](x)
\right\vert
+
t^{{1/2}} \left\vert \nabla^B\nabla
\int_0^t  R_{t-s}\left[H_{min}\left(\nabla^B g(s,\cdot)\right)+\ell_0(s,\cdot)
\right](x)ds\right\vert
\\[2mm]
& \leq C \Vert \nabla\phi\Vert_\infty
+ C t^{{1/2}}\int_0^t (t-s)^{-{1/2}}\|\nabla \ell_0(s,\cdot)\|_\infty ds
+ Ct^{{1/2}} \| g\|_{C^{0,2,B}_{1/2}}
\\[2mm]
& \leq C\left[ \Vert \nabla\phi\Vert_\infty+2t\Vert \nabla\ell_0\Vert_\infty
+t^{{1/2}} \Vert g \Vert_{C^{0,2,B}_{{1/2}}}\right],
\label{eq:stimadersecnew2}
\end{align}
with the constant $C$ (that may change from line to line) given by the quoted estimates.
Hence, for $g\in B_T(0,R)$, we get, for given $C_1>0$,
\begin{align}
\nonumber
\Vert \calc (g) \Vert_{C^{0,2,B}_{{1/2}}}
&\le
C_1\left[\Vert\phi\Vert_{C^1_b} +T+ T\Vert\ell_0\Vert_{C^1_b}\right]
+TL \Vert g \Vert_{C^{0,1,B}_{0}}+(ML+C)T^{1/2}\Vert g \Vert_{C^{0,2,B}_{{1/2}}}
\\
\label{eq:dersecstima1}
&
\le
C_1\left[\Vert\phi\Vert_{C^1_b} +T+ T\Vert\ell_0\Vert_{C^1_b}\right]
+\rho(T)R
\end{align}
where we define
\begin{equation}\label{eq:rhodef}
\rho(T):= TL+(ML+C)T^{1/2}.
\end{equation}
Now take $g_1,g_2 \in \Sigma^2_{T,{1/2}}$.
Arguing as in the above estimates we have, for every $(t,x)\in [0,T]\times \calh$,
\begin{align*}
\vert \calc (g_1)(t,x)- \calc(g_2)(t,x) \vert&=\left\vert \int_0^t R_{t-s}\left[H_{min}\left(\nabla^B g_1(s,\cdot)\right)
-H_{min}\left(\nabla^B g_2(s,\cdot)\right)\right](x)ds\right\vert\\[2mm]
&\leq t L\Vert g_1-g_2 \Vert_{C^{0,1,B}_{0}}
\end{align*}
\begin{align*}
&\vert \nabla\calc (g_1)(t,x)- \nabla\calc(g_2)(t,x) \vert =
\left\vert \nabla\int_0^t  R_{t-s}\left[H_{min}\left(\nabla^B g_1(s,\cdot)\right)
-H_{min}\left(\nabla^B g_2(s,\cdot)\right)\right](x)ds
\right\vert
\\
&\le M\int_0^t \left\Vert\nabla H_{min}\left(\nabla^B g_1(s,\cdot)\right)
\nabla\nabla^B g_1(s,\cdot)-
\nabla H_{min}\left(\nabla^B g_2(s,\cdot)\right)\nabla\nabla^B g_2(s,\cdot)\right\Vert_\infty
ds
\\
& \leq  2ML t^{{1/2}}\left[\Vert g_1-g_2 \Vert_{C^{0,1,B}_{0}}
\Vert g_1 \Vert_{C^{0,2,B}_{{1/2}}}+\Vert g_1-g_2 \Vert_{C^{0,2,B}_{{1/2}}} \right]
\end{align*}
and, using \myref{eq:der-seconda-conv-Sigma},
\begin{align*}
&t^{{1/2}}\vert \nabla^B\nabla\calc (g_1)(t,x)- \nabla^B\nabla\calc(g_2)(t,x) \vert =
t^{{1/2}} \left\vert \nabla^B\nabla\int_0^t  R_{t-s}\left[H_{min}\left(\nabla^B g_1(s,\cdot)\right)
-H_{min}\left(\nabla^B g_2(s,\cdot)\right)\right](x)ds\right\vert
\\[2mm]
& \leq  t^{{1/2}}M L \beta\left({1/2} , {1/2}\right)
\left[\Vert g_1-g_2 \Vert_{C^{0,1,B}_{0}}
\Vert g_1 \Vert_{C^{0,2,B}_{{1/2}}}+\Vert g_1-g_2 \Vert_{C^{0,2,B}_{{1/2}}} \right].
%
\end{align*}
Hence, for $g_1,g_2\in B_T(0,R)$, we have, recalling \myref{eq:rhodef} and the way $C$ is found
in \myref{eq:stimadersecnew2},
\begin{align}
\nonumber
\left\Vert \calc (g_1)-\calc(g_2)\right\Vert _{C^{0,2,B}_{{1/2}}([0,T]\times \calh)  }
&\leq L\left(T+ MT^{{1/2}}\left(2+\beta\left({1/2} , {1/2}\right)\right)
(1+R) \right)
\left\Vert g_1-g_2\right\Vert _{C^{0,2,B}_{{1/2}}  }
\\[2mm]
&\le \rho(T)(1+R)\left\Vert g_1-g_2\right\Vert _{C^{0,2,B}_{{1/2}}  },
\label{eq:dersecstima2}
\end{align}
Now, by \myref{eq:dersecstima1} and \myref{eq:dersecstima2},
choosing any $R>C_1\left[\Vert\phi\Vert_{C^1_b} +T+ T\Vert\ell_0\Vert_{C^1_b}\right]$
we can find $T_0$ sufficiently small so that
$\rho(T_0)<1/2$ and so, thanks to \myref{eq:dersecstima1}
and \myref{eq:dersecstima2}, $\calc$ is a contraction in $B_{T_0}(0,R)$.
Let then $w$ be the unique fixed point of $\calc$ in $B_{T_0}(0,R)$:
it must coincide with $v(T-\cdot,\cdot)$ for $t\in [0,T_0]$.
This procedure can be iterated arriving to cover the whole interval $[0,T]$
if we give an apriori estimate for the norm $\Vert w \Vert_{C^{0,2,B}_{{1/2}}}$.
By the last statement of Theorem \ref{esistenzaHJB} we already have an apriori estimate for
$\Vert w \Vert_\infty+\Vert \nabla^B w \Vert_\infty$.
To get the estimate for $\nabla w$ and $\nabla^B\nabla w $
we use the \myref{eq:stimadersecnew1} and \myref{eq:stimadersecnew2}
where we put $w$ in place of $\calc g$ and $g$.
From the first line of \myref{eq:stimadersecnew2} and \ref{eq:der-seconda-conv-Sigma}
we get
$$
\Vert \nabla^B \nabla w(t,\cdot)\Vert_\infty
\leq C\left[ \Vert \nabla\phi\Vert_\infty+2T\Vert \nabla\ell_0\Vert_\infty \right]
+L\int_0^t (t-s)^{-{1/2}}
\Vert \nabla^B \nabla w(s,\cdot)\Vert_\infty ds
$$
which, thanks to the Gronwall Lemma  (see \cite{Henry}, Subsection 1.2.1, p.6)
give the apriori estimate for $\nabla^B \nabla w$. Then from the second line of
\myref{eq:stimadersecnew1} we get
$$
\Vert \nabla w(t,\cdot)\Vert_\infty
\le  M\left[ \Vert \nabla\phi\Vert_\infty + T\Vert \nabla\ell_0\Vert_\infty\right]
+M L\int_0^t \Vert \nabla^B \nabla w(s,\cdot)\Vert_\infty ds
$$
which gives the apriori estimate for $\nabla w$ using the previous one for $\nabla^B \nabla w$.
Estimate \myref{stimanablav^2bis} follows by repeating the same arguments above
but replacing $\nabla^B \nabla $ with $\nabla^2$.

We now prove (ii).
Let $v$ be the mild solution of
(\ref{HJBformale}) and, for all $\eps \in ]0, T[$, $x \in \calh$, call $\phi^\eps (x) =
v(T-\eps,x)$. Then $v$ is the unique mild solution, on $[0,T-\eps]\times \calh$,
of the equation (for $t \in [0,T]$, $x \in H$)
\begin{equation}
\label{eq:HJBmildeps}
v(t,x) =   R_{T-\eps-t}  \phi^\eps (x)
+ \int_{t}^{T-\eps}
R_{s-t}\left[H_{min}(\nabla^B v(s,\cdot ))+\ell_0(s,\cdot) \right](x)ds
\end{equation}
This fact can be easily seen by applying the semigroup property of $R_t$
(see e.g. \cite{G1} Lemma 4.10 for a completely similar result).

Now, by Theorem \ref{esistenzaHJB}, $\phi^\eps$ is continuously differentiable,
so we can apply part (i) of this theorem to \myref{eq:HJBmildeps}
getting the required $C^2$ regularity.
Estimates
\myref{stimanablavreg}-\myref{stimanablav^2reg}-\myref{stimanablav^2bisreg}
follows using estimates
\myref{stimanablav}-\myref{stimanablav^2}-\myref{stimanablav^2bis}
with $\phi^\eps$ in place of $\phi$ and then using the arbitrariness of $\eps$ and applying
\myref{eq:stimavmainteo} to estimate $\phi^\eps$ in term of $\phi$.
\qed

\section{Verification Theorem and Optimal Feedbacks}\label{sec-verifica}

The aim of this section is to provide a verification theorem and the existence of optimal feedback controls for our problem. This in particular will imply that the mild solution $v$ of the HJB equation (\ref{HJBformale}) built in Theorem \ref{esistenzaHJB} is equal to the
value function $V$ of our optimal control problem.

The main tool needed to get the wanted results is an identity (often called ``{\em fundamental identity}'', see equation (\ref{relfond})) satisfied by the solutions of the HJB equation. When the solution is smooth enough
(e.g. it belongs to $UC_{b}^{1,2}\left(  \left[0,T\right]  \times H\right)$)
such identity is easily proved using the Ito's Formula. Since in our case the value
function does not possess this regularity, we proceed by approximation.
Due to the features of our problem (lack of smoothing and of the structure condition)
the methods of proof used in the literature do not apply here.
We will discuss the main issues along the way.

\subsection{$\calk$-convergence and approximations of solutions of the HJB equation}

We first introduce the notion of $\calk$-convergence, following\cite{Ce95} and \cite{G1}.
\begin{definition}\label{k-conv}
A sequence $(f_n)_{n\geq 0}\in C_b(\calh)$
is said to be $\calk$-convergent to a function $f\in C_b(H)$ (and we shall write
$f_n\overset{\calk}{\rightarrow} f$ or $f=\calk-\lim_{n\rightarrow\infty}f_n$) if for any compact set
$\calk\subset\calh$
$$
\sup_{n\in\N}\Vert f_n\Vert_\infty<+\infty \quad {\rm and } \quad
\lim_{n\rightarrow\infty}\sup_{x\in\calk}\vert f(x)-f_n(x)\vert =0.
$$
\noindent
Similarly, given $I\subseteq \R$, a sequence
$(f_n)_{n\geq 0}\in C_b(I \times \calh)$ is said to be $\calk$-convergent to a function $f\in C_b(I\times \calh)$ (and we shall write again
$f_n\overset{\calk}{\rightarrow} f$ or $f=\calk-\lim_{n\rightarrow\infty}f_n$) if for any compact set
$\calk\in\calh$ and for any $(t,x)\in I\times\calk $ we have
$$
\sup_{n\in\N}\Vert f_n\Vert_\infty<+\infty \quad {\rm and } \quad
\lim_{n\rightarrow\infty}\sup_{(t,x)\in I\times\calk}\vert f(t,x)-f_n(t,x)\vert=0.
$$
\end{definition}


Now we recall the definition (given in  \cite{DP3}, beginning of Chapter 7)
of \textit{strict solution}
of a family of Kolmogorov equations.
Consider the following forward Kolmogorov equation with unknown $w$:
\begin{equation}\label{eq:Kforward}
  \left\{\begin{array}{l}\dis
\frac{\partial w(t,x)}{\partial t}=\frac{1}{2} Tr \;GG^*\; \nabla^2w(t,x)
+ \< Ax,\nabla w(t,x)\>
+\calf(t,x)
\qquad t\in [0,T],\,
x\in \calh,\\
\\
\dis w(0,x)=\phi(x).
\end{array}\right.
\end{equation}
where
the functions
$\calf:[0,T]\times \calh \to \R$ and $\phi:\calh\to\R$ are bounded and continuous.

\begin{definition}\label{df:strictandpistrong}
By \textit{strict solution} of the Kolmogorov equation (\ref{eq:Kforward}) we mean
a function $w$ such that
\begin{equation}\label{propr-strict-sol}
 \left\lbrace
\begin{array}{l}
w \in C_b([0,T]\times \calh)\quad  {\rm and}\quad w(0,x)=\phi(x)
\\
w(t,\cdot) \in UC^2_b(\calh), \; \forall t\in [0,T];
\\
w(\cdot, x)\in C^{1}([0,T]),\;\forall x\in D(A)\; \hbox{and $w$ satisfies (\ref{eq:Kforward}).}
\end{array}
\right.
\end{equation}
\end{definition}

Now we prove a key approximation lemma.

\begin{lemma}\label{lm:approximation}
Let Hypothesis \ref{ipotesibasic} and \ref{ipotesicostoconcreto} hold.
Let also (\ref{eq:hpdebreg}) or (\ref{eq:hpdebregbis}) hold.
Let $v$ be the mild solution of the HJB equation (\ref{HJBformale})
and set $w(t,x)=v(T-t,x)$ for $(t,x)\in [0,T]\times \calh$.
Then there exist three sequences of functions $(\bar\phi_n)$, $(\bar \ell_{0,n})$ and $(\calf_n)$ such that, for all $n\in \N$,
    \begin{equation}\label{eq:defapproxphif}
        \bar\phi_n \in C_c^\infty(\R^n), \qquad
        \bar \ell_{0,n} \in C_c^\infty([0,T]\times \R^n), \qquad
        \calf_n  \in C_c^\infty([0,T]\times \calh) \cap \Sigma_{T,{1/2}}
    \end{equation}
and
    \begin{equation}\label{eq:convapproxphif}
        \bar\phi_n \rightarrow \bar\phi ,
        \qquad         \bar \ell_{0,n} \rightarrow \bar \ell ,
        \qquad\calf_n \to H_{min} (\nabla^B w)
    \end{equation}
in the sense of $\calk$-convergence.
Moreover, defining $\phi_n(x)=\bar\phi_n(x_0)$, $\ell_{0,n}(s,x)=\bar\ell_{0,n}(s,x_0)$
and
\begin{equation}\label{v_n-mild}
w_n(t,x):= R_t \phi_n + \int_0^t R_{t-s}[\calf_n (s,\cdot)+\ell_{0,n}(s,\cdot)] (x)ds
\end{equation}
the following hold:
\begin{itemize}
  \item $w_n \in UC^{1,2}_b ([0,T]\times \calh) \cap \Sigma_{T,{1/2}}$,
  \item $w_n$ is a strict solution of
  \myref{eq:Kforward} with $\phi_n$ in place of $\phi$ and
$\calf_n + \ell_{0,n}$ in place of $\calf $,
  \item we have, in the sense of $\calk$-convergence
  (the first in $[0,T]\times \calh$,
  the second in $(0,T]\times \calh$),
  \begin{equation}\label{eq:convapproxphifsol}
        w_n \rightarrow w ,
        \qquad t^{{1/2}}\nabla^B w_n \to t^{{1/2}}\nabla^B w.
    \end{equation}
\end{itemize}
\end{lemma}
\dim
We divide the proof in three steps.

{\bf Step 1: choosing the three approximating sequences.}
We choose $\bar\phi_n$ and $\bar\ell_{0,n}$ to be
the standard approximations by convolution of $\bar \phi$ and $\bar \ell_0$.
To define $\calf_n$ we observe first that, since $w \in \Sigma^1_{T,{1/2}}$, then the function
$$
(t,x) \mapsto  \calf(t,x):= H_{min} (\nabla^B w(t,x))
$$
has the property that there exist $f:[0,T)\times \R^n \to \R$, continuous and bounded, such that
$$
\calf (t,x)=t^{-{1/2}} f(t,(e^{tA}x)_0).
$$
We then let $f_n$ be the approximation by convolution of $f$ and define
$$
\calf_n (t,x)=t^{-{1/2}} f_n(t,(e^{tA}x)_0).
$$
{\bf Step 2: proof that $w_n\in UC^{1,2}_b([0,T]\times\calh) \cap \Sigma^1_{T,{1/2}}$ and
that it is a strict solution.}
The fact that $w_n\in \Sigma_{T,{1/2}}$ follows immediately from \myref{v_n-mild},
Proposition \ref{lemmaderhpdeb} and Lemma \ref{lemma reg-convpercontr}-(i).
Differentiability with respect to the variable $x$ follows
arguing as in the proof of formula (\ref{eq:derivateprimelisce})
in Lemma \ref{lemma-reg-R_t}: indeed differentiating twice
we get that, for all $h,k \in \calh$
\begin{align*}
& \<\nabla^2 R_{t}\left[\phi_n\right]  (x)h,k\>_\calh=
R_{t}\left[\<\nabla^2 \phi_n e^{tA}h,e^{tA}k\>_\calh\right ]  (x).
\end{align*}
Similarly we have, for the convolution term containing $\calf_n$,
\begin{align*}
& \<\nabla \int_0^t R_{t-s}\left[\calf_n(s,\cdot)\right]  (x),h\>_\calh
\\[2mm]
 &  =\lim_{\alpha\rightarrow 0}\frac{1}{\alpha}
 \int_0^t\left[\int_{\calh}\calf_n\left(s,z+e^{(t-s)A}(x+\alpha h)\right)
 \caln\left(0,Q_{t-s}\right)  (dz)
 -\int_{\calh}
 \calf_n\left(s,z+e^{(t-s)A}x\right)  \caln\left(  0,Q_{t-s}\right)  (dz) \,ds\right]
 \\[2mm]
 &  =\int_0^t\int_{\calh}\lim_{\alpha\rightarrow 0}\frac{1}{\alpha}
 \left[\calf_n\left(s,z+e^{(t-s)A}(x+\alpha h)\right)-\calf_n\left(s,z+e^{(t-s)A}x\right)\right]
 \caln\left(0,Q_{t-s}\right)  (dz)\,ds\\[2mm]
 &  =\int_0^t\int_{\calh}\<\nabla \calf_n\left(s,z+e^{(t-s)A}x\right), e^{(t-s)A}h\>_\calh
 \caln\left(0,Q_{t-s}\right)  (dz)\,ds\\[2mm]
 &=\int_0^tR_{t-s}\left[\<\nabla \calf_n(s,\cdot) ,e^{(t-s)A}h\>_\calh\right ]  (x)\,ds,
\end{align*}
and also, arguing in the same way,
\begin{align*}
& \<\nabla^2\int_0^t R_{t}\left[\calf_n(s,\cdot)\right]  (x)\,ds\, h,k\>\calh=
\int_0^tR_{t-s}\left[\<\nabla^2 \calf_n(s,\cdot) e^{(t-s)A}h,e^{(t-s)A}k\>_\calh\right ]  (x)\,ds.
\end{align*}
The  convolution term involving $\ell_{0,n}$ is treated exactly in the same way.
The proof that $w_n$ is differentiable with respect to time and that
$w_{nt}\in UC_b([0,T]\times\calh) $ is completely analogous
to what is done in \cite[Theorems 9.23 and 9.25]{DP1})
for homogeneous Kolmogorov equations and we omit it\footnote{The proof
in the nonhomogeneous case can be found in the forthcoming book \cite[Section 4.4]{FabbriGozziSwiech}.}.
By Theorem 5.3 in \cite{CeGo}, see also Theorem 7.5.1 in \cite{DP3} for Kolmogorov equations,
we finally conclude that $w_n$
is a strict solution to equation \ref{eq:Kforward}.

{\bf Step 3: proof of the convergences.}
First we prove that the sequences $(w_n)$ and $(t^{1/2}\nabla^B w_n)$
are bounded uniformly with respect to $n$. Indeed,
by (\ref{v_n-mild}) and by the properties of convolutions,
\begin{align*}
&\vert w_n(t,x)\vert
\leq \Vert \bar\phi_n\Vert_\infty +
\int_0^t\sup_{x\in\calh}\left[\vert\calf_n (s,x)\vert+\vert
\ell_{0,n}(s,x)\vert\right]\,ds
\\
&\leq
\Vert \bar\phi\Vert_\infty + \int_0^t\sup_{y\in\R^n}\left[s^{-{1/2}} \vert
f_n (s,y)\vert+\vert
\bar \ell_{0,n}(s,y)\vert\right]\,ds \leq
\Vert \bar\phi\Vert_\infty + \int_0^t\sup_{y\in\R^n}\left[ s^{-{1/2}}\vert f (s,y)\vert+\vert
\bar \ell_{0}(s,y)\vert\right]\,ds
\end{align*}
Moreover, using Proposition \ref{lemmaderhpdeb} and
\myref{eq:derBconvnew} for $\psi=identity$, we get
\begin{equation*}
 \vert t^{1/2} \nabla^B w_n(t,x)\vert
 \leq C\Vert \bar\phi\Vert_\infty
 +Ct^{1/2}\int_0^ts^{-{1/2}}(t-s)^{-{1/2}}\sup_{y\in\R^n}\left(\vert
   f\left(s, y\right)\vert+  \vert\bar \ell_{0}\left(s, y\right)
 \vert\right)
  ds,
\end{equation*}
for a suitable $C>0$.
Now, with similar computations, we prove the convergences. Indeed,
\begin{align*}
 \vert v_n(t,x)-v(t,x)\vert&\leq \int_\calh\vert \bar\phi_n\left( (e^{tA}z)_0+(e^{tA}x)_0\right)-
\bar\phi\left( (e^{tA}z)_0+(e^{tA}x)_0\right)\vert\caln(0,Q_{t})(dz) \\
& +
\int_0^t\int_\calh\left[\vert s^{-{1/2}} \calf_n \left(s,(e^{sA}z)_0+(e^{tA}x)_0\right)
-
s^{-{1/2}} \calf \left(s,(e^{sA}z)_0+(e^{tA}x)_0\right)\vert\right.\\
&\left.+\vert
\ell_{0,n}\left(s,(e^{sA}z)_0+(e^{tA}x)_0\right)-
\ell_0\left(s,(e^{sA}z)_0+(e^{tA}x)_0\right)
\vert\right]\caln(0,Q_{t-s})(dz)\,ds.
\end{align*}
Since, for every compact set $\calk \subset \calh$ the set
$\{(e^{tA}x)_0, \; t \in [0,T],\, x \in \calk\}$ is compact in $\R^n$,
then by the Dominated Convergence Theorem we get that for any compact set $\calk \subset \calh$
\begin{equation}
\label{eq:convvnnew}
\sup_{(t,x)\in[0,T]\times\calk} \vert v_n(t,x)-v(t,x)\vert\rightarrow 0.
\end{equation}
Moreover, using again Propositions \ref{lemmaderhpdeb}
and \ref{lemma reg-convpercontr}-(i),
for a suitable $C>0$, we get
\begin{align*}
 &\vert  \nabla^B v_n(t,x)- \nabla^B v(t,x)\vert\\
&\leq C
\int_\calh\left( \bar\phi_n\left( (e^{tA}z)_0+(e^{tA}x)_0\right)-
\bar\phi\left( (e^{tA}z)_0+(e^{tA}x)_0\right)\right)\<(Q^0_{t})^{-{1/2}}
\left( e^{tA} Bh\right)_0, (Q^0_{t})^{-{1/2}}z_0\>_{\R^n}\caln(0,Q_{t})(dz)\\
  & +C\int_0^t
  \int_\calh  \left(\left\vert
 s^{-{1/2}}  f_n\left(s, (e^{sA}z)_0+(e^{tA}x)_0\right)
 -s^{-{1/2}}  f\left(s, (e^{sA}z)_0+(e^{tA}x)_0\right)
 \right\vert\right.\\
  &\left. +\left\vert \bar \ell_{0,n}\left(s, (e^{sA}z)_0+(e^{tA}x)_0\right)
 -\bar \ell_{0}\left(s, (e^{sA}z)_0+(e^{tA}x)_0\right)\right\vert\right)\\
 &\quad \left \vert
 \<(Q^0_{t-s})^{-{1/2}}
 \left( e^{(t-s)A} Bh\right)_0, (Q^0_{t-s})^{-{1/2}}z_0\>_{\R^n}
  \right\vert\caln(0,Q_{t-s})(dz) ds.
\end{align*}
By Proposition \ref{lemmaderhpdeb} we know that for suitable $C>0$,
\[
 \int_\calh\vert\<(Q^0_{t})^{-{1/2}}
\left( e^{tA} Bh\right)_0, (Q^0_{t})^{-{1/2}}z_0\>_{\R^n}\vert^2\caln(0,Q_{t})(dz)
\leq C\Vert (Q^0_{t})^{-{1/2}}
\left( e^{tA} Bh\right)_0\Vert^2\leq \frac{C}{t}.
\]
Hence, applying Cauchy-Schwartz inequality we get
\begin{align*}
&\vert  \nabla^B v_n(t,x)- \nabla^B v(t,x)\vert\\
&\leq C t^{-{1/2}}
 \vert\int_\calh\vert \bar\phi_n\left( (e^{tA}z)_0+(e^{tA}x)_0\right)-
 \bar\phi_n\left( (e^{tA}z)_0+(e^{tA}x)_0\right)\vert^2\caln(0,Q_{t})(dz)\vert^{1/2}\\
  &+\int_0^t (t- s)^{-{1/2}}
 \left( \int_\calh  \left(\left\vert
 s^{-{1/2}}  f_n\left(s, (e^{sA}z)_0+(e^{tA}x)_0\right)
 -s^{-{1/2}}  f\left(s, (e^{sA}z)_0+(e^{tA}x)_0\right)
 \right\vert\right.\right.\\
  &\left.\left.+ \left\vert \bar \ell_{0,n}\left(s, (e^{sA}z)_0+(e^{tA}x)_0\right)
 -\bar \ell_{0}\left(s, (e^{sA}z)_0+(e^{tA}x)_0\right)\right\vert\right)^2\caln(0,Q_{t-s})(dz)
 \right)^{1/2}.
\end{align*}
Applying the Dominated Convergence Theorem as for the proof of
\myref{eq:convvnnew} we get the final claim
\[
\sup_{x\in(0,T]\times\calk} \vert t^{1/2}\nabla^B v_n(t,x)-t^{1/2}\nabla^B v(t,x)\vert\rightarrow 0, \qquad  \hbox{for any compact set $\calk\subset\calh$.}
\]
\qed
\medskip

Notice that, using the terminology of \cite{CeGo,G1}, the above result
imply that a mild solution (\ref{eq:Kforward}) is also a {\em $\calk$-strong solution}.
In general, in an infinite dimensional Hilbert space $H$, existence of $\calk$-strong solutions is not a routine application of the theory of evolution equations, as the operator $\call$ formally introduced in \ref{eq:ell}
is not the infinitesimal generator of a strongly continuous semigroup in the Banach space $UC_b(H)$. To overcome this difficulty in the already mentioned paper
\cite{CeGo} the theory of weakly continuous (or $\calk$-continuous)
semigroups has been used.

\begin{remark}\label{rm:defpistrictstrong}
The approximation results proved just above is needed to prove the fundamental identity,
(see next Proposition \ref{prop rel fond}) which is the key point to get the verification theorem and the existence of optimal feedback controls.
The idea is to apply Ito's formula to the approximating sequence $w_n$ composed with
the state process $Y$ and then to pass to the limit for $n\to +\infty$
(see e.g. \cite{G1} or \cite[Section 4.4]{FabbriGozziSwiech}).
However in the literature the approximating sequence is taken more regular,
i.e. the $w_n$ are required to be classical solutions
(see e.g. \cite[Section 6.2, p.103]{DP3}) of \myref{eq:Kforward}.
This in particular means that $\nabla w_n\in D(A^*)$ and this fact is crucial
since it allows to apply Ito's formula without requiring that the
state process $Y$ belongs to $D(A)$, which would be a too strong requirement.

In our case the used approximating procedure does not give rise in general to
functions $w_n$ with $\nabla w_n\in D(A^*)$.
Indeed for our purposes we need that the approximants of the data $\phi,\,\ell_0, \,\calf$ remain all in the space $\Sigma_{T,{1/2}}$; without this, since we only have
``partial'' smoothing,
it is not clear at all how to prove the convergence of the derivative $\nabla^B w_n$
(which is needed when we pass to the limit to prove the fundamental identity in next subsection).
Hence, in particular, since we need that, for all $n\in\N$, $\calf_n\in\Sigma^1_{T,{1/2}}$, and since
$\calf$ is written in terms of $f$, we approximate $f$
by $f_n$, and this procedure gives the approximants $\calf_n$ of $\calf$.
In this way $\calf_n\in \Sigma^1_{T,{1/2}}$ but $\nabla \calf_n \notin D(A^*)$.

Summing up, we are only able to find approximating strict solutions
and not classical solutions. Since the state process $Y$ does not belong to $D(A)$
this fact will force us to introduce suitable regularizations $Y_k$ of it
(see the proof of Proposition \ref{prop rel fond}).
\end{remark}

\begin{remark}\label{rm:defpistrictstrongbis}
Calling $\ell_n:=\ell_{0,n}+\calf_n -H_{min}(\nabla^Bw_n)$
is not difficult to see that the sequence
$w_n$ is a strict solution of the approximating HJB equation
\begin{equation*}
  \left\{\begin{array}{l}\dis
\frac{\partial w(t,x)}{\partial t}=
\frac{1}{2}Tr \;GG^*\nabla^2w(t,x)
+ \< Ax,\nabla w(t,x)\>_\calh
+H_{min} (\nabla^B w(t,x)) +\ell_n(t,x),\qquad t\in [0,T],\,
x\in \calh,\\
\\
\dis w(0,x)=\phi_n(x),
\end{array}\right.
\end{equation*}
This means, with the terminology used e.g. in \cite{G1}, that $w$ is a $\calk$-strong solution
of \myref{HJBformale}. We do not go deeper into this since here we use the approximation
only as a tool to solve our stochastic optimal control problem.
\end{remark}
\begin{remark}\label{rm:piconv}
In view of Remark \ref{remark:costoconcreto}-(i)
Lemma \ref{lm:approximation} can be easily generalized to the case when the data $\phi$
and $\ell_0$ are not bounded but satisfy a polynomial growth condition in the variable $x$
as from \myref{eq:polgrowthphil0}.

Concerning the generalization of Remark \ref{remark:costoconcreto}-(ii)
still it possible to extend the results of Lemma \ref{lm:approximation} to the case when $\phi$
and $\ell_0$ are only measurable. In this case
the approximations would take place in the sense of the $\pi$-convergence, which is weaker than the $\calk$-convergence and towards the $\calk$-convergence has also the disadvantage of being not metrizable. For more on the notion of $\pi$-convergence the reader can see
\cite[Section 6.3]{DP3} (see also \cite{EthierKurtz} and \cite{PriolaStudia}).
\end{remark}

\subsection{The Fundamental Identity and the Verification Theorem}

Now we finally go back to the control problem of Section \ref{section-statement}.
We rewrite here for the reader convenience the state equation (\ref{eq-contr-rit}),
\begin{equation*}
\left\{
\begin{array}
[c]{l}
dy(s)  =a_0 y(s) ds+b_0 u(s) ds +\int_{-d}^0b_1(\xi)u(s+\xi)d\xi+\sigma dW_s
,\text{ \ \ \ }s\in[t,T] \\
y(t)  =y_0,\\
u(\xi)=u_0(\xi), \quad \xi \in [-d,0),
\end{array}
\right.
\end{equation*}
and its abstract reformulation (\ref{eq-astr}),
\begin{equation*}
\left\{
\begin{array}
[c]{l}
dY(s)  =AY(s) ds+Bu(s) ds+GdW_s
,\text{ \ \ \ }s\in[t,T] \\
Y(t)=x=(x_0,x_1).
\end{array}
\right.
\end{equation*}
Similarly the cost functional in (\ref{costoconcreto}) is
\begin{equation*}
J(t,x,u)=\E \int_t^T \left(\bar\ell_0(s,x(s))+\ell_1(u(s))\right)\;ds +\E  \bar\phi(x(T))
\end{equation*}
and is rewritten as (see (\ref{costoconcreto1}))
\begin{equation*}
J(t,x;u)=\E \int_t^T \left(\ell_0(s,Y(s))+\ell_1(u(s))\right)\;ds +\E  \phi(Y(T)).
\end{equation*}
We notice that throughout this subsection and the following one, in order to avoid further technical difficulties,  we keep the probability space $(\Omega,\calf,\P)$
fixed. Nothing would change if we work in the weak formulation, where the probability space can change (see e.g. \cite{YongZhou99}[Chapter 2] and
\cite{FabbriGozziSwiech}[Chapter 2] for more on strong and weak formulations
in finite and infinite dimension, respectively).
We first prove the fundamental identity.
\begin{proposition}\label{prop rel fond}
Let Hypotheses \ref{ipotesibasic} and \ref{ipotesicostoconcreto} hold.
Let also (\ref{eq:hpdebreg}) or (\ref{eq:hpdebregbis}) hold.
Let $v$ be the mild solution of the HJB equation (\ref{HJBformale})
according to Definition \ref{defsolmildHJB}.
Then for every $t\in[ 0,T] $ and $x\in
H$, and for every admissible control $u$, we
have the fundamental identity
\begin{equation}\label{relfond}
 v(t,x)
=J(t,x;u)+\E\int_t^T \left[H_{min}(\nabla^B v(s,Y(s)))
- H_{CV}(\nabla^B v(s,Y(s));u(s))\right]\,ds.
\end{equation}
\end{proposition}
\dim
Take any admissible state-control couple $(Y(\cdot),u(\cdot))$,
and let $v_n(t,x):=w_n(T-t,x)$ where $(w_n)_n$ is the approximating sequence
of strict solutions defined
in Lemma \ref{lm:approximation}. We want to apply the Ito formula to
$v_n(t, Y(t))$. Unfortunately, as mentioned in Remark \ref{rm:defpistrictstrong},
this is not possible since the process
$Y(t)$ does not live in $D(A)$. So we approximate it as follows.
Set, for $k \in \N$, sufficiently large,
\begin{equation}\label{eq:Ykdef}
Y_k(s;t,x)=k(k-A)^{-1}Y(s;t,x).
\end{equation}
The process $Y_k$ is in $D(A)$, it converges to $Y$ ($\P$-a.s. and $s\in [t,T]$ a.e.)
and it is a strong solution\footnote{Here we mean strong in the probabilistic sense
and also in the sense of \cite{DP1}, Section 5.6.}
of the Cauchy problem
\begin{equation*}
\left\{
\begin{array}
[c]{l}
dY_k(s)  =AY_k(s) ds+B_ku(s) ds+G_kdW_s
,\text{ \ \ \ }s\in [t,T] \\
Y_k(t)=x_k,
\end{array}
\right.
\end{equation*}
where $B_k=k(k-A)^{-1}B$, $G_k=k(k-A)^{-1}G$ and $x_k=k(k-A)^{-1}x$.
Now observe that, thanks to \myref{resolvent}, \myref{Bnotbdd},
the operator $B_k$ is continuous, hence we can apply Dynkin's formula (see e.g.
\cite[Section 1.7]{FabbriGozziSwiech} or \cite[Section 4.5]{DP1})
to $v_{n}(s, Y_k(s))$ in the interval $[t,T]$, getting
\begin{equation}\label{Dynkinv^nk}
\E v_n(Y_k(T)) - v_n(t,x_k)
=
\E\int_t^T
\left[v_{nt}(s,Y_k(s))+  \frac{1}{2}Tr \;GG^*\nabla^2v(s,Y_k(s))
+ \< AY_k(s)+B_ku(s),\nabla v_n(s,Y_k(s))\>_\calh\right]ds.
\end{equation}
Using the Kolmogorov equation (\ref{eq:Kforward}), whose strict solution is $w_n$,
we then write
\begin{equation}\label{quasirelfondv^nk}
\E\phi_n(Y_k(T)) - v_n(t,x_k)=\E\int_t^T \left[\calf_n (s,Y_k(s))+ \ell_{0,n}(s,Y_k(s))+\<B_k u(s),\nabla v_n(s,Y_k(s))\>_{\R^m} \right]ds
\end{equation}
We first let $k\rightarrow\infty$ in \myref{quasirelfondv^nk}.
Since $\ell_{0,n}$ and $\nabla v_n$ are
bounded functions and since $\calf_n(s,x)$ has a singularity of type $s^{-{1/2}}$ with
respect to time and is bounded with respect to $x$, we can apply
the Dominated Convergence Theorem to all terms but the last getting
\begin{equation}\label{quasirelfondv^n}
\E\phi_n(Y(T)) - v_n(t,x)=\E\int_t^T [\calf_n (s,Y(s))+\ell_{0,n}(s,Y(s))]
+\lim_{k \to + \infty}\E\int_t^T\<B_k u(s),\nabla v_n(s,Y_k(s))\>_\calh]ds.
\end{equation}
Concerning the last term we observe first that, by \myref{resolvent} and
\myref{Bnotbdd} we have $B_k \to B$ in the sense that, for all
$z \in \R^n \times C([-d,0];\R^n)$, (writing down here $C$ for
$C([-d,0];\R^n)$ and $C^*$ for $C^*([-d,0];\R^n)$)
$$
\<B_k u, z\>_\calh =\<B_k u, z\>_{\<\R^n \times C^*,\R^n \times C\>)}
 \rightarrow \<B u, z\>_{
\<\R^n \times C^*,\R^n \times C\>)}
\quad as \; k \to + \infty.
$$
This in particular imply, by the Banach-Steinhaus theorem, that
$\{B_k u\}_k$ is uniformly bounded in $\R^n \times C^*([-d,0];\R^n)$.
Now we use the fact that $\nabla v_n(s,x) \in \R^n \times C([-d,0];\R^n)$
(see \myref{eq:nablaperSigma})
to rewrite the integrand of the last term of \myref{quasirelfondv^n} as
$$
\<B_k u(s),\nabla v_n(s,Y_k(s))-\nabla v_n(s,Y(s))\>_{
\<\R^n \times C^*,\R^n \times C\>}
+\<B_k u(s),\nabla v_n(s,Y(s))\>_{
\<\R^n \times C^*,\R^n \times C\>}
$$
Thanks to what said above the first term goes to $0$ as $k \to + \infty$ and is dominated while the second term is also dominated and converges to
$\<B u(s),\nabla v_n(s,Y(s))\>_{\<\R^n \times C^*,\R^n \times C\>}$
which, thanks to \myref{eq:nablaperSigma}, is equal to
$\< u(s),\nabla^B v_n(s,Y(s))\>_{\R^m}$ (both convergences are clearly
$\P$-a.s. and $s\in [t,T]$ a.e.). Hence
$$
\lim_{k \to + \infty}\E\int_t^T\<B_k u(s),\nabla v_n(s,Y_k(s))\>_\calh]ds
=
\E\int_t^T\<u(s),\nabla^B v_n(s,Y(s))\>_{\R^m}]ds
$$
Now we let $n\rightarrow\infty$. By Lemma \ref{lm:approximation}, we know that
$v_n(t,x)\rightarrow v(t,x)$
and $(T-t)^{1/2} \nabla^Bv_n(t,x)\rightarrow (T-t)^{1/2}\nabla^B v(t,x)$
pointwise. Moreover $v_n(t,x)$, $(T-t)^{1/2}\nabla^Bv^n(t,x)$
are uniformly bounded, so that, by dominated convergence, we get
\[
\E\int_t^T \<u(s),\nabla^Bv_n(s,Y(s))\>_{\R^m}]\,ds\rightarrow
\E\int_t^T \<u(s),\nabla^Bv(s,Y(s))\>_{\R^m}]\,ds.
\]
The convergence
\[
 \E\phi_n(x(T)) -\E\int_t^T  [\calf_n(s,Y(s))+ \ell_{0,n}(s,Y(s))]ds
 \rightarrow
\E\phi(x(T)) -\E\int_t^T [H_{min}(\nabla^B v(s,Y(s)))+ \ell_0(s,Y(s))]ds
\]
follows directly by the construction of the approximating sequences $(\phi_n)_n$
$(\calf_n)_n$ and $(\ell_{0,n})_n$.
Then, adding and subtracting $\E\int_t^T  \ell_1(u(s))ds$
and letting $n\rightarrow\infty$ in (\ref{quasirelfondv^n})
we obtain
\begin{align*}
 v(t,x)&
=\E\phi(Y(T))+\E\int_t^T [\ell_0(s,Y(s))+ \ell_1(u(s))]ds\\
&+\E\int_t^T \left[H_{min}(\nabla^B v(s,Y(s)))-
H_{CV}(\nabla^B v(s,Y(s));u(s))\right]
\,ds
\end{align*}
which immediately gives the claim.
\qed


\begin{remark}
One may wonder why we approximate the process $Y$ with $Y_k$ as in \myref{eq:Ykdef}
instead of using the Yosida approximants $A_k$ of $A$ as, e.g., in
\cite[p.144]{DP3}. The reason is that we need that $Y_k$ belongs to $D(A)$,
which is not guaranteed if we use Yosida approximants.
A similar procedure is used, in a different context, in
the book \cite{DP1}, in the proof of Theorem 7.7, p. 203.
\end{remark}

We can now pass to prove our Verification Theorem i.e.
a sufficient condition of optimality given in term of the
mild solution $v$ of the HJB equation (\ref{HJBformale}).

\begin{theorem}
\label{teorema controllo}
Let Hypotheses \ref{ipotesibasic} and \ref{ipotesicostoconcreto} hold.
Let also (\ref{eq:hpdebreg}) or (\ref{eq:hpdebregbis}) hold.
Let $v$ be the mild solution of the HJB equation (\ref{HJBformale})
whose existence and uniqueness is proved in (\ref{esistenzaHJB}).
Then the following holds.
\begin{itemize}
 \item For all $(t,x)\in [0,T]\times \calh$ we have
$v(t,x) \le V(t,x)$, where $V$ is the value function
defined in (\ref{valuefunction}).
\item
Let $t\in [0,T]$ and $x\in \calh$ be fixed.
If, for an admissible control $u^*$, we
have, calling $Y^*$ the corresponding state,
$$
u^*(s)\in \arg\min_{u\in U}H_{CV}(\nabla^B v(s,Y^*(s);u)
$$
for a.e. $s\in [t,T]$, $\P$-a.s., then the pair $(u^*,Y^*)$ is
optimal for the control problem starting from $x$ at time $t$
and $v(t,x)=V(t,x)=J(t,x;u^*)$.
\end{itemize}
\end{theorem}
\dim The first statement follows directly by \eqref{relfond} due to the
negativity of the integrand.
Concerning the second statement, we immediately see that, when $u=u^*$
(\ref{relfond}) becomes
$v(t,x)=J(t,x;u^*)$.
Since we know that for any admissible control $u$
\[
 J(t,x;u)\geq V(t,x) \ge v(t,x),
\]
the claim immediately follows.
\qed

\subsection{Optimal feedback controls and $v=V$}
\label{sec:contr-feedback}
We now prove the existence of optimal feedback controls. Under the Hypotheses
of Theorem \ref{teorema controllo} we define, for $(s,y)\in [0,T)\times \calh$,
the {\em feedback map}
\begin{equation}\label{defdiPsi}
\Psi(s,y):=\arg \min_{u\in U} H_{CV}(\nabla^Bv(s,y);u),
\end{equation}
where, as usual, $v$ is the solution of the HJB equation \myref{HJBformale}.
Given any $(t,x)\in [0,T)\times \calh$,
the so-called Closed Loop Equation (which here is, in general, an inclusion)
is written, formally, as
\begin{equation}\label{cleinclusion}
\left\{
\begin{array}
[c]{l}
dY(s) \in AY(s) ds+B\Psi\left(s,Y(s)\right) ds+GdW_s
,\text{ \ \ \ }s\in [t,T) \\
Y(t)=x.
\end{array}
\right.
\end{equation}
First of all we have the following straightforward corollary whose proof is immediate from Theorem \ref{teorema controllo}.

\begin{corollary}
\label{cr:optimalfeedback}
Let the assumptions of Theorem \ref{teorema controllo} hold true.
Let $v$ be the mild solution of \myref{HJBformale}.
Fix $(t,x)\in [0,T)\times H$ and assume that, on $[t,T)\times H$, the map $\Psi$
defined in (\ref{defdiPsi}) admits a measurable selection
$\psi:[t,T)\times H \to \Lambda$ such that the Closed Loop Equation
\begin{equation}
\label{eq:CLEselection}
\left \{
\begin{array}{l}
d Y(s) = AY(s)d s+B\psi\left(s,Y(s)\right) ds+GdW_s
,\text{ \ \ \ }s\in [t,T) \\
Y(t)=x.
\end{array}
\right.
\end{equation}
has a mild solution $Y_\psi(\cdot;t,x)$  (in the sense of \cite[p.187]{DP1}).
Define, for $s \in [t,T)$, $u_\psi (s)=\psi(s,Y_\psi(s;t,x))$.
Then the couple
$(u_\psi(\cdot),Y_\psi(\cdot;t,x))$ is optimal at
$(t,x)$ and $v(t,x)=V(t,x)$.
If, finally, $\Psi(t,x)$ is always a singleton and the mild solution
of \myref{eq:CLEselection} is unique,
then the optimal control is unique.
\end{corollary}

We now give sufficient conditions to verify the assumptions of Corollary
\ref{cr:optimalfeedback}.
First of all define
 \begin{equation}\label{defdigammagrandebis}
\Gamma(p):=\left\{ u\in U: \<p,u\>+\ell_1(u)= H_{min }(p)\right\}.
\end{equation}
Then, clearly, we have $\Psi(t,x)=\Gamma(\nabla^B v(t,x))$.
Observe that, under mild additional conditions on $U$ and $\ell_1$
(for example taking $U$ compact or $\ell_1$ of superlinear growth),
the set $\Gamma$ is nonempty for all $p \in \R^m$.
If this is the case then, by \cite{AubFr}, Theorems 8.2.10 and
8.2.11, $\Gamma$ admits a measurable selection, i.e. there exists
a measurable function $\gamma:\R^m \rightarrow U$ with
$\gamma(z)\in \Gamma(z)$ for every $z\in \R^m$.
Since $H_{min }$ is Lipschitz continuous, then $\Gamma$, and so
$\gamma$, must be uniformly bounded. In some cases
studied in the literature this is enough to find an optimal feedback
but not in our case (read on this the subsequent Remark \ref{rm:discfeed}-(ii)).
Hence to prove existence of a mild solution of the closed loop
equation \myref{eq:CLEselection}, as requested in Corollary
\ref{cr:optimalfeedback}, we need more regularity of the feedback
term $\psi(s,y)=\gamma(\nabla^B v(s,y))$.
Beyond the smooth assumptions on the coefficients required in Theorem \ref{lemma-stimev}, which give the regularity of
$\nabla^B v(t,x)$, we need the following assumption about the map $\Gamma$.

\begin{hypothesis}\label{hp:lipsel}
The set-valued map $\Gamma$ defined in (\ref{defdigammagrandebis}) is
always non empty; moreover it admits a Lipschitz continuous selection $\gamma$.
\end{hypothesis}

\begin{remark}
\label{rm:discfeed}
\begin{itemize}
  \item[(i)]
 Notice that the above Hypothesis is verified if we assume that $\ell_1:\R^m\rightarrow \R$ is differentiable with an invertible derivative and that $(\ell_1')^{-1}$ is Lipschitz continuous.
Indeed in this case, see (\ref{psi1}),
the infimum of $H_{CV}$ is achieved at
\[
u=(\ell_1')^{-1}(z), \text{ so that }\Gamma_0(z)=(\ell_1')^{-1}(z).
\]
\item[(ii)]
The problem of the lack of regularity of the feedback law
is sometimes faced (see e.g. in \cite{FT2}) by formulating the optimal control problem in the weak sense (see e.g. \cite{FlSo} or \cite{YongZhou99}, Section 4.2)
and then using Girsanov Theorem to prove existence, in the weak sense, of a mild solution of (\ref{eq:CLEselection}) when the map $\psi$ is only measurable and bounded.
This is not possible here due to the already mentioned absence of the structure condition in the controlled state equation (i.e. the control affects the system not only through the noise).
\end{itemize}
\end{remark}

Taking the selection $\gamma$ from Hypothesis \ref{hp:lipsel}
we consider the closed loop equation
\begin{equation}\label{cle}
\left\{
\begin{array}
[c]{l}
dY(s)  =AY(s) ds+B\gamma\left(\nabla^B v(s,Y(s))\right) ds+GdW_s
,\text{ \ \ \ }s\in[t,T] \\
Y(t)=x=(x_0,x_1),
\end{array}
\right.
\end{equation}
and we have the following result.
\begin{theorem}\label{teo su controllo feedback}
Assume that Hypotheses \ref{ipotesibasic},
\ref{ipotesicostoconcreto}, \ref{ipotesicostoconcretobis} and \ref{hp:lipsel} hold true.
Fix any $(t,x)\in [0,T)\times H$.
Let also (\ref{eq:hpdebreg}) or (\ref{eq:hpdebregbis}) hold.
Then the closed loop equation (\ref{cle}) admits a unique mild solution
$Y_\gamma(\cdot;t,x)$ (in the sense of \cite[p.187]{DP1})
and setting
$$
u_\gamma(s)=\gamma\left(\nabla^{B} v(s,Y_\gamma(s;t,x) )\right), \quad s \in [t,T]
$$
we obtain an optimal control at $(t,x)$. Moreover $v(t,x)=V(t,x)$.
\end{theorem}

\noindent {\bf Proof.} Thanks to Corollary \ref{cr:optimalfeedback}
it is enough to prove the existence and uniqueness of the mild solution of (\ref{cle}). We apply a fixed point theorem to the following integral form of
(\ref{cle}):
\begin{equation}\label{clemild}
Y(s)= e^{(s-t)A}x+
\int_t^s e^{(s-r)A}G\,dW_r
+\int_{t}^s e^{(s-r)A}B\gamma(\nabla^{B}v(r ,\overline{X}_r)))dr.
\end{equation}
By Hypothesis \ref{ipotesicostoconcretobis} and Theorem \ref{lemma-stimev}
we get that the mild solution $v$ of the HJB equation (\ref{HJBformale})
is differentiable, with bounded derivative. Moreover, since, again by
Theorem \ref{lemma-stimev},
$v$ admits the second order derivative $\nabla^B\nabla v$, and
$\nabla^B\nabla v=\nabla\nabla^B v$,
we deduce that $t^{1/2}\nabla^B v(t,\cdot)$ is
Lipschitz continuous, uniformly with respect to $t$.
Using this Lipschitz property we can solve (\ref{clemild}) by a fixed point argument. Since the argument to do this is straightforward
we only show how to estimate the difficult term in \myref{clemild}. We have
\begin{align*}
\int_{t}^s& \vert e^{(s-r)A}B\left(\gamma(\nabla^{B}
v(r ,\overline{X}(r)))-\gamma(\nabla^{B}
v(r ,\overline{Y}(r)))\right)\vert_\calh dr
\leq \int_{t}^s C\vert\gamma(\nabla^{B}
v(r ,\overline{X}(r)))-\gamma(\nabla^{B}
v(r ,\overline{Y}(r)))\vert_{\R^m}] dr
\\
&
\leq C\int_{t}^s \vert\nabla^{B}
v(r ,\overline{X}(r)))-\nabla^{B}
v(r ,\overline{Y}(r))\vert_{\R^m}] dr
\leq C\int_{t}^s r^{-{1/2}}\vert \overline{X}(r)-\overline{Y}(r)\vert_\calh dr
\end{align*}
where $C$ is a constant that can change its values from line to line.
\qed

We devote our final result to show that the identification $v=V$ can be done,
using an approximation procedure, also
in cases when we do not know if optimal feedback controls exist.


\begin{theorem}\label{teo:v=V}
Let Hypotheses \ref{ipotesibasic}, \ref{ipotesicostoconcreto} hold.
Let also (\ref{eq:hpdebreg}) or (\ref{eq:hpdebregbis}) hold.
Moreover let Hypotheses \ref{ipotesicostoconcretobis}-(ii)
and \ref{hp:lipsel} hold and let
$\phi$ and $\ell_0$ be uniformly continuous.
Then $v=V$.
%
\end{theorem}
\dim
We approximate $\ell_0$ and $\phi$ by approximating $\bar\ell_0$ and $\bar\phi$
with standard approximants $\bar\ell_{0,n}$ and $\bar\phi_n$ built
by convolutions.
We set
\begin{equation}\label{costo_n}
J_n(t,x;u)=\E \int_t^T \left(\ell_{0,n}(s,Y(s))+\ell_1(u(s))\right)\;ds +\E  \phi_n(Y(T))
\end{equation}
and call $w_n$ the mild solution of the HJB equation
($\call$ is given by \myref{eq:ell})
\begin{equation}\label{HJBformale-n}
  \left\{\begin{array}{l}\dis
\frac{\partial w(t,x)}{\partial t}=\call [w(t,\cdot)](x) +\ell_{0,n}(t,x)+
H_{min} (\nabla^B w(t,x)),\qquad t\in [0,T],\,
x\in \calh,\\
\\
\dis w(0,x)=\phi_n(x).
\end{array}\right.
\end{equation}
In particular $w_n$ satisfies the integral equation
\begin{equation}
  w_n(t,x) =R_{t}[\phi_n](x)+\int_0^t R_{t-s}[
H_{min}(
\nabla^B w_n(s,\cdot)+\ell_{0,n}(s,\cdot)
](x)\; ds,\qquad t\in [0,T].\
x\in \calh,\label{solmildHJB-forwardn}
\end{equation}
By Theorem \ref{teo su controllo feedback}
calling $v_n(t,x)=w_n(T-t,x)$ we have
\begin{equation}\label{identif-v_n}
  v_n(t,x) =V_n(t,x):=\inf_{u\in \calu} J_n(t,x;u).
\end{equation}
and there exists an optimal feedback control $u_n(s)=\psi_n(s,Y(s))$.
Moreover, by Lemma \ref{lm:approximation} we know that
$$
v_n(t,x) \overset{\calk}{\rightarrow} v(t,x).
$$
Now it is enough to prove that $V_n(t,x) \to V(t,x)$ pointwise.
Given $\eps >0$, we have, for $n$ large enough,
\begin{align*}
 V_n(t,x)&=\inf_{u\in\calu}
\left[\E \int_t^T \left(\ell_{0,n}(s,Y(s))+\ell_1(u(s))\right)\;ds +\E  \phi_n(Y(T))\right]\\
&=\inf_{u\in\calu}\left[\E \int_t^T \left(\ell_{0}(s,Y(s))+\ell_1(u(s))\right)\;ds +\E \phi(Y(T))\right.\\
&\left.+\E \int_t^T \left(\ell_{0,n}(s,Y(s))-\ell_{0}(s,Y(s))\right)\;ds
+\E\left[\phi_n(Y(T)) -\phi(Y(T))\right]\right]\\
&\leq\inf_{u\in\calu}\left[\E \int_t^T \left(\ell_{0}(s,Y(s))+\ell_1(u(s))\right)\;ds +\E \phi(Y(T))\right]
+\varepsilon,
\end{align*}
where the last passage follows by the Dominated Convergence Theorem, and since $\phi$
and $\ell_{0,n}$ are uniformly continuous.
We have shown that
\[
 V_n(t,x)\leq V(t,x)+\varepsilon.
\]
Exchanging the role of $V_n$ and $V$ we also find that the reverse inequality holds true. Hence $V_n \to V$ pointwise an the claim follows.
\qed


\begin{remark}\label{rm:unif-cont-feedback}
 In Theorem \ref{teo:v=V} we have assumed further uniform continuity on the data.
When $U$ is compact the result still remain true if the data are only continuous.
\end{remark}



\begin{thebibliography}{99}

 \bibitem{AubFr} {J.P. Aubin, H. Frankowska},
    {\textit{Set-valued analysis}},
      {Systems \& Control: Foundations \& Applications},
Vol. {2}, {Birkh\"auser Boston Inc.}, {Boston, MA},     {1990}

\bibitem{BardiDolcetta} M. Bardi, I. Capuzzo-Dolcetta, \textit{Optimal control and viscosity solutions of Hamilton-Jacobi-Bellman equations. With appendices by Maurizio Falcone and Pierpaolo Soravia.}
{Systems \& Control: Foundations \& Applications. Birkh\"{a}user Boston, Inc., Boston, MA, 1997.}


\bibitem{BDDM07} A. Bensoussan, G. DaPrato, M. Delfour, S. Mitter
\textit{Representation and Control of infinite dimensional systems}, 2
nd edition, Birkh?auser, 2007 xxvi + 575 pp. ISBN 978-0-8176-4461-1.

\bibitem{PhamBruder} B. Bruder and H. Pham,
\textit{Impulse control problem on finite horizon with execution delay}.
Stochastic Process. Appl. 119 (2009), no. 5, 1436--1469.

\bibitem {CDP1}P. Cannarsa and G. Da Prato, \textit{Second order
Hamilton-Jacobi equations in infinite dimensions, }SIAM J. Control Optim, 29,
2, (1991), pp. 474-492.

\bibitem {CDP2}P. Cannarsa and G. Da Prato, \textit{Direct solution of a
second order Hamilton-Jacobi equations in Hilbert spaces, S}tochastic Partial
Differential Equations and Applications, (1992), pp. 72-85.


\bibitem{CannarsaSoner} P. Cannarsa, H. M. Soner, \textit{Generalized one-sided estimates for solutions of Hamilton-Jacobi equations and applications.} Nonlinear Anal. 13 (1989), no. 3, 305--323.

\bibitem{Ce95} S. Cerrai, \textit{Weakly continuous semigroups in the space of functions with polynomial growth.} Dynam. Systems Appl. 4 (1995), no. 3, 351--371.

\bibitem {CeGo} S. Cerrai and F. Gozzi,
 \textit{Strong solutions of Cauchy problems
associated to weakly continuous semigroups},
   Differential Integral Equations (1995), no. 3, 465--486.


\bibitem{ChowMenaldi}
P.L. Chow, J.L. Menaldi,
\textit{Infinite-dimensional Hamilton-Jacobi-Bellman equations in Gauss-Sobolev spaces.} Nonlinear Anal. 29 (1997), no. 4, 415--426.

\bibitem {DP1}G. Da Prato and J. Zabczyk,\textit{\ Stochastic equations in
infinite dimensions. Second Edition.}
Encyclopedia of Mathematics and its Applications 152,
Cambridge University Press, \textit{2014.}

\bibitem {DaPratoZabczyk95}G. Da Prato and J. Zabczyk,\textit{\ Regular densities of invariant measures in
Hilbert spaces. } J. Funct. Anal. 130 (1995), no. 2, 427-449.

\bibitem {DP3}G. Da Prato and J. Zabczyk,\textit{\ Second order partial
differential equations in Hilbert spaces}. London Mathematical Society Note
Series, 293, Cambridge University Press, Cambridge, 2002.

\bibitem {EthierKurtz} S. N. Ethier and T. G. Kurtz,
\textit{Markov processes. Characterization and convergence.} Wiley Series in Probability and Mathematical Statistics: Probability and Mathematical Statistics. John Wiley \& Sons, Inc., New York, 1986.

\bibitem{FabbriFederico} G. Fabbri, S. Federico
\textit{On the Infinite-Dimensional Representation of Stochastic Controlled Systems with Delayed Control in the Diffusion Term.} Mathematical Economic Letters 2 (2014), no. 3-4, 33--44

\bibitem{FabbriGozziSwiech} G. Fabbri, F. Gozzi, A. Swiech \textit{Stochastic Optimal Control in Infinite Dimensions: Dynamic Programming and HJB Equations}. Forthcoming.

\bibitem{FabbriGozzi08} G. Fabbri and F. Gozzi, \textit{Solving Optimal
Growth Models with Vintage Capital: the Dynamic Programming Approach.}
 J. Econom. Theory 143 (2008), no. 1, 331-373.


\bibitem{FedFinSto} S. Federico, \textit{A stochastic control problem with delay
arising in a pension fund model. }Finance Stoch. 15 (2011), no. 3, 421--459

\bibitem{FedGolGozSICON1} S. Federico, B. Goldys, F. Gozzi
\textit{HJB equations for the optimal control of differential equations with delays and state constraints, I: regularity of viscosity solutions.} SIAM J. Control Optim. 48 (2010), no. 8, 4910--4937.

\bibitem{FedGolGozSICON2} S. Federico, B. Goldys, F. Gozzi
 \textit{HJB equations for the optimal control of differential equations with delays and state constraints, II: Verification and optimal feedbacks.} SIAM J. Control Optim. 49 (2011), no. 6, 2378--2414.


\bibitem{FedTacSicon} S. Federico and E. Tacconi, \textit{Dynamic Programming for
optimal control problems with delays in the control variable}.
SIAM J. Control Optim. 52 (2014), no. 2, 1203-1236.


\bibitem{FlaZan} F. Flandoli and G. Zanco, \textit{An infinite-dimensional approach to path-dependent Kolmogorov's equations}. Preprint, arXiv:1312.6165. To appear in Annals of Probability.

\bibitem{FlSo} W. H. Fleming, H. M. Soner,
\textit{Controlled Markov processes and viscosity solutions.} Applications
of Mathematics 25. Springer-Verlag, 1993.


\bibitem {FT2}M. Fuhrman and G. Tessitore, \textit{Non linear Kolmogorov
equations in infinite dimensional spaces: the backward stochastic differential
equations approach and applications to optimal control. }Ann. Probab. 30
(2002), no. 3, 1397--1465.

\bibitem {FTGgrad}M. Fuhrman and G. Tessitore, \textit{Generalized directional gradients, backward stochastic differential equations and mild solutions of semilinear parabolic equations}. Appl. Math. Optim. 51 (2005), no. 3, 279--332.



\bibitem {GGSPA}B. Goldys and F.Gozzi, \textit{Second order parabolic HJ
equations in Hilbert spaces and stochastic L}$^{2}$ \textit{approach,
} Stochastic processes and Applications, 116 (2006) 1932--1963.


\bibitem {G1}F. Gozzi, \textit{Regularity of solutions of second order
Hamilton-Jacobi equations in Hilbert spaces and applications to a control
problem, }(1995) Comm Partial \ Differential Equations 20, pp. 775-826.

\bibitem {G2}F. Gozzi, \textit{Global regular solutions of second order
Hamilton-Jacobi equations in Hilbert spaces with locally Lipschitz
nonlinearities}, (1996) J. Math. Anal. Appl. 198, pp. 399-443.


\bibitem {GM}F. Gozzi and C. Marinelli, \textit{Stochastic optimal control of delay equations arising in advertising models}.
    Stochastic partial differential equations and applications - VII,  133-148, Lect. Notes Pure Appl. Math., 245, Chapman $\&$ Hall/CRC, Boca Raton, FL, 2006.


\bibitem{GMSJOTA} F. Gozzi, C. Marinelli, S. Savin \textit{On controlled linear diffusions with delay in a model of optimal advertising under uncertainty with memory effects}.
J. Optim. Theory Appl. 142 (2009), no. 2, 29--321.

\bibitem{Henry} D. Henry, \textit{Geometric Theory of Semilinear Parabolic Equations}.
Lecture Notes in Mathematics, 840, Springer-Verlag, Berlin-New York, 1981.


\bibitem{IkedaWatanabe} N. Ikeda, S. Watanabe, \textit{Stochastic differential equations and diffusion processes. Second edition.} North-Holland Mathematical Library, 24. North-Holland Publishing Co., Amsterdam; Kodansha, Ltd., Tokyo, 1989.

\bibitem{Kol-Sha} V.B. Kolmanovskii, L.E. Shaikhet
\textit{Control of Systems with Aftereffect,} Springer.
Translation of Mathematical Monographs, Vol. 157, AMS, Providence (1996).


\bibitem{LarssenRisebro03}
B.~Larssen and N.~H. Risebro.
\textit{ When are {HJB}-equations in stochastic control of delay systems
  finite dimensional?}
Stochastic Anal. Appl., 21 (3): 643--671, 2003.


\bibitem {Lu}A. Lunardi, \textit{Schauder estimates for a class of degenerate
elliptic and parabolic operators with unbounded coefficients in }$R^{n}%
$\textit{.} Ann. Scuola Norm. Sup. Pisa Cl. Sci. (4) 24 (1997), no. 1, 133--164


\bibitem {Mas}F. Masiero, \textit{Semilinear Kolmogorov equations and applications to stochastic optimal control}.  Appl. Math. Optim.  51  (2005),  no. 2, 201-250.

\bibitem{Mas-inf-or} F. Masiero, \textit{Infinite horizon stochastic optimal control problems with degenerate noise and elliptic equations in Hilbert spaces.} Appl. Math. Optim. 55 (2007), no. 3, 285-326.


\bibitem{MasBanach} F. Masiero, \textit{Stochastic optimal control problems and parabolic equations in Banach spaces.} SIAM J. Control Optim. 47 (2008), no. 1, 251-300.

\bibitem{PaPe} E. Pardoux, S. Peng,
\textit{Backward stochastic differential equations and quasilinear
parabolic partial differential equations,} in: Stochastic
partial differential equations and their applications, eds. B.L.
Rozowskii, R.B. Sowers, 200-217, Lecture Notes in Control Inf.
Sci. 176, Springer, 1992.

\bibitem{PriolaTesi} E. Priola, \textit{Partial differential equations with infinitely
many variables,} PhD. Thesis, (1999).

\bibitem {PriolaStudia} E. Priola, \textit{On a class of Markov type semigroups in
spaces of uniformly continuous and bounded functions,} Studia Math. 136 (3) (1999), pp. 271-295.

\bibitem {Seidman}T.I. Seidman, \textit{How violent are fast controls?,
I,} Math. of Control, Signals and Syst. 1 (1988), pp. 89-95.

\bibitem {SY}T.I. Seidman and J. Yong, \textit{How violent are fast controls?,
II,} Math. of Control, Signals and Syst. 9 (1997), pp. 327-340.

\bibitem {VK}R. B. Vinter and R. H. Kwong, \textit{\ The infinite time quadratic control problem for linear systems
with state and control delays: an evolution equation approach}. SIAM J. Control Optim., 19
(1):139-153, 1981.

\bibitem{YongZhou99} J. Yong and X.Y. Zhou
\textit{Stochastic Control and Hamilton-Jacobi-Bellman equations.} Applications of Mathematics (New York), 43. Springer-Verlag, New York, 1999.

\bibitem {Z}J. Zabczyk\textit{, Mathematical Control Theory: an Introduction,
} Systems \& Control: Foundations \& Applications. Birkh\"{a}user Boston,
Inc., Boston, MA, 1992.

\end{thebibliography}
\end{document}